\documentclass[preprint,fleqn,3p]{elsarticle}
\usepackage{amsmath}
\usepackage{amssymb,amsfonts,bm}
\usepackage{amsthm}
\usepackage{dsfont}
\usepackage{color,graphicx,caption,subcaption}
\usepackage{hyperref}
\usepackage{xspace}
\usepackage{fullpage}
\usepackage{placeins}
\usepackage{enumitem}
\usepackage{url}
\usepackage{verbatim}
\usepackage{amsthm}
\usepackage{commath}
\usepackage{multirow}
\usepackage{algorithm2e}
\usepackage{tabularx}

\newcommand{\gm}{\gamma}

\newcommand\bx{\mathbf{x}}

\newcommand{\p}{\partial}
\newcommand{\og}{\omega}

\newcommand{\fl}[2]{\frac{#1}{#2}}

\newcommand{\be}{\begin{equation}}
\newcommand{\ee}{\end{equation}}
\newcommand{\ba}{\begin{array}}
\newcommand{\ea}{\end{array}}
\def\bea{\begin{eqnarray}}
\def\eea{\end{eqnarray}}
\def \beas{\begin{eqnarray*}}
\def \eeas{\end{eqnarray*}}

\newcommand{\lela}{\left \langle}
  \newcommand{\rira}{\right \rangle}
 
\newtheorem{exmp}{Example}[section]

\renewcommand{\Re}{\operatorname{Re}}

\renewcommand{\norm}[1]{\left\| #1 \right\|}

\newcommand{\C}{\mathbb C}

\DeclareMathOperator*{\argmin}{arg\,min}

\newcommand\dsp{\displaystyle}
\newcommand\xbf{\mathbf{x}}

\def\<{\langle}
\def\>{\rangle}

\def\ba{{\bf a}}

\newcommand{\eps}{\varepsilon}
\newcommand{\R}{\mathbb R}

\begin{document}

\begin{frontmatter}

\title{Efficient spectral computation of the stationary states of rotating Bose-Einstein condensates
by the preconditioned nonlinear conjugate gradient method}
\author[1]{Xavier Antoine}
\ead{xavier.antoine@univ-lorraine.fr}

\address[1]{Institut Elie Cartan de Lorraine, Universit\'e de Lorraine, Inria Nancy-Grand Est, SPHINX Team,
  F-54506 Vandoeuvre-l\`es-Nancy Cedex, France.}

\author[2]{Antoine Levitt}
 \ead{antoine.levitt@inria.fr}
\address[2]{Inria Paris, F-75589 Paris Cedex 12, Universit\'e Paris-Est, CERMICS (ENPC), F-77455 Marne-la-Vall\'ee}

\author[1,3]{Qinglin Tang}
\ead{tqltql2010@gmail.com}

\address[3]{Laboratoire de Math\'ematiques Rapha\"el Salem, Universit\'e de Rouen,
Technop\^{o}le du Madrillet, 76801 Saint-Etienne-du-Rouvray, France.}



\begin{abstract}
  We propose a preconditioned nonlinear conjugate gradient method
  coupled with a spectral spatial discretization scheme for computing
  the ground states (GS) of rotating Bose-Einstein condensates (BEC),
  modeled by the Gross-Pitaevskii Equation (GPE).  We first start by
  reviewing the classical gradient flow (also known as {\it imaginary
    time} (IMT)) method which considers the problem from the PDE
  standpoint, leading to numerically solve a dissipative equation.
  Based on this IMT equation, we analyze the forward Euler (FE),
  Crank-Nicolson (CN) and the classical backward Euler (BE) schemes
  for linear problems and recognize classical power iterations,
  allowing us to derive convergence rates.  By considering the
  alternative point of view of minimization problems, we propose the
  preconditioned gradient (PG) and conjugate gradient (PCG) methods
  for the GS computation of the GPE.  We investigate the choice of the
  preconditioner, which plays a key role in the acceleration of the
  convergence process. The performance of the new algorithms is tested
  in 1D, 2D and 3D.  We conclude that the PCG method outperforms all the
  previous methods, most particularly for 2D and 3D fast rotating
  BECs, while being simple to implement.
  \end{abstract}


\end{frontmatter}

\tableofcontents

\section{Introduction}
Bose-Einstein Condensates (BECs) were first predicted theoretically by
S.N. Bose and A. Einstein, before being realized experimentally in
1995
\cite{ISI:A1995RJ02900031,ISI:A1995RQ99000002,DalfovoReview,ISI:A1995TF68400001}. This
state of matter has the interesting feature that macroscopic quantum
physics properties can emerge and be observed in laboratory
experiments.  The literature on BECs has grown extremely fast over the
last 20 years in atomic, molecular, optics and condensed matter
physics, and applications from this new physics are starting to appear
in quantum computation for instance \cite{ByrnesWenYamamoto}. Among
the most important directions, a particular attention has been paid
towards the understanding of the nucleation of vortices
\cite{ISI:000168187300038,ISI:000188785200003,ISI:000168623500006,
  ISI:000085005700004,ISI:000082705600005,ISI:000172400600002,ISI:000279003500019}
and the properties of dipolar gases
\cite{BaoCaiReview,BaoDipolar} or multi-components BECs
\cite{bao2004groundmulti,bao2011ground,BaoCaiReview}. At temperatures
$T$ which are much smaller than the critical temperature $T_{c}$, the
macroscopic behavior of a BEC can be well described by a condensate
wave function $\psi$ which is solution to a Gross-Pitaevskii Equation
(GPE).  Being able to compute efficiently the numerical solution of
such a class of equations is therefore extremely useful. Among the
most crucial questions are the calculations of stationary states,
i.e. ground/excited states, and of the real-time dynamics
\cite{Antoine2013,AntDubBookChapter,BaoCaiReview}.

To fully analyze a representative and nontrivial example that can be extended to
other more general cases, we consider in this paper a BEC that can be
modeled by the rotating (dimensionless) GPE.  In this setting, the computation of a ground state of a
$d$-dimensional BEC takes the form of a constrained minimization
problem: find $\phi \in L^{2}(\R^{d})$ such that
\begin{align*}
  \phi \in \argmin_{\norm{\phi} = 1} E(\phi),
\end{align*}
where $E$ is the associated energy functional.
Several approaches can be developed for computing the stationary state
solution to the rotating GPE.  For example, some techniques are based
on appropriate discretizations of the continuous normalized gradient
flow/imaginary-time formulation
\cite{Adhikari1,AntoineDuboscqKrylovJCP,
  AntDubBookChapter,BaoCaiReview,BaoDuSISC,ImaginaryBaye,Cerimele1,Chiofalo1,ZengZhang}, leading to various
iterative algorithms. These approaches are general and can be applied
to many situations (dipolar interactions, multi-components GPEs...).
We refer for instance to the recent freely distributed Matlab solver GPELab
that provides the stationary states computation \cite{ABGPELabI} (and
real-time dynamics \cite{ABGPELabII}) for a wide variety of GPEs based
on the so-called BESP (Backward Euler pseudoSPectral)
scheme
\cite{AntoineDuboscqKrylovJCP,AntDubBookChapter,BaoCaiReview,BaoDuSISC}
(see also Sections \ref{ReviewSection} and
\ref{sectionNumericalResultats}).  Other methods are related to the
numerical solution of the nonlinear eigenvalue problem
\cite{DionCances,jeng2013two} or on {optimization techniques under
constraints}
\cite{BaoMinimization,CaliariMinimization,DanailaHecht1,Danaila2010}. As
we will see below in Section \ref{ReviewSection}, some
connections exist between these approaches.  Finally, a regularized
Newton-type method was proposed recently in \cite{weizhu2015}.

Optimization problems with orthogonal or normalization constraints
also occur in different branches of computational science. An
elementary but fundamental example is the case of a quadratic energy,
where solving the minimization problem is equivalent to finding an
eigenvector associated with the lowest eigenvalue of the symmetric matrix
representing the quadratic form. A natural generalization is a class
of orthogonalized minimization problems, which for a quadratic energy
reduce to finding the $N$ first eigenvectors of a matrix. Many
problems in electronic structure theory are of this form, including
the celebrated Kohn-Sham and Hartree-Fock models
\cite{cances2003computational, saad2010numerical}. Correspondingly, a
large amount of effort has been devoted to finding efficient
discretization and minimization schemes. A workhorse of these
approaches is the nonlinear preconditioned conjugate gradient method,
developed in the 80s \cite{payne1992iterative}, as well as several
variants of this (the Davidson algorithm, or the LOBPCG method
\cite{knyazev2001toward}).

Although similar, there are significant differences between the
mathematical structure of the problem in electronic structure theory
and the Gross-Pitaevskii equation. In some respects, solving the GPE
is easier: there is only one wavefunction (or only a few for
multi-species gases), and the nonlinearity is often local (at least
when dipolar effects are not taken into account), with a simpler
mathematical form than many electronic structure models. On the other
hand, the Gross-Pitaevskii equation describes the formation of
vortices: the discretization schemes must represent these very
accurately, and the energy landscape presents shallower minima,
leading to a more difficult optimization problem.

In the present paper, we consider the constrained nonlinear conjugate gradient
method for solving the rotating GPE (Section \ref{sec:defs}) with a
pseudospectral discretization scheme (see Section \ref{SectionPS}). This approach
provides an efficient and robust way to solve the minimization
problem.  Before introducing the algorithm, we review in Section \ref{ReviewSection}
the discretization of the gradient
flow/imaginary-time equation by standard schemes (explicit/implicit
Euler and Crank-Nicolson schemes).  This enables us to make some
interesting and meaningful connections between these approaches and
some techniques related to eigenvalue problems, such as the
power method. In sections \ref{sectionGradient} and
\ref{sectionConjugateGradient}, we introduce the
projected preconditioned gradient (PG) and preconditioned conjugate
gradient (PCG) methods for solving the minimization problem on
the Riemannian manifold defined by the spherical constraints. In particular, we
provide some formulae to compute the stepsize arising in such
iterative methods to get the energy decay assumption fulfilled.
 The stopping criteria and convergence analysis are discussed in
sections \ref{sectionStopping} and \ref{sectionConvergence}. We then
investigate the design of preconditioners (section
\ref{sectionPreconditioners}). In particular, we propose a new simple
symmetrical combined preconditioner, denoted by $P_{\textrm{C}}$.  In
Section \ref{sectionNumericalResultats}, we consider the numerical
study of the minimization algorithms for the 1D, 2D and 3D GPEs
(without and with rotation).  We first propose in section
\ref{Section1DNumerics} a thorough analysis in the one-dimensional
case. This shows that the PCG approach with combined preconditioner
$P_{\textrm{C}}$ and pseudospectral discretization, called
PCG$_{\textrm{C}}$ method, outperforms all the other approaches, most
particularly for very large nonlinearities. In sections
\ref{Section2DNumerics} and \ref{Section3DNumerics}, we confirm these
properties in 2D and 3D, respectively, and show how the PCG$_{\rm C}$
algorithm behaves with respect to increasing the rotation speed.
Finally, section \ref{SectionConclusion} provides a conclusion.

\section{Definitions and notations}
\label{sec:defs}
The problem we consider is the following: find
$\phi \in L^{2}(\R^{d})$ such that
\begin{align}
  \label{eq:minimization_pb}
  \phi \in \argmin_{\norm{\phi} = 1} E(\phi).
\end{align}
We write $\displaystyle \norm{\phi} = \int_{\R^{d}} |\phi|^{2}$ for the standard
$L^2$-norm and the energy functional $E$ is defined by
\begin{align*}
  E(\phi) &= \int_{\R^d}\left[ \frac{1}{2}|\nabla\phi|^2+V(\bx)|\phi|^2+\frac{\eta}{2}|\phi|^4-\og\phi^*L_z\phi \right],
\end{align*}
where $V$ is an external potential, $\eta$ is the nonlinearity strength,
 $\og$ is the rotation speed, and $L_z=i(y \p_x-x \p_y)$ is the
angular momentum operator.

A direct computation of the gradient of the energy leads to
\begin{align*}
  \nabla E(\phi) &= 2H_{\phi} \phi, \quad\text{with} \quad H_{\phi} = -\fl{1}{2}\Delta+V+\eta|\phi|^2- \og L_{z}
\end{align*}
the mean-field Hamiltonian. We can compute the second-order derivative as
\begin{align*}
\displaystyle  \frac 1 2 \nabla^{2} E(\phi)[f, f] = \lela f, H_{\phi}f \rira +
                                       \eta \Re
                                       \lela \phi^{2}, f^{2}\rira.
\end{align*}
We introduce $
\mathcal S = \{\phi \in L^{2}(\R^{d}), \norm \phi = 1\}$ as the spherical
manifold associated to the normalization constraint. Its tangent space
at a point $\phi \in \mathcal S$ 
is $T_{\phi}\mathcal S = \{h \in L^{2}(\R^{d}), \Re \lela \phi, h\rira = 0\}$,
and the orthogonal projection $M_{\phi}$ onto this space is given by $M_{\phi} h =
h - \Re \lela \phi, h\rira \phi$.

The Euler-Lagrange equation (first-order necessary condition)
associated with the problem \eqref{eq:minimization_pb} states that, at
a minimum $\phi \in \mathcal S$, the projection of the gradient on the
tangent space is zero, which is equivalent to
\begin{align*}
  H_{\phi} \phi = \lambda \phi,
\end{align*}
where $\lambda = \lela H_{\phi} \phi, \phi\rira$ is the Lagrange
multiplier associated to the spherical constraint, and is also known
as the chemical potential. Therefore, the minimization problem can be
seen as a nonlinear eigenvalue problem. The second-order necessary
condition states that, for all $h \in T_{\phi}\mathcal S$,
$$\frac 1 2 \nabla^{2}E(\phi)[h,h] - \lambda \norm{h}^{2} \geq 0.$$ For a
linear problem ($\eta = 0$) and for problems where the nonlinearity
has a special structure (for instance, the Hartree-Fock model), this
implies that $\lambda$ is the lowest eigenvalue of $H_{\phi}$ (a
property known as the \textit{Aufbau} principle in electronic
structure theory). This property is not satisfied here.

\section{Discretization}
\label{SectionPS}
To find a numerical solution of the minimization problem, the function
$\phi \in L^{2}(\R^{d})$ must be discretized. The presence of vortices
in the solution imposes strong constraints on the discretization,
which must be accurate enough to resolve fine details. Several
discretization schemes have been used to compute the solution to the GPE, including
 high-order finite difference schemes or finite
element schemes with adaptive meshing strategies \cite{DanailaHecht1,Danaila2010}. Here, we consider a
standard pseudo-spectral discretization based on Fast Fourier
Transforms (FFTs)
\cite{AntoineDuboscqKrylovJCP,AntDubBookChapter,BaoDuSISC,ZengZhang}.

We truncate the wave function $\phi$ to a square domain $[-L, L]^{d}$,
with periodic boundary conditions, and discretize $\phi$ with the same
even number of grid points $M$ in any dimension. These two conditions can
of course be relaxed to a different domain size $L$ and number of
points $M$ in each dimension, at the price of more complex
notations. We describe our scheme in 2D, its extension to other
dimensions being straightforward.
We introduce a uniformly
sampled grid:
$\mathcal{D}_{M}:=\{ \mathbf{x}_{k_1,k_2}=(x_{k_1},y_{k_2})
\}_{(k_1,k_2)\in \mathcal{O}_{M}}$, with
$\mathcal{O}_{M}:=\{0,\ldots, M-1\}^{2}$, $x_{k_{1}+1}-x_{k_{1}}=y_{k_{2}+1}-y_{{k}_{2}}=h$, with mesh size
$h=2L/M$, $M$ an even number. We define the discrete Fourier
frequencies $ 
(\xi_p,\mu_q)$, with
$\xi_p =p{\pi}/{L}, -M/2 \leq p \leq M/2-1$, and
$\mu_q=q{\pi}/{L}, -M/2 \leq q \leq M/2-1$.  The pseudo-spectral
approximations $\widetilde{\phi }$ of the function $\phi $ in the $x$-
and $y$-directions are such that
\begin{equation*}
\begin{array}{l}
\displaystyle \widetilde{\phi }(t,x,y)= 
\frac{1}{M}\sum_{p=-M/2}^{M/2-1}\widehat{\widetilde{\phi }_{p}}(t,y)e^{i\xi_{p}(x+L)}, \hspace{2cm}
\displaystyle
  \widetilde{\phi}(t,x,y)=
   \frac{1}{M}\sum_{q=-M/2}^{M/2-1}\widehat{\widetilde{\phi}_{q}}(t,x)e^{i\mu_{q}(y+L)},
\end{array}
\end{equation*}
where $\widehat{\widetilde{\phi }_{p}}$ and $\widehat{\widetilde{\phi }_{q}}$ are respectively
 the Fourier coefficients in the $x$-  and $y$-directions
\begin{equation*}
\begin{array}{l}
\displaystyle \widehat{\widetilde{\phi }_{p}}(t,y)= \sum_{k_1=0}^{M-1}\widetilde{\phi }_{k_1}(t,y)e^{-i\xi_{p}(x_{k_1}+L)}, \hspace{2cm}
\displaystyle \widehat{\widetilde{\phi }_{q}}(t,x)= \sum_{k_2=0}^{M-1}\widetilde{\phi }_{k_2}(t,x)e^{-i\mu_{q}(y_{k_2}+L)}.
\end{array}
\end{equation*}
The following notations are used: $\widetilde{\phi }_{k_1}(t,y)=
\widetilde{\phi }(t,x_{k_1},y)$
and $\widetilde{\phi }_{k_2}(t,x)=\widetilde{\phi }(t,x,y_{k_2})$.
In order to  evaluate the operators, we introduce the matrices
\begin{equation*}\label{physeva}
\mathbb{I}_{k_1,k_2}:= \delta_{k_1,k_2}, \hspace{1cm} [[V]]_{k_1,k_2}:=V(\mathbf{x}_{k_1,k_2}), \hspace{1cm} 
 [[ |\phi |^{2}  ]]_{k_1,k_2} = |\phi _{k_1,k_2}|^{2},
\end{equation*}
for $(k_1,k_2)\in\mathcal{O}_{M}$, and  $\delta_{k_1,k_2}$ being the Dirac delta symbol which is
 equal to $1$ if and only if $k_1=k_2$ and $0$ otherwise.
We also need the operators $[[\partial^{2}_{x}]]$, $[[\partial^{2}_{y}]]$, $y [[\partial_{x}]]$ and $x [[\partial_{y}]]$
which are applied to the approximation $\widetilde{\phi }$ of $\phi $,
for $(k_1,k_2)\in\mathcal{O}_{M}$,
 \begin{equation}\label{FFTdiscret}
 \begin{array}{l}
\dsp \partial^{2}_{x}\phi (\xbf_{k_1,k_2})\approx  ([[\partial^{2}_{x}]]\widetilde{\phi })_{k_1,k_2}:= -\frac{1}{M}
\sum_{p=-M/2}^{M/2-1}\xi^{2}_{p} \widehat{(\widetilde{\phi }_{k_2})}_{p} e^{ i \xi_p(x_{k_1}+L)},\\[5mm]
\dsp \partial^{2}_{y}\phi (\xbf_{k_1,k_2})\approx  ([[\partial^{2}_{y}]]\widetilde{\phi })_{k_1,k_2}:= -\frac{1}{M}
\sum_{q=-M/2}^{M/2-1}\mu^{2}_{q} \widehat{(\widetilde{\phi }_{k_1})}_{q} e^{ i \mu_q(y_{k_2}+L)},\\[5mm]
\dsp  (x\partial_{y}\phi )(\xbf_{k_1,k_2})
\approx   (x [[\partial_{y}]]\widetilde{\phi })_{k_1,k_2}:=\frac{1}{M}
\sum_{q=-M/2}^{M/2-1} ix_{k_1} \mu_{q}\widehat{(\widetilde{\phi }_{k_1})}_{q} e^{ i \mu_q(y_{k_2}+L)},\\[5mm]
\dsp  (y \partial_{x}\phi )(\xbf_{k_1,k_2})
\approx   (y [[\partial_{x}]]\widetilde{\phi })_{k_1,k_2}:=\frac{1}{M}
\sum_{p=-M/2}^{M/2-1} iy_{k_2} \xi_{p}\widehat{(\widetilde{\phi }_{k_2})}_{p} e^{ i \xi_p(x_{k_1}+L)}.
\end{array}
\end{equation}
By considering the operators
from $\mathbb{C}^{N}$ ($N=M^{2}$ (in 2D))
 to $ \mathbb{C}$  given by
$[[\Delta]]:=  [[\partial^{2}_{x}]] + [[\partial^{2}_{y}]]$
and 
$[[\mathbb{L}_{z}]]:= -i(x [[\partial_{y }]]-y [[\partial_{x}]])$,
we obtain the discretization of the gradient of the energy 
$$
  \nabla E(\phi) = 2H_{\phi} \phi,\quad \textrm{
  with }
H_{\phi} = 
 -\frac 1 2 [[\Delta]] + [[V]] + \eta [[|\phi|^2]] - \og [[L_z]].$$
We set $\phi :=(\tilde\phi (\xbf_{k_1,k_2}))_{(k_1,k_2)\in\mathcal{O}_{M}}$ as the
 discrete unknown vector in $\mathbb{C}^{N}$.
 For  conciseness,  we identify an array $\phi$
in
the vector space of 2D complex-valued arrays
 $\mathcal{M}_{M}(\mathbb{C})$ (storage according to the 2D grid) and
 the  reshaped  vector in $\mathbb{C}^{N}$.
Finally, the cost for evaluating the application of a 2D FFT is
$\mathcal{O}(N\log N)$.

In this discretization, the computations of Section \ref{sec:defs} are
still valid, with the following differences: $\phi$ is an element of
$\C^{N}$, the operators $\Delta$, $V$, $|\phi|^{2}$ and $L_{z}$
are $N \times N$ hermitian matrices, and the inner product is the
standard $\C^{N}$ inner product. In the following sections, we will
assume a discretization like the above one, and drop the brackets in
the operators for conciseness.

\section{Review and analysis of classical methods for computing the ground states of GPEs}
\label{ReviewSection}

Once an appropriate discretization is chosen, it remains to compute
the solution to the discrete minimization problem
\begin{align}
  \label{eq:discrete_minimization}
\phi \in \argmin_{\phi \in \C^{N}, \norm{\phi} = 1} E(\phi).
\end{align}

Classical methods used to find solutions of
\eqref{eq:discrete_minimization} mainly use the so-called
\emph{imaginary time equation}, which is formally obtained by
considering imaginary times in the time-dependent Schr\"odinger
equation.
Mathematically, this corresponds to the gradient flow associated with
the energy $E$ on the manifold $\mathcal S$:
\begin{align}
  \label{eq:imaginary_time}
  \partial_{t} \phi = - \frac{1}{2}M_{\phi} \nabla E(\phi)=- (H_{\phi} \phi - \lambda(\phi) \phi).
\end{align}
As is well-known, the oscillatory behavior of the eigenmodes of the
Schr\"odinger equation become dampening in imaginary time, thus
decreasing the energy. The presence of the Lagrange multiplier
$\lambda$, coming from the projection on the tangent space on
$\mathcal S$, ensures the conservation of norm: $\norm{\phi} = 1$
for all times.
This equation can be discretized in time and solved. However, since
$H_{\phi}$ is an unbounded operator, explicit methods encounter
CFL-type conditions that limit the size of their time step, and many
authors \cite{AntoineDuboscqKrylovJCP,ABGPELabI,BaoCaiReview,BaoDipolar,BaoDuSISC,ZengZhang}Ê use a backward-Euler
discretization scheme.

We are interested in this section in obtaining the asymptotic
convergence rates of various discretizations of this equation, to
compare different approaches. To that end, consider here the model
problem of finding the first eigenpair of a $N \times N$ symmetric
matrix $H$. We label its eigenvalues
$\lambda_{1} \leq \lambda_{2} \dots \leq \lambda_{N}$, and assume that
the first eigenvalue $\lambda_{1}$ is simple, so that
$\lambda_{1} < \lambda_{2}$. This model problem is a linearized
version of the full nonlinear problem. It is instructive for two
reasons: first, any good algorithm for the nonlinear problem must also
be a good algorithm for this simplified problem. Second, this model
problem allows for a detailed analysis that leads to tractable
convergence rates. These convergence rates allow a comparison between
different schemes, and are relevant in the asymptotic regime of the
nonlinear problem.

\begin{figure}[h!]
  \centering
  \includegraphics[width=0.3\textwidth]{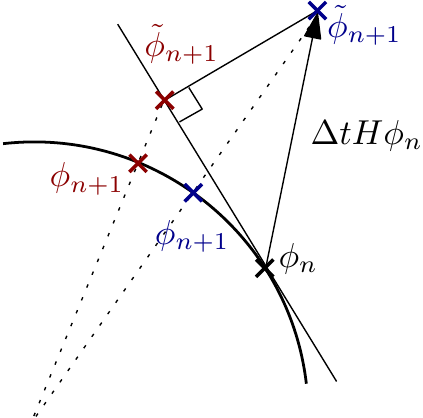}
  \caption{Update rule with the projected gradient, with the $\lambda$
    term (red) and with the unprojected gradient, without the
    $\lambda$ term (blue).}
  \label{fig:proj_grad}
\end{figure}

We consider the following discretizations of equation \eqref{eq:imaginary_time}:  Forward
Euler (FE), Backward Euler (BE), and Crank-Nicolson (CN) schemes:
\begin{align}
  \frac{\tilde \phi_{n+1}^{\text{FE}} - \phi_{n}}{\Delta t} &= -(H \phi_{n} -
                                           \lambda(\phi_{n})
                                           \phi_{n}),\\
  \frac{\tilde \phi_{n+1}^{\text{BE}} - \phi_{n}}{\Delta t} &= -(H \tilde \phi_{n+1}^{\text{BE}}
                                           - \lambda(\phi_{n})
                                           \phi_{n}),\label{BEeq}\\
  \frac{\tilde \phi_{n+1}^{\text{CN}} - \phi_{n}}{\Delta t} &= -\frac 1 2 (H\tilde \phi_{n+1}^{\text{CN}}
                                           - \lambda(\phi_{n})
                                                  \tilde \phi_{n+1}^{\text{CN}}) - \frac 1 2 (H \phi_{n} -
                                           \lambda(\phi_{n})
                                           \phi_{n}).
\end{align}
These discretizations all decrease the energy when $\Delta t>0$ is small
enough, but do not preserve the norm: the departure from normalization
is of order $\mathcal{O}(\Delta t^{2})$. Therefore, they are followed by a
projection step
\begin{align*}
  \phi_{n+1} &= \frac{\tilde \phi_{n+1}}{\norm{\tilde \phi_{n+1}}}.
\end{align*}
Note that some authors do not include the $\lambda$ term, choosing
instead to work with
\begin{align}
  \frac{\tilde \phi_{n+1}^{\text{FE}} - \phi_{n}}{\Delta t} &= -H \phi_{n},\label{FEscheme}\\
  \frac{\tilde \phi_{n+1}^{\text{BE}} - \phi_{n}}{\Delta t} &= -H \tilde \phi_{n+1}^{\text{BE}},\label{StdBEeq}\\
  \frac{\tilde \phi_{n+1}^{\text{CN}} - \phi_{n}}{\Delta t} &= -\frac 1 2 (H\tilde
                                                  \phi_{n+1}^{\text{CN}}
                                                  + H \phi_{n}).
\end{align}
These also yield schemes that decrease the energy. However, because of
the use of the unprojected gradient $H \phi$ instead of
$M_{\phi} (H \phi) = H \phi - \lambda \phi$, the departure from
normalization is $\mathcal{O}(\Delta t)$, instead of $\mathcal{O}(\Delta t^{2})$ for the
projected case. The difference between the two approaches is
illustrated in Figure \ref{fig:proj_grad}.

For the FE and BE methods, a simple algebraic
manipulation shows that one step of the method with the $\lambda$ term
is equivalent to one step of the method without the $\lambda$ term,
but with an effective $\Delta t$ modified as
\begin{align}
  \Delta_{t}^{\lambda, \text{FE}} &= \frac{\Delta t}{1 + \Delta t \lambda},\\
  \Delta_{t}^{\lambda, \text{BE}} &= \frac{\Delta t}{1 - \Delta t
                                    \lambda}. \label{eq:eq_BE_BElambda}
\end{align}
This is not true for the CN method, nor it is true when nonlinear
terms are included. However, even in this case, the difference between
including and not including the $\lambda$ term is $\mathcal{O}(\Delta t^{2})$,
and their behavior is similar. Since the analysis of unprojected
gradient methods is simpler, we focus on this here.

Then, these schemes can all be written in the form
\begin{align}
  \label{eq:power_iter}
  \phi_{n+1} &= \frac{A \phi_{n}}{\norm{A \phi_{n}}},
\end{align}
where the matrix $A$ is given by
\begin{align*}
  A^{\text{FE}} &= I - \Delta t H,\\
  A^{\text{BE}} &= (I + \Delta t H)^{-1},\\
  A^{\text{CN}} &= (I + \frac {\Delta t}2 H)^{-1}
                      (I - \frac{\Delta t} 2 H).
\end{align*}
The
eigenvalues $\mu_{i}$ of $A$ are related to the eigenvalues
$\lambda_{i}$ of $H$ by the following spectral transformation
\begin{align*}
  \mu_{i}^{\text{FE}} &= 1 - \Delta t \lambda_{i},\\
  \mu_{i}^{\text{BE}} &= \frac{1}{1 + \Delta t \lambda_{i}},\\
  \mu_{i}^{\text{CI}} &= \frac{1 - \frac {\Delta t}2 \lambda_{i}}{1 + \frac{\Delta t}2 \lambda_{i}}.
\end{align*}
We call this eigenvalue the amplification factor: if $\phi_{n}$ has
eigencomponents $c_{n,i} = \lela v_{i}, \phi_{n}\rira$ on the
eigenvector $v_{i}$ of $H$ and $\phi_{0}$ is normalized to $1$, then
the iteration \eqref{eq:power_iter} can be solved as
\begin{align*}
  c_{n,i} = \frac{\mu_{i}^{n}}{\sqrt{\sum_{i=1}^{N}
  |\mu_{i}^{n}c_{0,i}|^{2}}} \,c_{0,i}.
\end{align*}
This iteration converges towards the eigenvector associated to the
largest eigenvalue $\mu_{N}$ (in modulus) of $A$, if it is simple, with
convergence rate $\mu_{N-1}/\mu_{N}$ . This is
nothing but the classical power method for the computation of
eigenvalues, with a spectral transformation from $H$ to $A$. Therefore
we identify the FE method as a shifted power method, the BE method as
a shift-and-invert approach, and the CN uses a generalized Cayley
transform \cite{bai2000templates,saad2011numerical}.

From this we can readily see the properties of the different
schemes. We make the assumption that either $\lambda_{1}$ is positive
or that $\Delta t < \frac 1 {-\lambda_{1}}$ (for BE) and
$\Delta t < \frac 2 {-\lambda_{1}}$ (for CN). If this condition is not
verified, then the iteration will generally not converge towards an
eigenvector associated with $\lambda_{1}$ because another eigenvalue
than $\lambda_{1}$ will have a larger amplification factor. Under this
assumption, we see that the BE and CN converge unconditionally, while
FE only converges if
\begin{align*}
  \Delta t < \frac 2 {\lambda_{N}}.
\end{align*}
This is a CFL-like condition: when $H$ is the discretization of an
elliptic operator, $\lambda_{N}$ will tend to
infinity as the size of the basis increases, which will force FE to
take smaller time steps.

The asymptotic convergence rate of these methods is
$\frac{\mu_{N-1}}{\mu_{N}}$. While the FE method has a bounded
convergence rate, imposed by $\lambda_{1}, \lambda_{2}$ and
$\lambda_{N}$, the BE and CN methods can be made to have an
arbitrarily small convergence rate, by simply choosing $\Delta t$
arbitrarily close to $-\frac {1} {\lambda_{1}}$ (BE) or $-\frac {2} {\lambda_{1}}$
(CN). Since in practice $\lambda_{1}$ is unknown, it has to be
approximated, for instance by $\lambda(\phi_{n})$. This yields the
classical Rayleigh Quotient Iteration:
\begin{align*}
  \phi_{n+1} &= \frac{(H - \lela \phi_{n}, H \phi_{n}\rira)^{-1}
                      \phi_{n}}{\norm{(H - \lela \phi_{n}, H \phi_{n}\rira)^{-1}
                      \phi_{n}}},
\end{align*}
which is known to converge cubically. This iteration can also be seen
as a Newton-like method.

From the previous considerations, it would seem that the BE and CN are
infinitely superior to the FE method: even with a fixed stepsize, the
BE and CN methods are immune to CFL-like conditions, and with an
appropriately chosen stepsize, it can be turned into a
superlinearly-convergent scheme. The first difficulty with this
approach is that it is a linear strategy, only guaranteed to converge
when close to the ground state. As is always the case with Newton-like
methods, it requires a globalization strategy to be efficient and
robust in the nonlinear setting. The second issue, is that the BE and
CN method require the solution of the linear system
$(H - \lambda(\phi_{n})) \tilde \phi_{n+1} = \phi_{n}$.

The difficulty of solving the system depends on the discretization
scheme used. For localized basis schemes like the Finite Element
Method, $H$ is sparse, and efficient direct methods for large scale
sparse matrices can be used \cite{Saa03}. For the Fourier
pseudo-spectral scheme, which we use in this paper, $H$ is not sparse,
and only matrix-vector products are available efficiently through FFTs
(a \textit{matrix-free problem}). This means that the system has to be
solved using an iterative method. Since it is a symmetric but
indefinite problem, the solver of choice is MINRES
\cite{barrett1994templates}, although the solver BICGSTAB has been
used \cite{AntoineDuboscqKrylovJCP,AntDubBookChapter}. The number of
iterations of this solver will then grow when the grid spacing tends
to zero, which shows that BE also has a CFL-like limitation. However,
as is well-known, Krylov methods \cite{Saa03} only depend on the
square root of the condition number for their convergence, as opposed
to the condition number itself for fixed-point type methods
\cite{BaoDuSISC,ZengZhang}. This explains why BE with a Krylov solver
is preferred to FE in practice
\cite{AntoineDuboscqKrylovJCP,AntDubBookChapter}.

Furthermore, preconditioners can be used to reduce this number of
iterations \cite{AntoineDuboscqKrylovJCP,AntDubBookChapter}, e.g. with
a simple preconditioner (one that is diagonal either in real or in
Fourier space). This method is effective for many problems, but
requires a globalization strategy, as well as an appropriate selection
of parameters such as $\Delta t$ and the precision used to solve the
linear system \cite{AntoineDuboscqKrylovJCP,AntDubBookChapter}. Here,
we propose a method that has all the advantages of BE (robust,
Krylov-like dependence on the square root of the condition number,
ability to use a preconditioner), but is explicit, converges faster
than BE, and has no free parameter (no fine-tuning is necessary).

\section{The Preconditioned nonlinear Gradient (PG) and Conjugate Gradient (PCG) methods}

\subsection{The gradient method}\label{sectionGradient}
The previous approaches usually employed in the literature to compute
the ground states of the Gross-Pitaevskii equation are all based on implicit
discretizations of the imaginary-time equation
\eqref{eq:imaginary_time}. As such, these methods come from PDE theory
and lack the insight of minimization algorithms. The difficulty of
applying classical minimization algorithms comes from the spherical
constraints. However, the general theory of optimization algorithms on
Riemannian manifolds has been developed extensively in
\cite{absil2009optimization,edelman1998geometry}, where the authors
derive constrained analogues of gradient, conjugate gradient and
Newton algorithms. This is the approach we follow here, and employ a
preconditioned conjugate gradient method on the manifold $\mathcal
S$.

In this section, we work with an arbitrary symmetric positive definite
preconditioner $P$. The choice of $P$ will be discussed later in subsection \ref{sectionPreconditioners}.
The (projected, preconditioned) gradient method for the minimization of $E$ on
$\mathcal S$ is the update
\begin{align}
  \label{eq:grad_update_1}
    \tilde \phi_{n+1} &= \phi_{n} - \alpha_{n} P \left( H_{\phi_{n}} \phi_{n} -
                      \lambda_{n} \phi_{n}\right), \quad \phi_{n+1} = {\tilde \phi_{n+1}}/{\norm{\tilde \phi_{n+1}}},
\end{align}
where $\lambda_{n}=\lambda(\phi_n).$
We reformulate this equation as
\begin{align}
    \label{eq:grad_update_2}
  \phi_{n+1} &= \cos(\theta_{n}) \phi_{n} + \sin(\theta_{n})
               \frac{p_{n}}{\norm{p_{n}}}, \quad \text{with}\quad  p_{n} = {d_{n} - \Re\lela d_{n}, \phi_{n}\rira \phi_{n}},
\end{align}
where $d_{n} = -P r_{n}$ is the descent direction, equal to the
negative of the preconditioned residual
$r_{n} = H_{\phi_{n}} \phi_{n} - \lambda(\phi_{n})\phi_{n}$. The
equations \eqref{eq:grad_update_1} and \eqref{eq:grad_update_2} are
equivalent when $\theta_{n}$ or $\alpha_{n}$ is small enough, with a
one-to-one correspondance between $\theta_{n}$ and $\alpha_{n}$. To
first order, we have: $\alpha_{n} = \theta_{n} \norm{p_{n}}$.
 
   Without preconditioner, this method, summarized in
Algorithm \ref{alg:G}, is identical to the FE method (\ref{FEscheme}).

\RestyleAlgo{boxed}
\begin{algorithm}[h!]
  \caption{The gradient method\label{alg:G}}
  \While{not converged}{
    $\lambda_{n} = \lambda(\phi_{n})$
    
    $r_{n} = H_{\phi_{n}} \phi_{n} - \lambda_{n} \phi_{n}$
    
    $d_{n} = - P r_{n}$
    
    $p_{n} = {d_{n} - \Re\lela d_{n}, \phi_{n}\rira \phi_{n}}$
 
    $\theta_{n}  = \argmin_{\theta} E\left(\cos(\theta) \phi_{n} + \sin(\theta) {p_{n}}/{\norm{p_{n}}}\right)$

    $\phi_{n+1} = \cos(\theta_{n}) \phi_{n} + \sin(\theta_{n})
    {p_{n}}/{\norm{p_{n}}}$

    $n = n+1$
}
\end{algorithm}

To choose the parameter $\theta_{n}$, a number of strategies are
possible. We first show that, when $\theta_{n}$ is small enough, the
gradient method decreases the energy.

Expanding $\phi_{n+1}$ up to second-order in $\theta_{n}$, we obtain
\begin{align}
  \phi_{n+1}&= \left(1-\frac{\theta_{n}^{2}}{2}\right)\phi_{n} + \theta_{n} \frac{p_{n}}{\norm{p_{n}}}
              + \mathcal{O}(\theta_{n}^{3}),
\end{align}
and therefore
\begin{align}
  E(\phi_{n+1}) &= E(\phi_{n}) + \frac{\theta_{n}}{\norm{p_{n}}}\Re\lela \nabla
                  E(\phi_{n}), p_{n}\rira + \frac 1
                  2\frac{\theta_{n}^{2}}{\norm{p_{n}}^{2}}\left[\nabla^{2}E(\phi_{n})[p_{n},
                  p_{n}] - \lambda_{n} \norm{p_{n}}^{2}\right] + \mathcal{O}(\theta_{n}^{3}).
     \label{eq:expansion_energy}
\end{align}
We now compute the first-order variation
\begin{align*}
  \Re\lela \nabla E(\phi_{n}), p_{n}\rira &= \frac {\Re \lela \nabla
                                            E(\phi_{n}), d_{n} - \Re\lela d_{n},
              \phi_{n}\rira \phi_{n}\rira} {\norm{{d_{n} - \Re\lela d_{n},
                                            \phi_{n}\rira \phi_{n}}}}
  = \frac {\Re \lela r_{n}, d_{n}\rira} {\norm{{d_{n} - \Re\lela d_{n},
                                            \phi_{n}\rira \phi_{n}}}}\\
  &= -\frac {\lela r_{n}, P r_{n}\rira} {\norm{{d_{n} - \Re\lela d_{n},
                                            \phi_{n}\rira \phi_{n}}}}.
\end{align*}
Since $P$ was assumed to be positive definite, this term is always
negative so that the algorithm decreases the energy when $\theta_{n}$
is chosen small enough. Since $p_{n}$ is orthogonal to
$\phi_{n}$, the second-order term
$\nabla^{2}E(\phi_{n})[p_{n}, p_{n}] - \lambda_{n} \norm{p_{n}}^{2}$ is guaranteed to
be positive when $\phi_{n}$ is close to a minimizer by the second-order optimality conditions.

Therefore, a basic strategy is to choose $\theta_{n}$ fixed and small
enough so that the energy decreases. A better one is to choose
$\theta_{n}$ adaptively. For instance, we could perform the linesearch
\begin{align}
\label{Opt_Theta}
\theta_{n}  = \argmin_{\theta} E\left(\cos(\theta) \phi_{n} + \sin(\theta) \frac{p_{n}}{\norm{p_{n}}}\right).
\end{align}
Since $E(\theta)$ is not a quadratic function,  this is a nonlinear
one-dimensional minimization problem, generally requiring many
evaluations of $E(\theta)$ to converge to a minimum. However, many of the
computations for the evaluation of $E(\theta)$, including all that
require FFTs, can be pre-computed. Since the FFT step is the dominant
one in the computation of the energy, the evaluation of $E(\theta)$ at
many points is not much more costly than the evaluation at a single
point. Therefore it is feasible to use a standard one-dimensional
minimization routine.

Alternatively, we can obtain a simple and cheap approximation by
minimizing the second-order expansion of $E$ in $\theta_{n}$. We
expect this to be accurate when $\theta_{n}$ is small, which is the
case close to a minimizer. Minimizing \eqref{eq:expansion_energy} with
respect to $\theta_{n}$ yields
\begin{align}
  \label{eq:theta_opt}
  \theta_{n}^{\rm opt} &= \frac{-\Re\lela \nabla
                  E(\phi_{n}), p_{n}\rira\norm{p_{n}}}{\Re\left[\nabla^{2}E(\phi_{n})[p_{n},
                  p_{n}] - \lambda_{n} \right]}.
\end{align}
As we have seen, the numerator is always positive, and the denominator
is positive when $\phi_{n}$ is close enough to a minimizer. In our
implementation, we compute the denominator, and, if it is positive, we use
$\theta_{n}^{\rm opt}$ as a trial stepsize. If not, we use some
default positive value. If the energy of $\phi_{n+1}$ using this trial
stepsize is decreased, we accept the step. If the energy is not
decreased, we reject the step, decrease the trial stepsize, and try
again, until the energy is decreased (which is mathematically ensured
when $\theta_{n}$ is small enough). Alternatively, we can use Armijo
or Wolfe conditions as criterion to accept or reject the stepsize, or
even use the full line search \eqref{Opt_Theta}. The evaluation of the
energy at multiple values of $\theta$ do not require more Fourier
transforms than at only one point, but only more computations of the
nonlinear term, so that a full line search is not much more costly
than the heuristic outlined above. In our tests however, the heuristic
above was sufficient to ensure fast convergence, and a full line
search only marginally improved the number of iterations. Therefore, we simply
use the heuristic in the numerical results of Section
\ref{sectionNumericalResultats}.

Let us note that under reasonable assumptions on the structure of
critical points and on the stepsize choice, there are various results
on the convergence of this iteration to a critical point (see
\cite{absil2009optimization} and references therein).

\subsection{The conjugate gradient method}\label{sectionConjugateGradient}
The conjugate gradient method is very similar, but uses an update
rule of the form
\begin{align}
  \label{eq:constrained_CG}
  d_{n} &= - P r_{n} + \beta_{n} p_{n-1}
\end{align}
instead of simply $d_{n} = - P r_{n}$. This is justified when
minimizing unconstrained quadratic functionals, where the formula
\begin{align}
  \label{CG:FR}
    d_{n} = - P r_{n} + \beta_{n} d_{n-1}, \quad \text{with} \quad \beta_{n}  = \frac{\lela r_{n}, P r_{n}\rira}{\lela r_{n-1}, P
                r_{n-1}\rira},
\end{align}
yields the well-known PCG method to solve linear systems. For
nonlinear problems, different update formulas can be used, all
equivalent in the linear case. Equation (\ref{CG:FR}) is known as the
Fletcher-Reeves update. Another popular formula is the Polak-Ribi\`ere
choice $\beta = \max(\beta^{\textrm{PR}}, 0)$, where
\begin{align}
  \label{CG:PR}
  \beta^{\textrm{PR}} & = \frac{\lela r_{n} - r_{n-1}, P r_{n}\rira}{\lela r_{n-1}, P
                r_{n-1}\rira}.
\end{align}
We use $\beta = \max(\beta^{\textrm{PR}}, 0)$, which is equivalent to
restarting the CG method (simply using a gradient step) when
$\beta^{\textrm{PR}} < 0$ and is a standard choice in nonlinear CG
methods. For the justification of the CG method for constrained
minimization, see \cite{absil2009optimization,edelman1998geometry}.

\begin{algorithm}
  \caption{The conjugate gradient
    method}\label{alg:CG}
  \While{not converged}{
  $\lambda_{n} = \lambda(\phi_{n})$
  
  $r_{n} = H_{\phi_{n}}\phi_{n} - \lambda_{n} \phi_{n}$
  
  $\beta_{n} = {\lela r_{n} - r_{n-1}, P r_{n}\rira}/{\lela
    r_{n-1}, Pr_{n-1}\rira}$
  
  $\beta_{n} = \max(\beta_{n}, 0)$
  
  $d_{n} = - P r_{n} + \beta p_{n-1}$
  
  $p_{n} = d_{n} - \Re\lela d_{n}, \phi_{n}\rira \phi_{n}$
  
    $\theta_{n}  = \argmin_{\theta} E\left(\cos(\theta) \phi_{n} +
      \sin(\theta) {p_{n}}/{\norm{p_{n}}}\right)$
    
  $\phi_{n+1} = \cos(\theta_{n}) \phi_{n} +
  \sin(\theta_{n}){p_{n}}/{\norm{p_{n}}}$

      $n = n+1$
  }
\end{algorithm}

The CG algorithm is presented in Algorithm \ref{alg:CG}. In contrast
with the gradient algorithm, the quantity
$\Re\lela \nabla E(\phi_{n}), p_{n}\rira$ does not have to be
negative, and $p_{n}$ might not be a descent direction: even with a
small stepsize, the energy does not have to decrease at each step. To
obtain a robust minimization method, we enforce energy decrease to
guarantee convergence. Therefore, our strategy is to first check if
$p_{n}$ is a descent direction by computing
$\Re\lela \nabla E(\phi_{n}), p_{n}\rira$. If $p_{n}$ is not a descent
direction, we revert to a gradient step, which we know will decrease
the energy, else, we choose $\theta_{n}$ as in \eqref{eq:theta_opt},
and use the same stepsize control as in the gradient algorithm.

In our numerical tests, we observe that these precautions of checking
the descent direction and using a stepsize control mechanism are
useful in the first stage of locating the neighborhood of a
minimum. Once a minimum is approximately located, $p_{n}$ is always a
descent direction and the stepsize choice \eqref{eq:theta_opt} always
decreases the energy.

\subsection{Stopping criteria}
\label{sectionStopping}

A common way  to terminate the iteration (in the BE schemes) is to use the stopping criterion 
\be
\label{Stop_Max}
\phi_{ \rm err}^{n,\infty}:=\|\phi_{n+1}-\phi_{n}\|_{\infty}\le\eps.
\ee
This can be problematic because the minima are generally not
isolated but form a continuum due to symmetries (for instance,
complex phase or rotational invariance), and this criterion might be
too restrictive. A more robust one is based on the norm of the
symmetry-covariant residual
\be
\label{Stop_Residual}
r^{n,\infty}_{ \rm err}:=\|r_n\|_\infty=\|H_{\phi_n}\phi_n-\lambda_n\phi_n\|_\infty\le\eps.
\ee
or the symmetry-invariant energy difference
\be\label{Stop_Energy}
\mathcal{E}^{n}_{ \rm err}:=|E(\phi_{n+1})-E(\phi_{n})|\le\eps.
\ee

This third one converges more rapidly than the two previous ones: as
is standard in optimization, when $\phi^{*}$ is a minimum of the
constrained minimization problem and $\phi \in \mathcal S$, then
\begin{align*}
  E(\phi) - E(\phi^{*}) = \mathcal{O}(\|\phi - \phi^{*}\|^{2}).
\end{align*}
This is consistent with our results in Figure
\ref{fig:og0_2d_Diff_Beta_Comp_PG_PCG}.

In the current paper, we always use the energy based stopping
criterion (\ref{Stop_Energy}): for the 2D and 3D cases, a criteria
based on $\phi^{n,\infty}_{\rm err}$ or $r^{n,\infty}_{\rm err}$ can
lead to long computational times, most particularly for large
rotations $\og$, even without changing the energy (see the example in
subsection \ref{Section2DNumerics}).

\subsection{Convergence analysis}\label{sectionConvergence}
A full analysis of the convergence properties of our methods is
outside the scope of this paper, but we give in this section some elementary
properties, and heuristics to understand their asymptotic performance.

Based on the expansion of the energy \eqref{eq:expansion_energy} as a
function of $\theta$ for the gradient method, it is straightforward to
prove that, when $E$ is bounded from below and the stepsize $\theta$
is chosen optimally, the norm of the projected gradient $H_{\phi_{n}}
\phi_{n} - \lambda_{n} \phi_{n}$ converges to $0$. Convergence
guarantees for the conjugate gradient method are harder, but can still
be proven under a suitable restart strategy that ensures that the
energy always decreases fast enough (for instance, the Armijo rule).

With additional assumptions on the non-degeneracy of critical points,
we can even prove the convergence of $\phi_{n}$ to a critical point,
that will generically be a local minimum. However, the question of the
precise convergence speed of the gradient and conjugate gradient
algorithms we use is problematic, because of three factors: the
constraint $\norm{\phi} = 1$, the non-quadraticity of $E$, and the
presence of a preconditioner. To our knowledge, no asymptotically
optimal bound for this problem has been derived. Nevertheless, based
on known results about the convergence properties of the conjugate
gradient method for preconditioned linear systems on the one hand
\cite{Saa03}, and of gradient methods for nonlinear constrained
minimization \cite{absil2009optimization} on the other, it seems
reasonable to expect that the convergence will be influenced by the
properties of the operator
\begin{align}
  \label{op_precond}
  M = (1 - \phi \phi^{*}) P (\nabla^{2}E(\phi) - \lambda I)(1 - \phi
\phi^{*}),
\end{align}
where $\phi$ is the minimum and
$\lambda = \lambda(\phi)$. This operator admits $0$ as its lowest
eigenvalue, associated with the eigenvector $\phi$. It is reasonable
to expect that the convergence rate will be determined by a condition
number $\sigma$ equals to the ratio of the largest to the lowest non-zero
eigenvalue of this operator. As is standard for linear systems, we
also expect that the number of iterations to achieve a given tolerance
will behave like $\sqrt \sigma$ for the conjugate gradient algorithm,
and $\sigma$ for the gradient algorithm. As we will see in Section
\ref{sectionNumericalResultats}, this is verified in our tests.

The Hessian operator $\nabla^{2}E(\phi)$, which includes a Laplace
operator, is not bounded. Correspondingly, on a given discretization
domain, when the grid is refined, the largest eigenvalues of this
operator will tend to $+\infty$. For a linear meshsize $h$, the
eigenvalues of $\nabla^{2}E(\phi)$ will behave as
$\mathcal{O}(h^{-2})$. This is another instance of the CFL condition
already seen in the discretization of the imaginary time equation. The
Hessian $\nabla^{2}E$ also includes a potential term $V$, which is
often confining and therefore not bounded, such as the classical
harmonic potential $V(\textbf{x}) = |\textbf{x}|^{2}$, or more
generally confining potentials whose growth at infinity is like
$|\textbf{x}|^{p}$. for some $p > 0$ Thus, even with a fixed meshsize
$h$ on a domain $[-L,L]^{d}$, when $L$ is increased, so will the
largest eigenvalues of $\nabla^{2}E(\phi)$, with a
$\mathcal{O}(L^{p})$ growth. When a (conjugate) gradient method is
used without preconditioning, the convergence will be dominated by
modes associated with largest eigenvalues of $\nabla^{2}E$. This
appears in simulations as high-frequency oscillations and localization
at the boundary of the domain of the residual $r_{n}$.

To remedy these problems and achieve a good convergence rate, adequate
preconditioning is crucial.

\subsection{Preconditioners}\label{sectionPreconditioners}
We consider the question of building preconditioners $P$ for the
algorithms presented above. In the schemes based on the discretization
of the gradient flow, preconditioning is naturally needed when solving
linear systems by iterative methods. In the gradient and conjugate
gradient optimization schemes, it appears as a modification of the
descent direction to make it point closer to the minimum:
\begin{equation}
d_{n}:=-P(H_{\phi_{n}}\phi_{n}-\lambda_{n}\phi_{n}).
\end{equation}
In both cases, the preconditioning matrix should be an approximation of
the inverse of the Hessian matrix of the problem.

\paragraph{Kinetic energy preconditioner}
One of these approximations is to use only the kinetic energy term
\begin{align}\label{PDelta}
  P_{\Delta} &= (\alpha_{\Delta} - \Delta/2)^{-1},
\end{align}
where $\alpha_{\Delta}$ is a positive shifting constant to get an
invertible operator, and $I$ is the identity operator. This has been
called a ``Sobolev gradient'' in \cite{Danaila2010} because it is
equivalent to taking the gradient of the energy in the Sobolev
$H^{1}$-norm (with $\alpha_{\Delta}=1/2$).  In the framework of the
BESP scheme for the GPE with Krylov solver, a similar preconditioner
has been proposed in \cite{AntoineDuboscqKrylovJCP}, $\alpha_{\Delta}$
being the inverse of the time step $\Delta t$ of the semi-implicit
Euler scheme.  A closely-related variant is standard in plane-wave
electronic structure computation \cite{zhou2015preconditioning}, where
it is known as the Tetter-Payne-Allan preconditioner
\cite{teter1989solution}. This preconditioner is diagonal in Fourier
space and can therefore be applied efficiently in our pseudo-spectral approximation scheme.
On a fixed domain $[-L,L]^{d}$, the effect of this preconditioner is
to make the number of iterations independent from the spatial
resolution $h$, because $P \nabla^{2}E(\phi)$, seen as an operator on
the space of functions on $[-L,L]^{d}$, will be equal to the
identity plus a compact operator. This is supported by numerical
experiments in Section \ref{sectionNumericalResultats}.
However, this operator is not bounded in the full domain
$\R^{d}$. Therefore, as $L$ increases, so will the largest eigenvalues
of $P \nabla^{2}E(\phi)$. For a potential $V$ that grows at infinity
like $|\mathbf{x}|^{p}$, the largest eigenvalues of $P \nabla^{2}E(\phi)$ are
$\mathcal{O}(L^{p})$, resulting in an inefficient preconditioner. Similarly,
when $\eta$ is large, the nonlinear term becomes dominant, and the
kinetic energy preconditioner is inefficient.

The choice of $\alpha_{\Delta}$ is a compromise: if $\alpha_{\Delta}$
is too small, then the preconditioner will become close to indefinite,
which can produce too small eigenvalues in the matrix
\eqref{op_precond} and hamper convergence. If $\alpha_{\Delta}$ is too
big, then the preconditioner does not act until very large
frequencies, and large eigenvalues result.  {We found that a suitable
  adaptive choice, that has consistently good performance and avoids
  free parameters, is
\begin{equation}
\label{Opt_Alpha}
\alpha_\Delta=\tilde{\lambda}_n:=\int\bigg(\fl{1}{2}|\nabla \phi_n|^2+ V |\phi_n|^2+\eta |\phi_n|^4\bigg) d\bx >0
\end{equation}
which is a positive number that represents the characteristic energy
of $\phi_{n}$. We use this choice for our numerical simulations.}

\paragraph{Potential energy preconditioner}
Another natural approach is to use the potential energy term for the
preconditioner:
\begin{align}\label{PV}
  P_{V} &= (\alpha_{V} + V + \eta |\phi_{n}|^{2})^{-1}.
\end{align}
This preconditioner is diagonal in real space and can therefore be
applied efficiently. Dual to the previous case, this preconditioner
has a stable performance when the domain and $\eta$ are increased, but
deteriorates as the spatial resolution is increased. Such a
preconditioner has been used in \cite{AntoineDuboscqKrylovJCP} when
the gradient flow for the GPE is discretized through a BE scheme,
leading then to a Thomas-Fermi preconditioner. In this study, the
parameter $\alpha_{V}$ was $1/\Delta t$. As in the kinetic energy
case, we found it efficient to use $\alpha_{V} = \tilde{\lambda}_{n}$,
and we will only report convergence results for this choice of
parameter.

\paragraph{Combined preconditioner}
In an attempt to achieve a stable performance independent of the size
of the domain or the spatial resolution, we can define the combined
preconditioners
\begin{align}\label{combinedP}
P_{\textrm{C}_1}&=P_{V} P_{\Delta},\quad P_{\textrm C_2}=P_{\Delta}P_{V}
\end{align}
or a symmetrized version
\begin{align}
\label{combinedP_12}
  P_{\textrm{C}} &= 
   P_{V}^{1/2}
   P_{\Delta}
   P_{V}^{1/2}.
\end{align}
With these preconditioners, $P \nabla^{2}E(\phi)$ is
bounded as an operator on $L^{2}(\R^{d})$ (this can be proven by
writing explicitly its kernels in Fourier space and then using
Schur's test). However, we found numerically that this operator is not
bounded away from zero, and has small eigenvalues of size
$\mathcal{O}(L^{-p} + h^{2})$. Therefore, the conditioning deteriorates as
both the spatial resolution and the size of the domain increase.

In summary, for a spatial resolution $h$ and a domain size $L$, the
asymptotic condition numbers of the preconditioned Hessian with these
preconditioners are
\begin{equation}
\begin{array}{l}
  \kappa_{\Delta} = \mathcal{O}(L^{p}),\\
  \kappa_{V} = \mathcal{O}(h^{-2}),\\
  \kappa_{\textrm C} = \mathcal{O}\left(\frac 1 {L^{-p} + h^{2}}\right) = \mathcal{O}(\min(L^{p}, h^{-2})).
\end{array}
\end{equation}
Therefore, the combined preconditioners act asymptotically as the best
of both the kinetic and potential preconditioners. However, they might
not be more efficient in the pre-asymptotic regime and require
additional Fourier transforms.
\paragraph{Computational efficiency}

The application the operator $P_{V}$ is almost free (since it only
requires a scaling of $\phi$), but the naive application of
$P_{\Delta}$ requires a FFT/IFFT pair. However, since we apply the
preconditioners after and before an application of the Hamiltonian,
we can reuse FFT and IFFT computations, so that the application of
$P_{\Delta}$ does not require any additional Fourier transform.
Similarly, the use of $P_{\rm C_{1}}$ and $P_{\rm C_{2}}$ only require
one additional Fourier transform per iteration, and that of the
symmetrized version $P_{\rm C}$ two.

In summary, the cost in terms of Fourier transforms per iteration for
the rotating GPE model is
 \begin{itemize}
\item no preconditioner: 3 FFTs/iteration (get the Fourier transform
  of $\phi$, and two IFFTs to compute $\Delta \phi$ and $L_{z} \phi$ respectively),
\item $P_\Delta$ or $P_V$: 3 FFTs/iteration,
\item non-symmetric combined $P_{\rm C_{1}}$ or $P_{\rm C_{2}}$: 4 FFTs/iteration,
\item symmetric combined $P_{\rm C}$: 5 FFTs/iteration.
\end{itemize}
Note that this total might be different for another type of GPE model e.g. when
a nonlocal dipole-dipole interaction is included \cite{AntDubBookChapter,BJTZ2015}.

As we will see in Section \ref{sectionNumericalResultats}, all
combined preconditioners have very similar performance, but the
symmetrized one might be more stable in some circumstances. A
theoretical explanation of these observations, and in particular of the
effect of a non-symmetric preconditioner is, to the best of our
knowledge, still missing.
\section{Numerical results}
\label{sectionNumericalResultats}

{
We first introduce some notations.  When we combine one of the
preconditioners $P_\nu$
($\nu=I$, $\Delta$, $V$, $\textrm{C}$, ${\rm C_1}$, ${\rm C_2}$)
\eqref{PDelta}-\eqref{combinedP_12} with the gradient algorithm (Alg:
\ref{alg:G}), we denote the resulting methods by PG$_{\nu}$.  Similarly, we denote
by PCG$_{\nu}$ if the preconditioned
conjugate gradient algorithm (Alg: \ref{alg:CG}) was applied. In the
following, we denote by $\#\rm iter$ the number of global iterations
for an iterative algorithm to get the converged solution with an
\textit{a priori} tolerance $\eps$ with respect to the stopping
criterion (\ref{Stop_Energy}).}

Concerning the BESP schemes \eqref{BEeq} and \eqref{StdBEeq}, at each
outer iteration $n$, one needs to solve an implicit system with the
operator $( 1/\Delta t+H_{\boldsymbol{\phi}_{n}})$. We use a Krylov
subspace iterative solver (MINRES here) with one of the
preconditioners $P_\nu$ ($\nu=I, \Delta,V,\textrm{C}$)
\eqref{PDelta}-\eqref{combinedP} to accelerate the number of inner
iterations \cite{AntoineDuboscqKrylovJCP}.  The preconditioned BESP
schemes is then denoted by BE$_{\nu}$, according to the chosen
preconditioner. The number of iterations reported is equal to the sum
of the inner iterations over the outer iterations.
 
In the following numerical experiments, we consider two types of
trapping potential $V(\bx)$: the harmonic plus lattice potential
\cite{BaoCaiReview} \be
\label{Lattice_Poten}
V(\bx)=V_d^0(\bx)+\left\{
\begin{array}{l}
\kappa_x\sin^2(q_xx^2),	\\[0.4em]
\sum_{\nu=x,y}\kappa_\nu\sin^2(q_\nu\nu^2), \\[0.4em]
\sum_{\nu=x,y,z}\kappa_\nu\sin^2(q_\nu\nu^2),
\end{array}
\right.
\quad
{\rm with}
\quad
V_d^0(\bx)=
\left\{
\begin{array}{ll}
\gm_x^2x^2,& d=1,	\\[0.4em]
\sum_{\nu=x,y}\gm_\nu\nu^2,&	d=2, \\[0.4em]
\sum_{\nu=x,y,z}\gm_\nu\nu^2,& d=3,
\end{array}
\right.
\ee
and the harmonic plus quartic potential for $d=2,3$ \cite{DanailaHecht1,Danaila2010,ZengZhang}
\be
\label{Quart_Poten}
V(\bx)=(1-\alpha)V_2^0(\bx) + \fl{\kappa\, (x^2+y^2)^2}{4} +
\left\{
\begin{array}{ll}
0, &d=2,	\\[0.4em]
\gm_z^2\,z^2, & d=3.
\end{array}
\right.
\ee
Moreover, unless stated otherwise, we  take the initial data as the
Thomas Fermi approximation \cite{AntoineDuboscqKrylovJCP, BaoCaiReview}:
  \be\label{TF_App}
\phi_0=\fl{\phi^{\rm TF}_g}{\|\phi^{\rm TF}_g\|},\qquad
{\rm with}\qquad
\phi^{\rm TF}_g(\bx)=\left\{\begin{array}{ll}
 \sqrt{\big(\mu_g^{\rm TF}-V(\bx)\big)/\eta}, & V(\bx)<\mu_{g}^{\rm TF},  \\[0.4em]
  0,&  {\rm otherwise},
\end{array}
\right.
\ee
where
 \be
\label{TF_Rdius}
\mu_g^{\rm TF}=\fl{1}{2}\left\{
\begin{array}{ll}
(3\eta\gm_x)^{2/3},& d=1,  \\
(4\eta\gm_x\gm_y)^{1/2},&d=2,  \\
(15\eta\gm_x\gm_y\gm_z)^{2/5},&d=3.
\end{array}
\right.
\ee
The algorithms were implemented  in Matlab (Release 8.5.0).

\subsection{Numerical results in 1D}
\label{Section1DNumerics}
Here, $V(\bx)$ is chosen as the harmonic plus lattice potential
\eqref{Lattice_Poten} with $\gm_x=1$, $k_x=25$ and $q_x=\fl{\pi}{2}.$
The computational domain and mesh size are respectively denoted as
$\mathcal{D}=[-L, L]$ and $h$.  In addition, to compare with the
common existing method BESP, we choose the stopping criteria
\eqref{Stop_Energy} with $\eps=10^{-14}$ all through this section.
For {\rm BESP}, we choose $\Delta t=0.01$ unless specified otherwise,
and fix the error tolerance for the inner loop to $10^{-10}$. Other
values of the error tolerance were also tried, but this choice was
found to be representative of the performance of BESP.

\begin{figure}[h!]
\centerline{
\includegraphics[height=4.7cm,width=5.2cm]{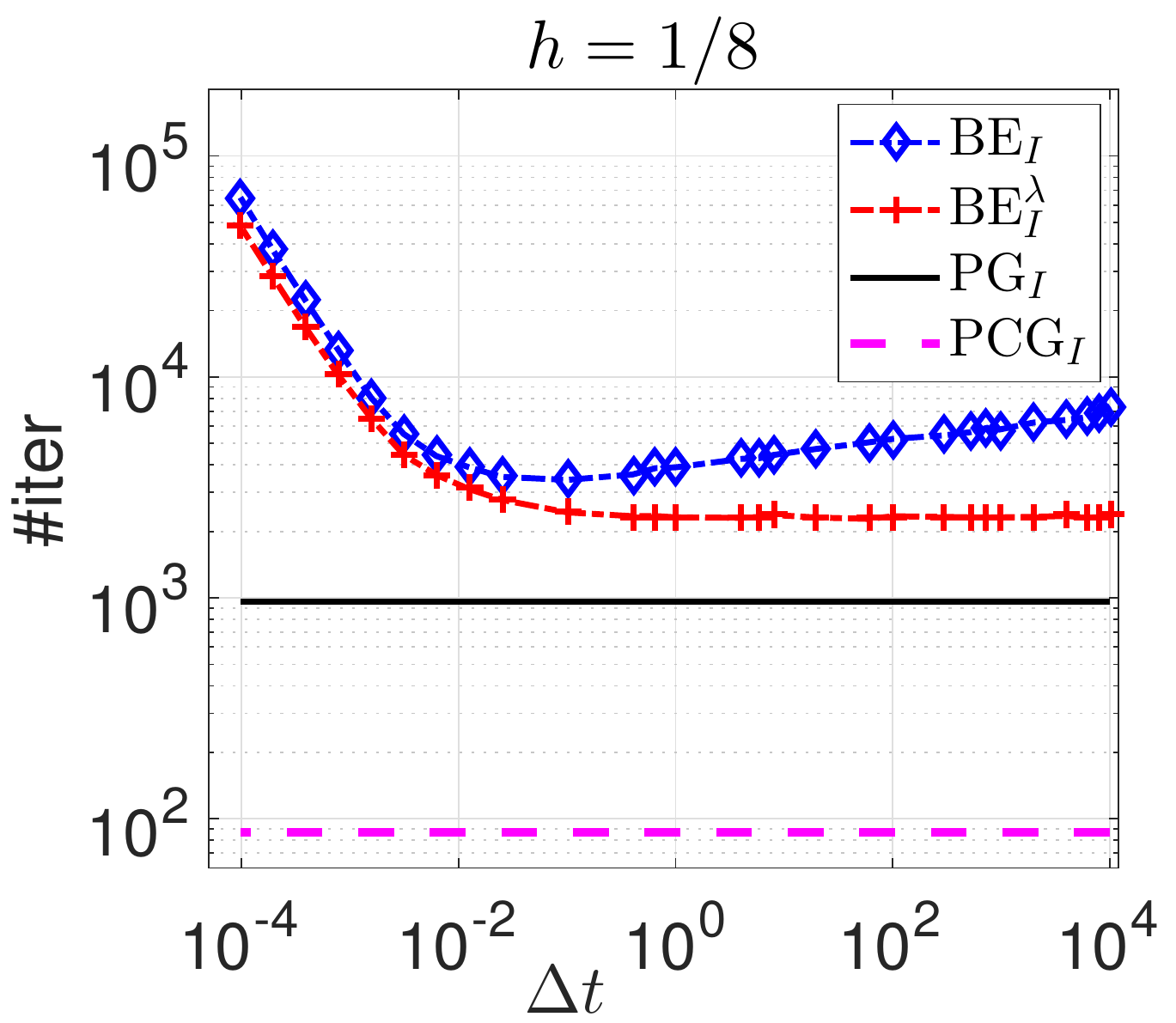}\quad
\includegraphics[height=4.7cm,width=5.2cm]{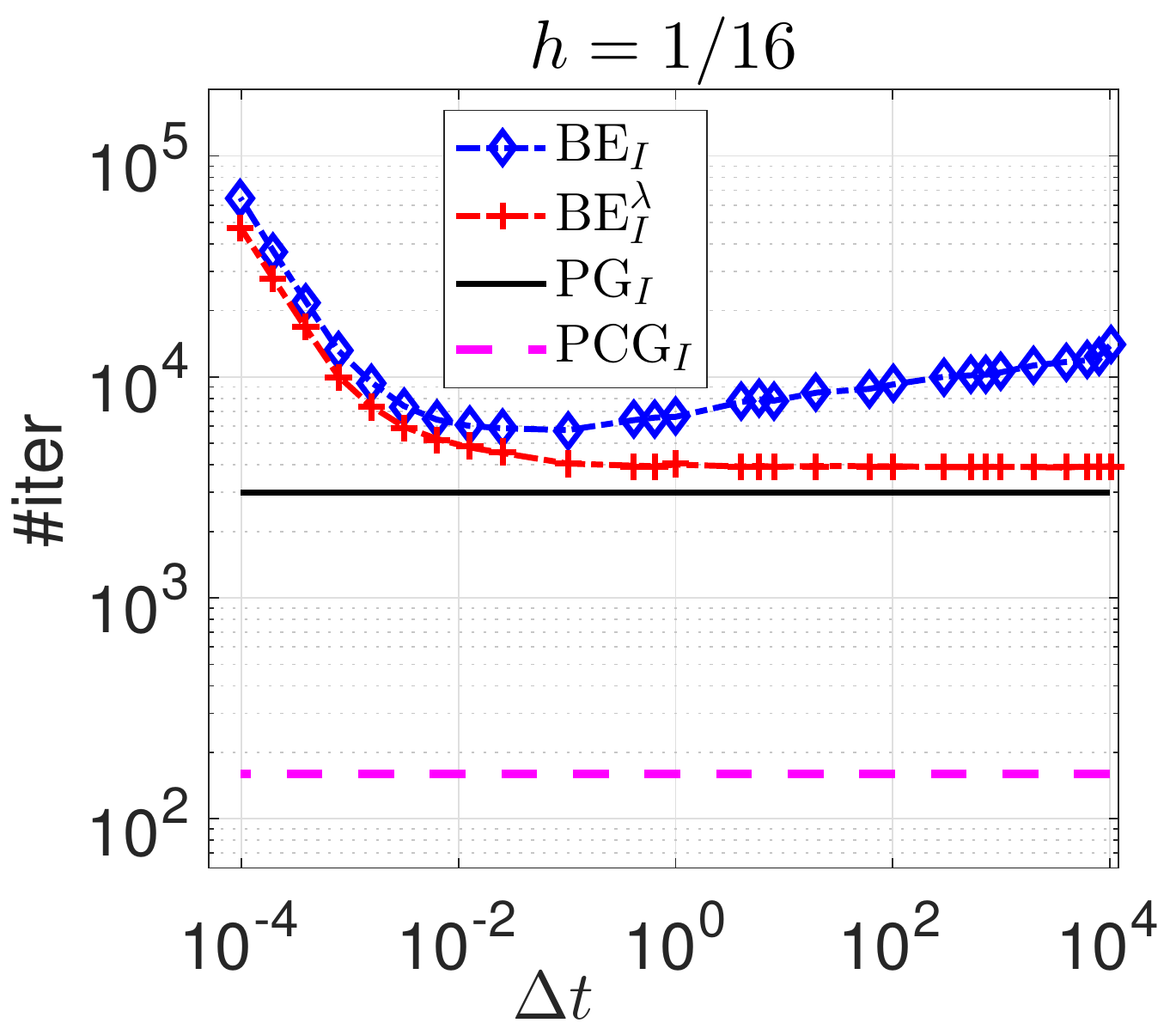}\quad
\includegraphics[height=4.7cm,width=5.2cm]{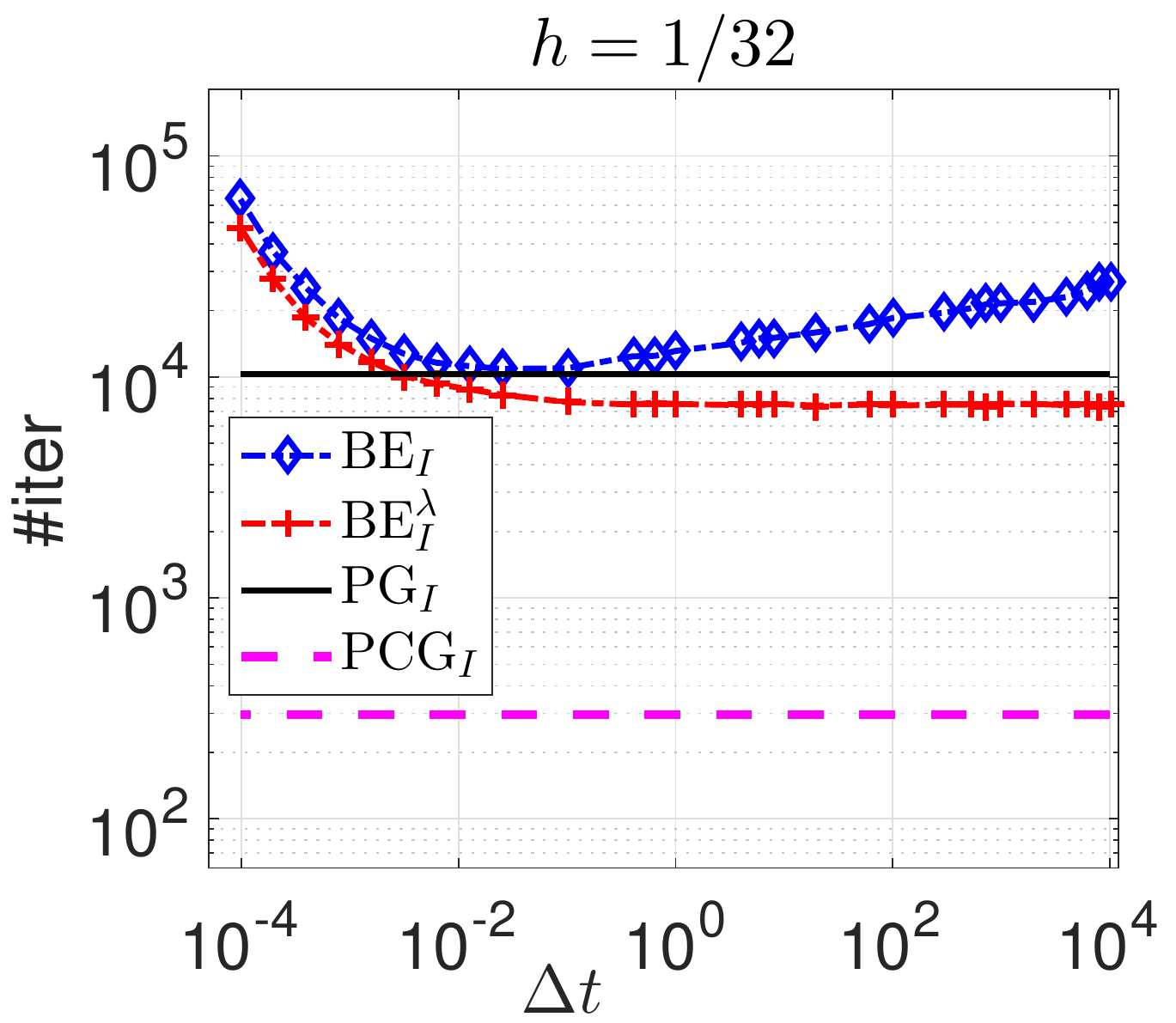}
}
\vspace{-0.2cm}
\caption{Example \ref{eg:1D_Compare_Solver}. Number of iterations to
  converge for different methods and different stepsizes, with
  different values of the discretization parameter $h$.}
\label{fig:1D_Compair_Solvers}
\end{figure}

\begin{exmp}
  \label{eg:1D_Compare_Solver}
  \rm We first compare the performance of various solvers without
  preconditioning in a simple case. We choose $L=16$ and $\eta=250$,
  and varying mesh sizes. We compare the method BE$_{I}$ given by
  \eqref{StdBEeq}, BE$_{I}^{\lambda}$ given \eqref{BEeq}, and the
  gradient and conjugate gradient algorithms in Figure
  \ref{fig:1D_Compair_Solvers}. The difference between methods
  BE$_{I}$ and BE$_{I}^{\lambda}$ is the inclusion of the chemical
  potential in the discretized gradient flow: we showed in Section
  \ref{ReviewSection} that both were equivalent for linear problems up
  to a renormalization in $\Delta t$. We see here that this conclusion
  approximately holds even in the nonlinear regime ($\eta \neq 0$),
  with both methods performing very similarly until $\Delta t$ becomes
  large, at which point the BE$_{I}^{\lambda}$ effectively uses a
  constant stepsize (see \eqref{eq:eq_BE_BElambda}), while the large
  timestep in BE$_{I}$ makes the method inefficient. In this case,
  $\lambda_{1}$ is positive, so that both methods converge to the
  ground state even for a very large $\Delta t$. Overall we see that
  the optimum number of iterations is achieved for a value of
  $\Delta t$ of about $0.01$, which we keep in the following tests to
  ensure a fair comparison. We also use the BE$^{\lambda}$ variant in
  the following tests.

  For modest values of the discretization parameter $h$, the Backward
  Euler methods are less efficient than the gradient method (which can
  be interpreted as a Forward Euler iteration with adaptive
  stepsize). As $h$ is decreased, the conditioning of the problem
  increases as $h^{-2}$. The gradient/Forward Euler method is limited
  by its CFL condition, and its number of iterations grows like
  $h^{-2}$, as can readily be checked in Figure
  \ref{fig:1D_Compair_Solvers}.  The Backward Euler methods, however,
  use an efficient Krylov solver that is only sensitive to the square
  root of the conditioning, and its number of iterations grows only
  like $h^{-1}$. Therefore it become more efficient than the
  gradient/Forward Euler method.

  The conjugate gradient method is always more efficient than the
  other methods by factors varying between one and two orders of
  magnitude. Its efficiency can be attributed to the combination of
  Krylov-like properties (as the Backward Euler method, its iteration
  count displays only a $h^{-1}$ growth) and optimal
  stepsizes.
\end{exmp}

\begin{figure}[h!]
\centerline{
\includegraphics[height=5.cm,width=7.5cm]{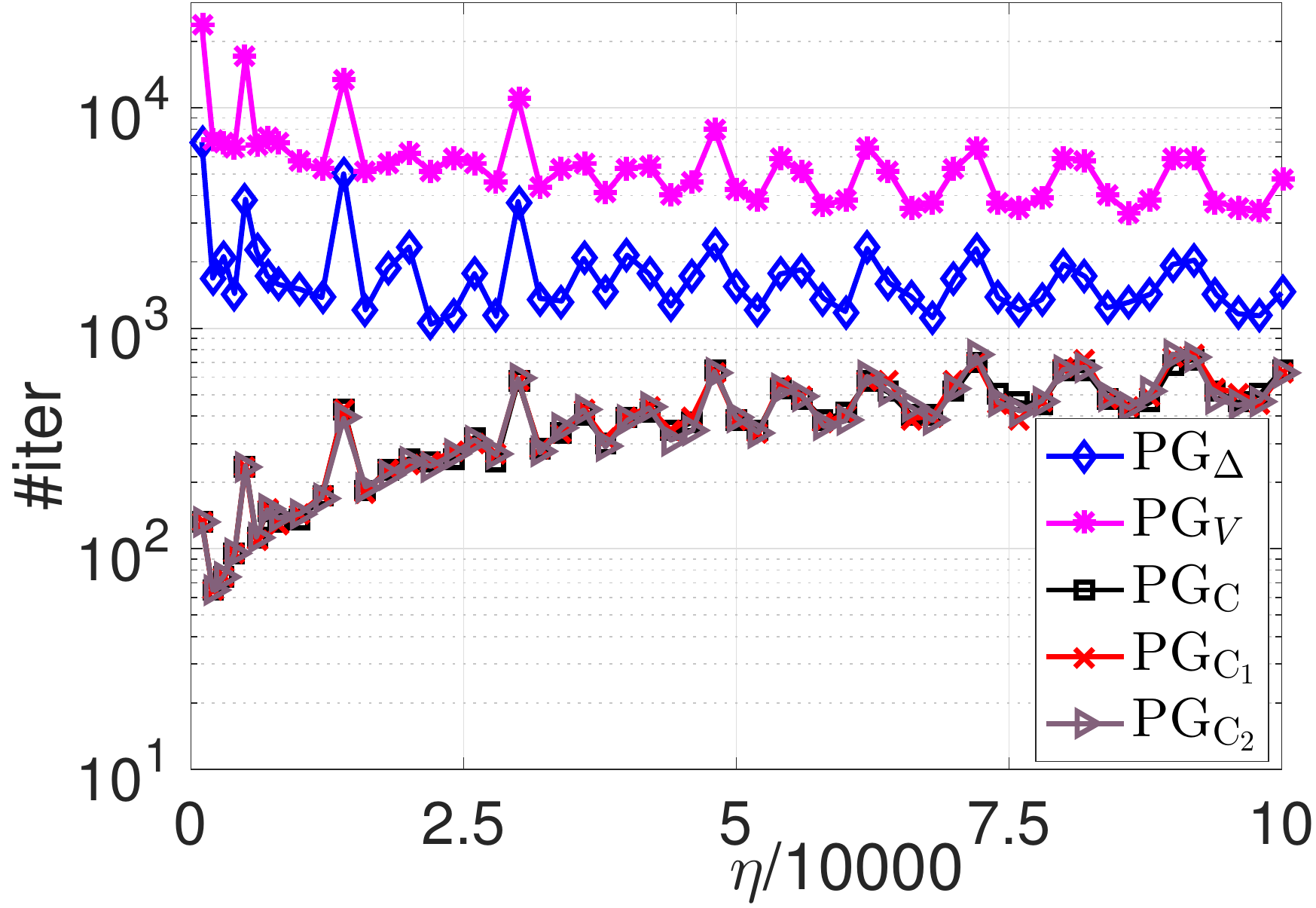}\qquad
\includegraphics[height=5.cm,width=7.5cm]{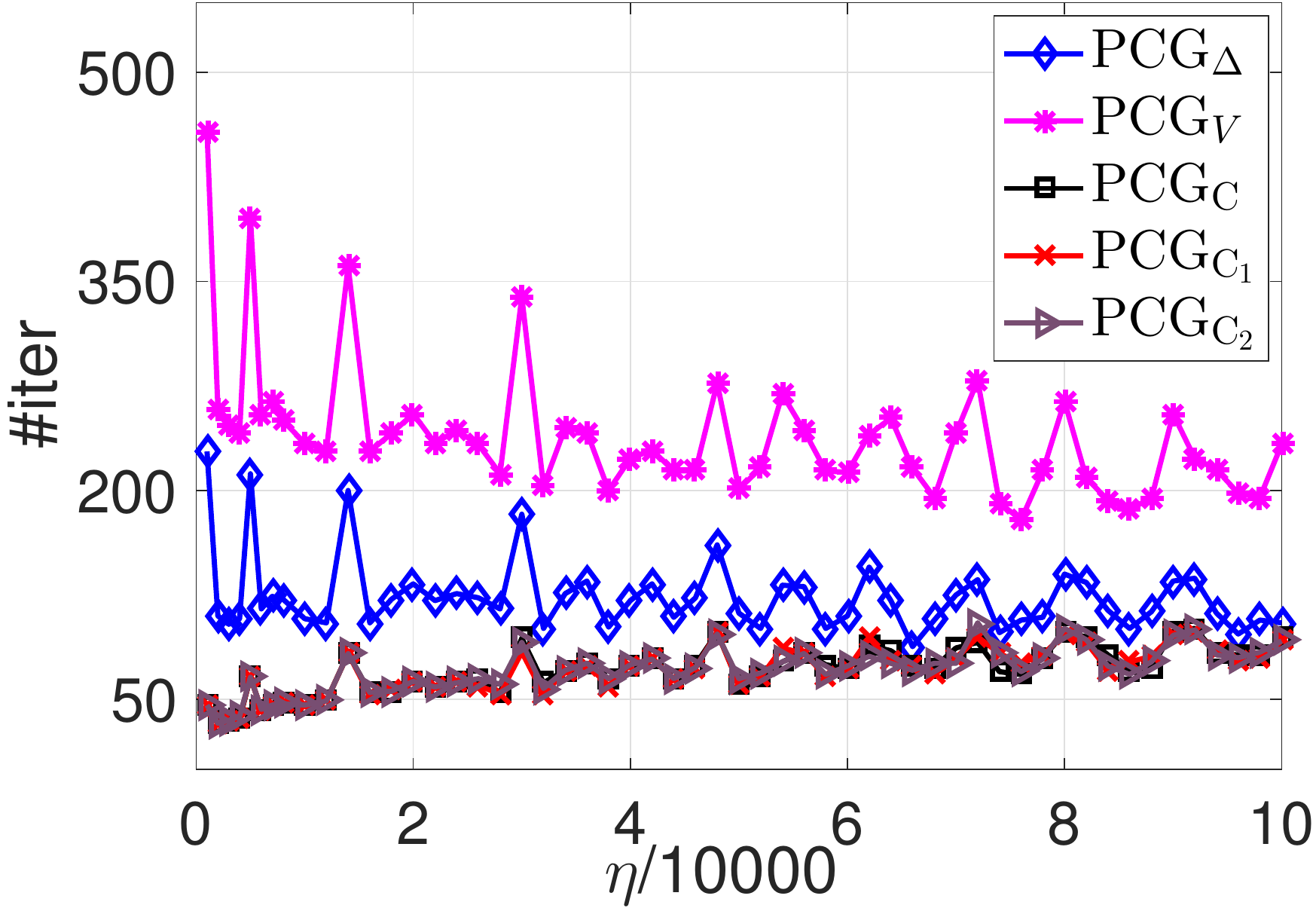}
}
\vspace{-0.2cm}
\caption{Example \ref{eg:1D_PCG_Diff_Precon}. Number of iterations of PG$_\nu$  and
 PCG$_\nu$ ($\nu=\Delta,V, {\rm C}, {\rm C_1},  {\rm C_2}$) to converge,
  for different nonlinear strengths $\eta$.}
\label{fig:1D_Com_PCG_Diff_LaplacianPre}
\end{figure}

\begin{exmp}
  \label{eg:1D_PCG_Diff_Precon} {\rm We  compare now the performance
    of the (conjugate) gradient method with different preconditioners.
    To this end, we consider the algorithms PG$_{\nu}$ and PCG$_{\nu}$
    with $\nu=\Delta,V,{\rm C}, {\rm C_1}, {\rm C_2}.$ The
    computational parameters are chosen as $L=128$ and $h=\fl{1}{64}$,
    respectively.  Fig. \ref{fig:1D_Com_PCG_Diff_LaplacianPre} shows
    the iteration number $\# \rm iter$ for these schemes and different
    values of the nonlinearity strength $\eta$.  From this figure and
    other numerical results not shown here, we can see that: (i) For
    each fixed preconditioner, the PCG schemes works better than the
    PG schemes; (ii) the combined preconditioners all work equally
    well, and bring a reduction in the number of iteration, at the
    price of more Fourier transforms.
  }

\end{exmp}

\begin{figure}[h!]
  \centering
  \begin{tabular}[h!]{c|m{4cm}m{4cm}m{4cm}}
    &\hspace{2cm}BE&\hspace{2cm}PG&\hspace{2cm}PCG\\
    \hline\\
    $P_{I}$ &
              \includegraphics[height=3.6cm,width=4cm,angle=0]{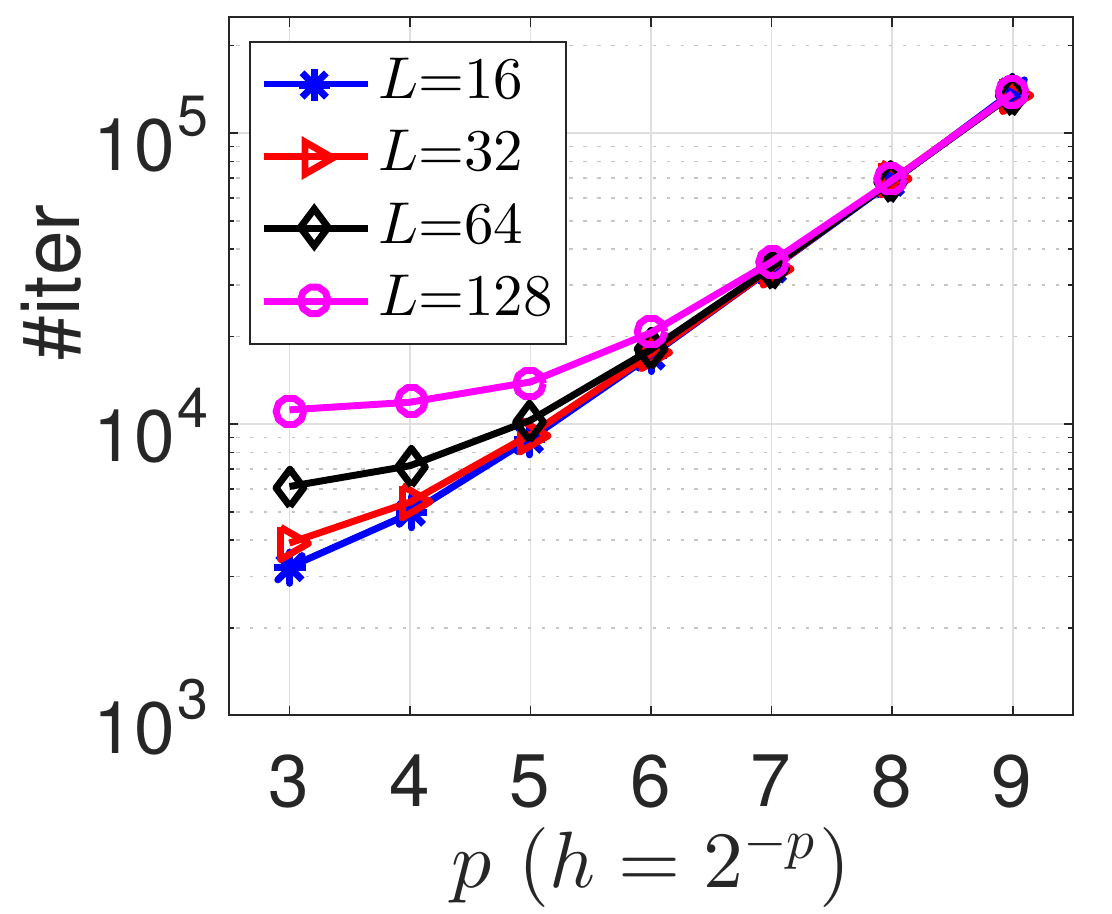}
                   &
                     \includegraphics[height=3.6cm,width=4cm]{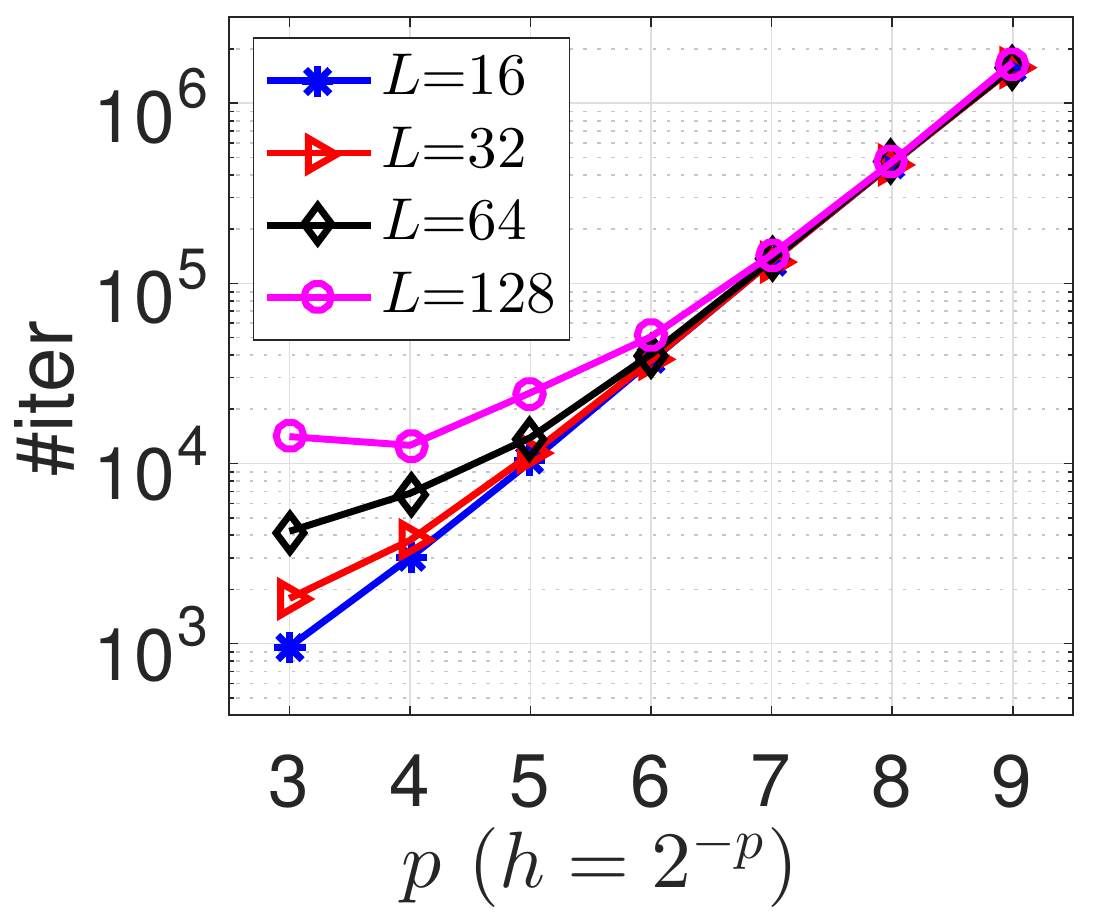}&
\includegraphics[height=3.6cm,width=4cm]{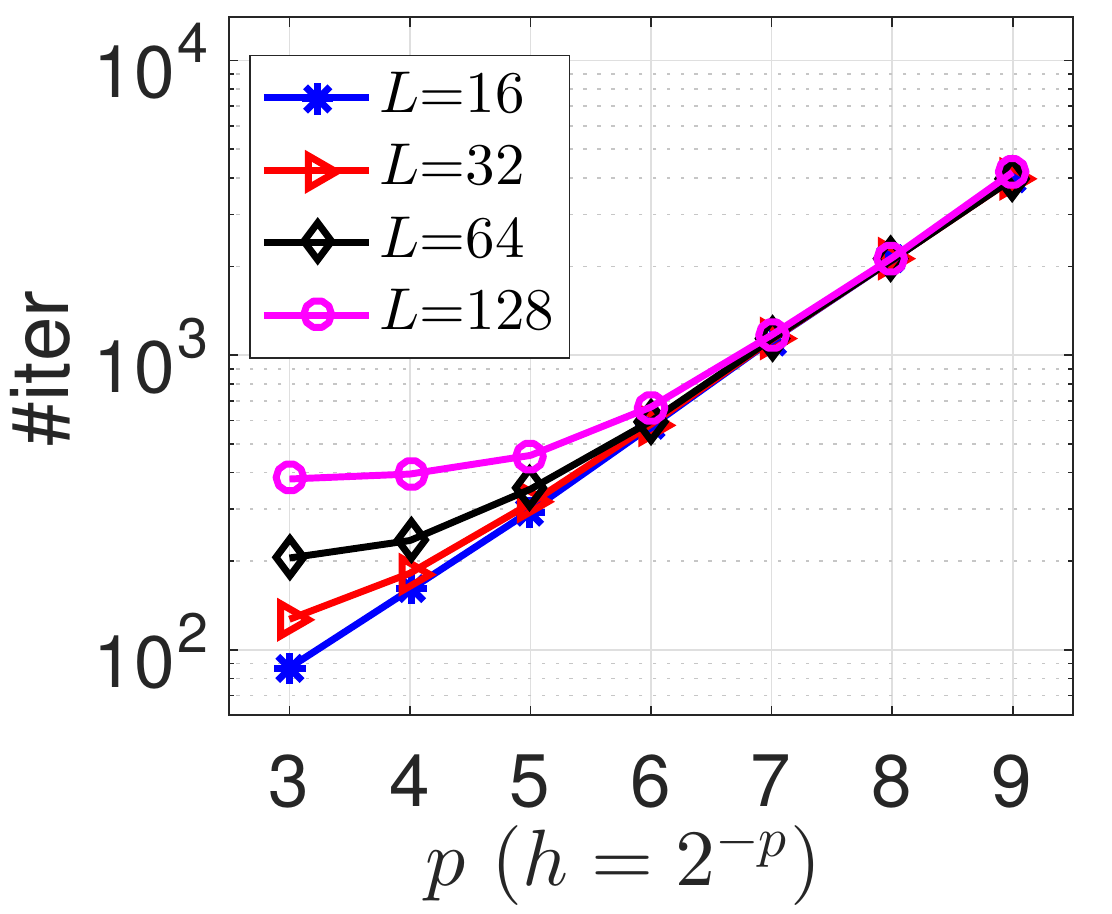}\\
$P_{V}$&\includegraphics[height=3.6cm,width=4cm]{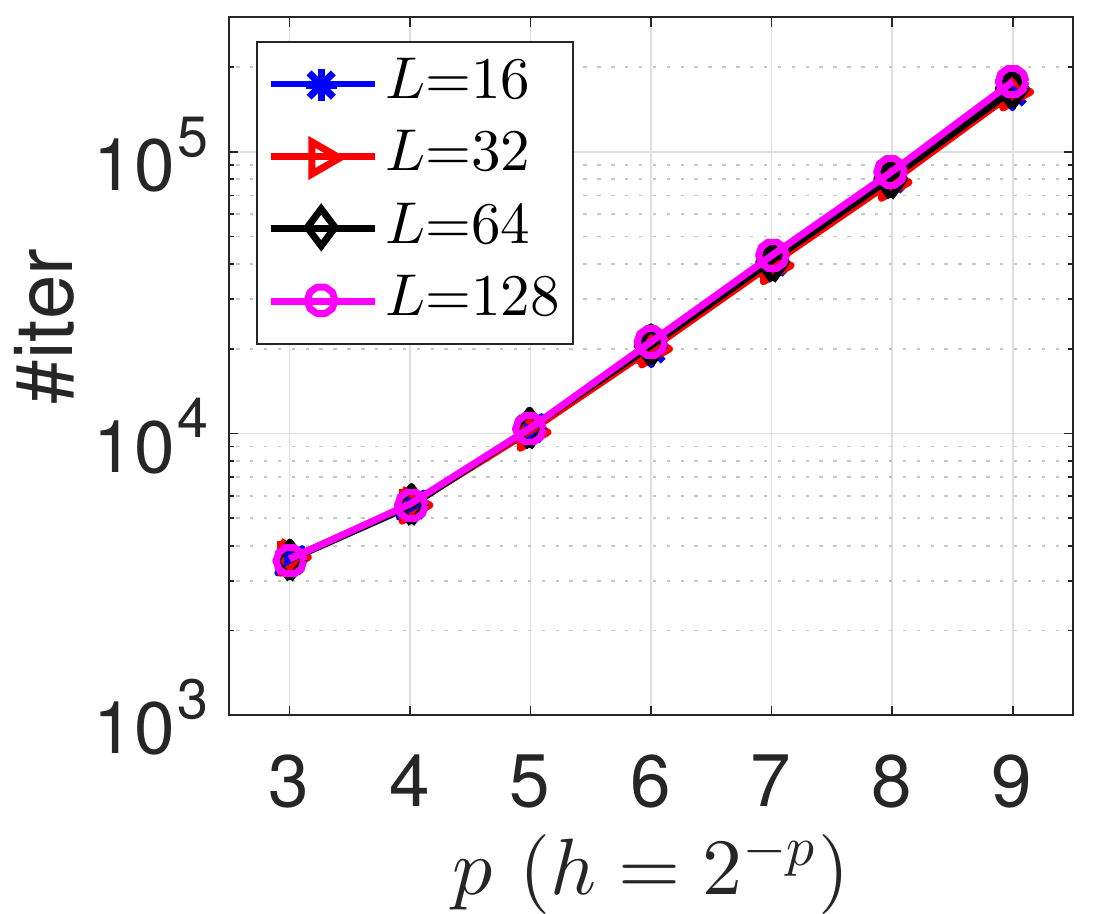}&
\includegraphics[height=3.6cm,width=4cm]{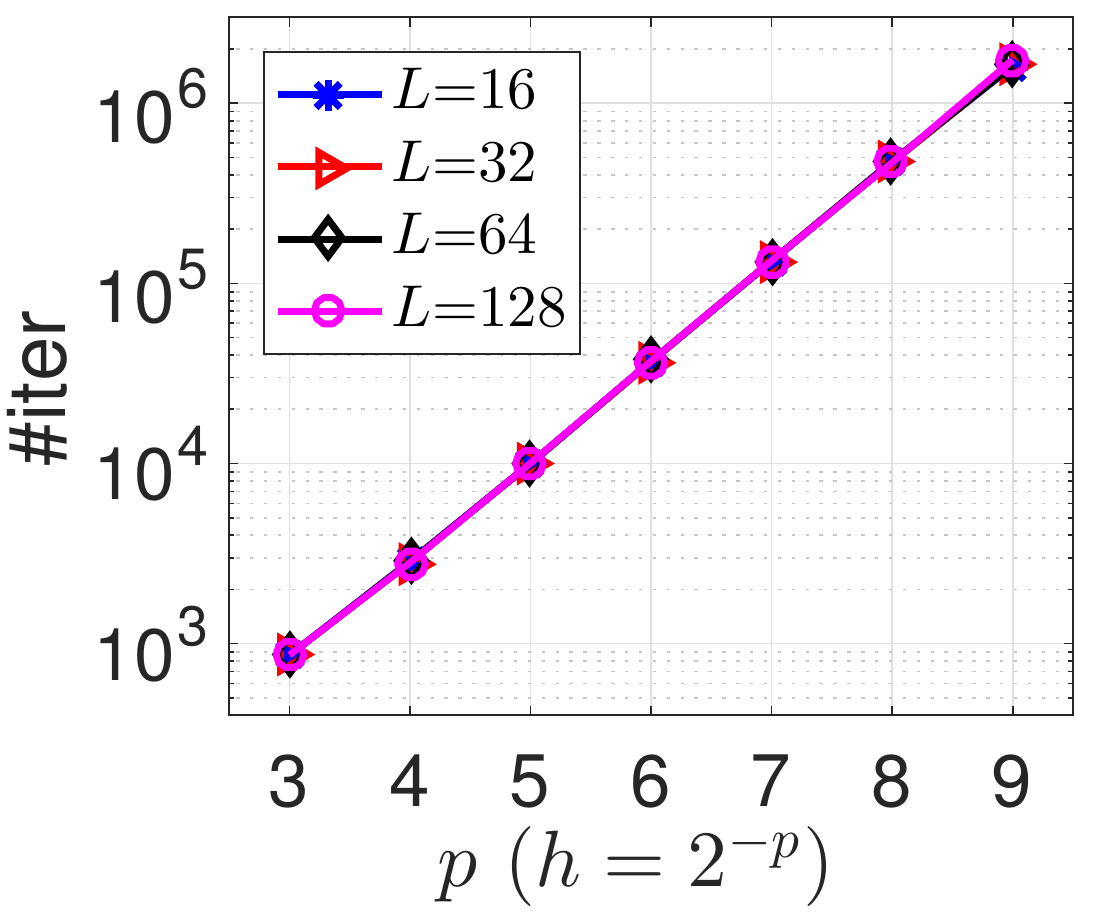}&
\includegraphics[height=3.6cm,width=4cm]{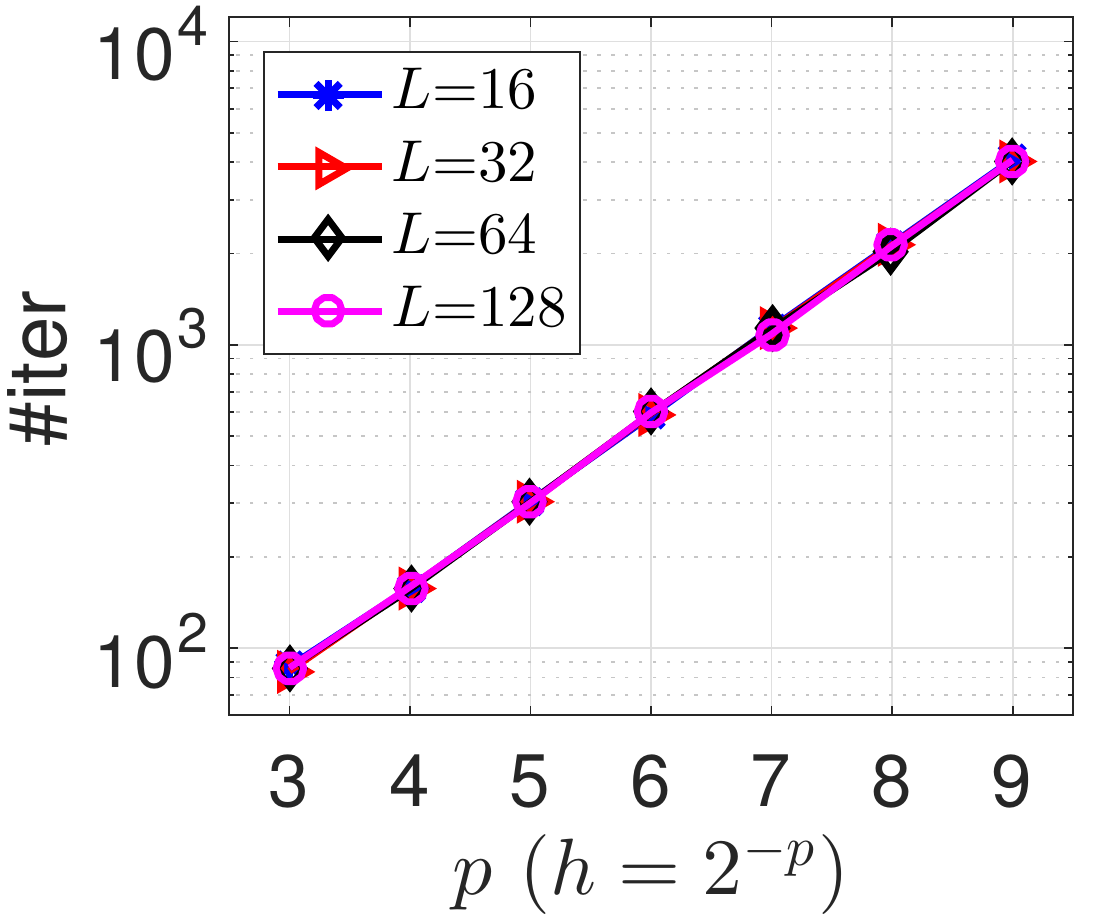}\\
$P_{\Delta}$&\includegraphics[height=3.6cm,width=4cm]{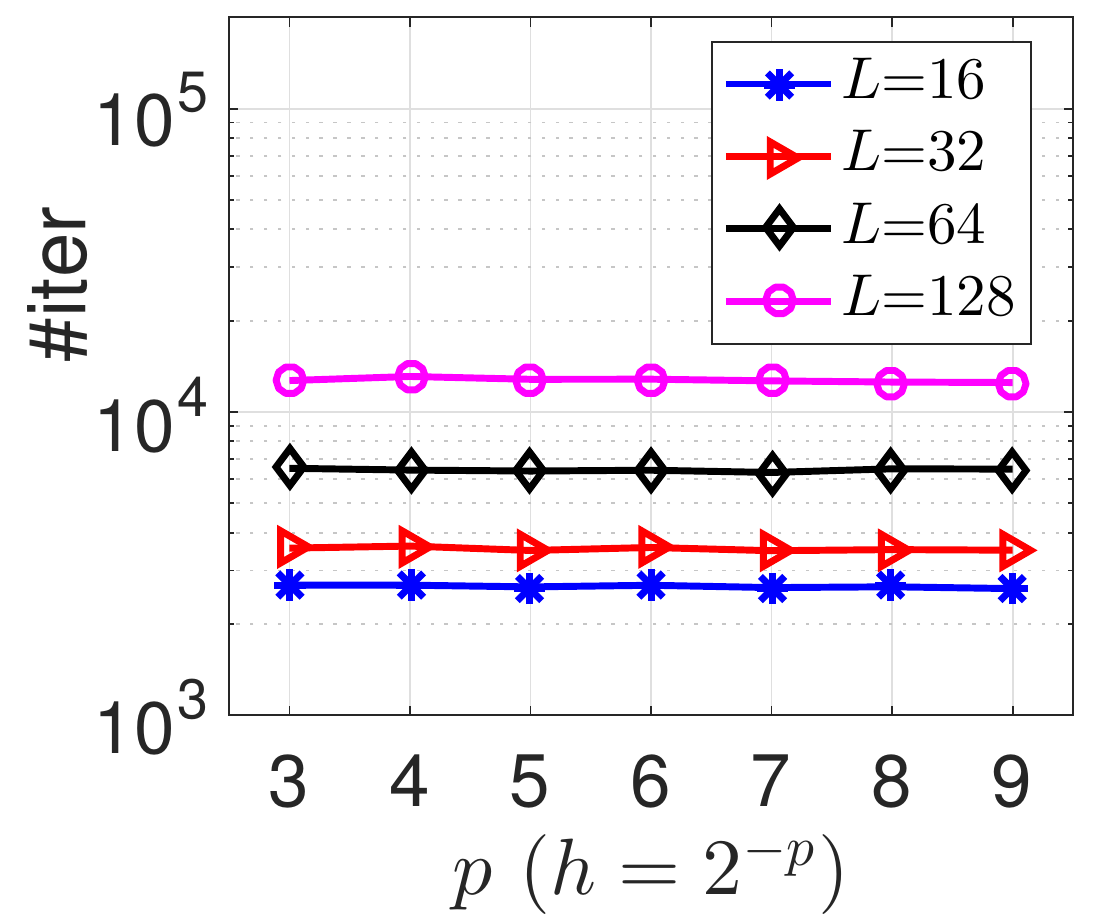}&
\includegraphics[height=3.6cm,width=4cm]{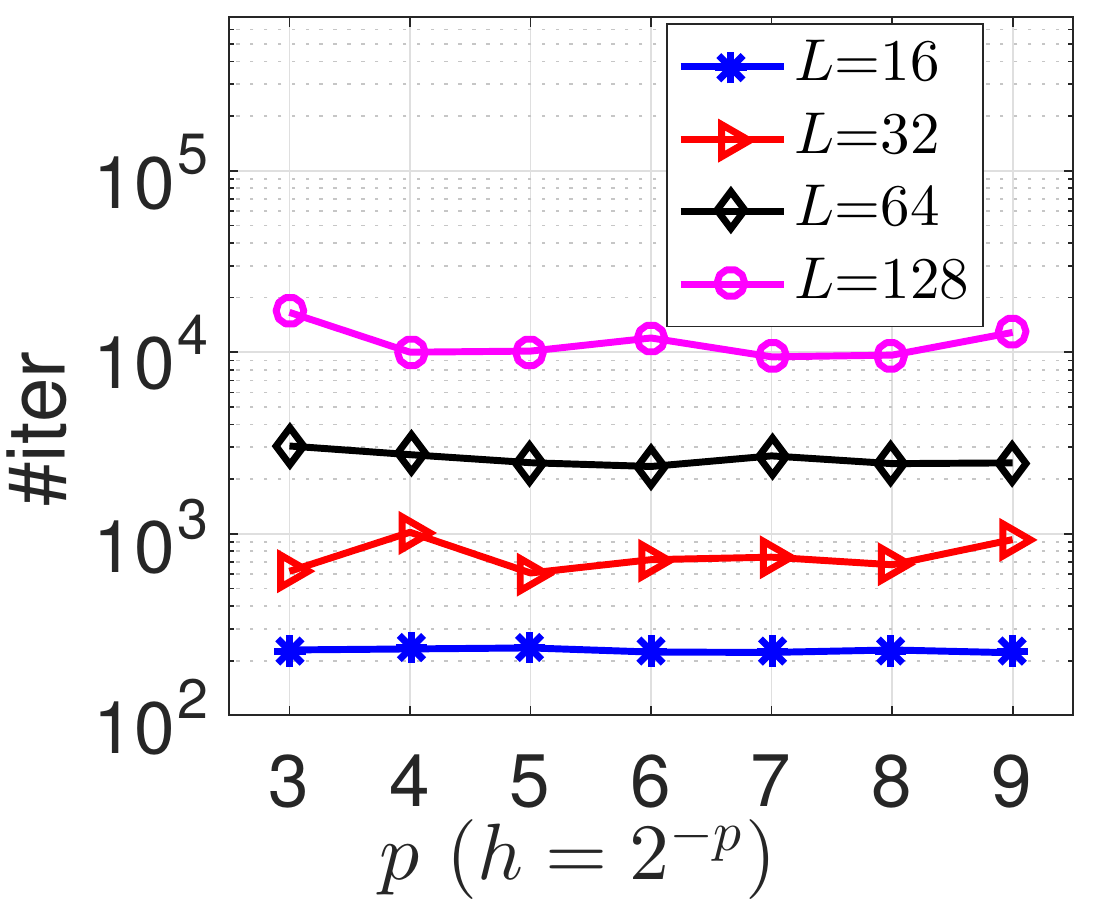}&
\includegraphics[height=3.6cm,width=4cm]{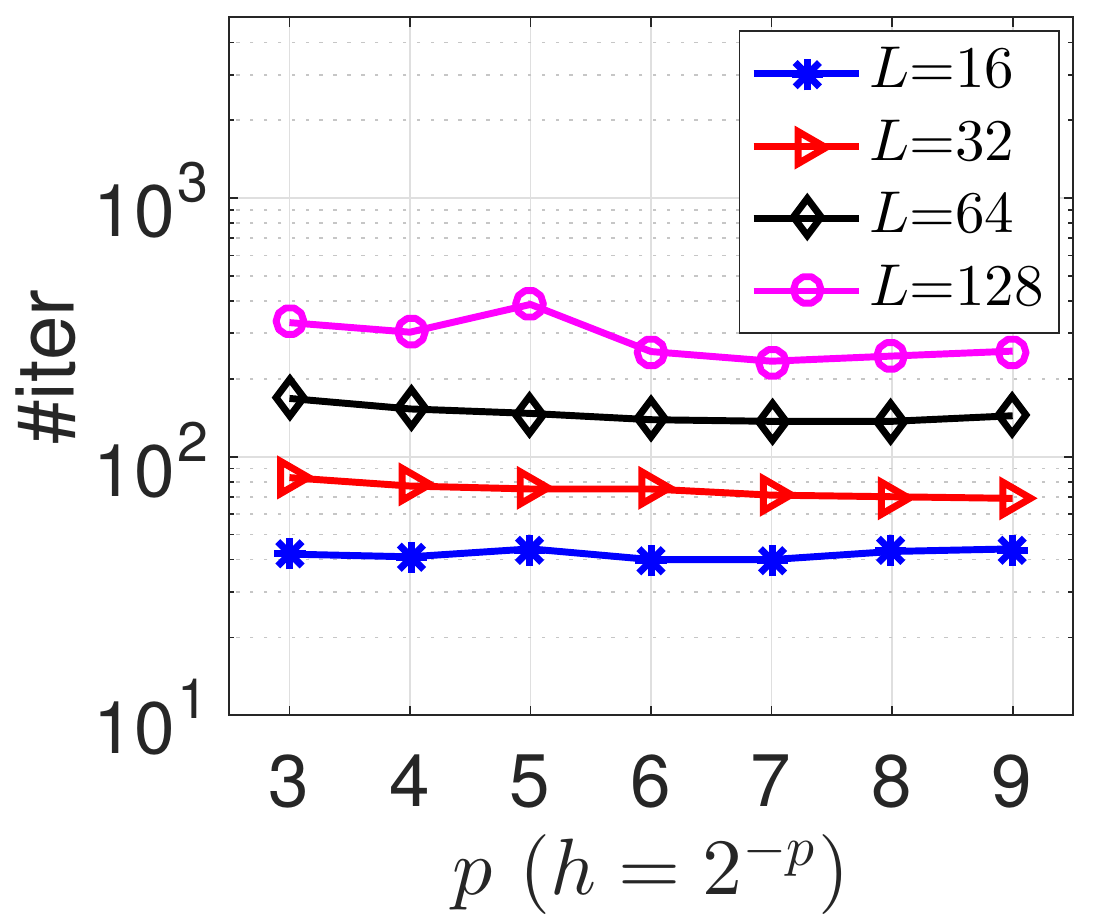}\\
$P_{\rm C}$&\includegraphics[height=3.6cm,width=4cm]{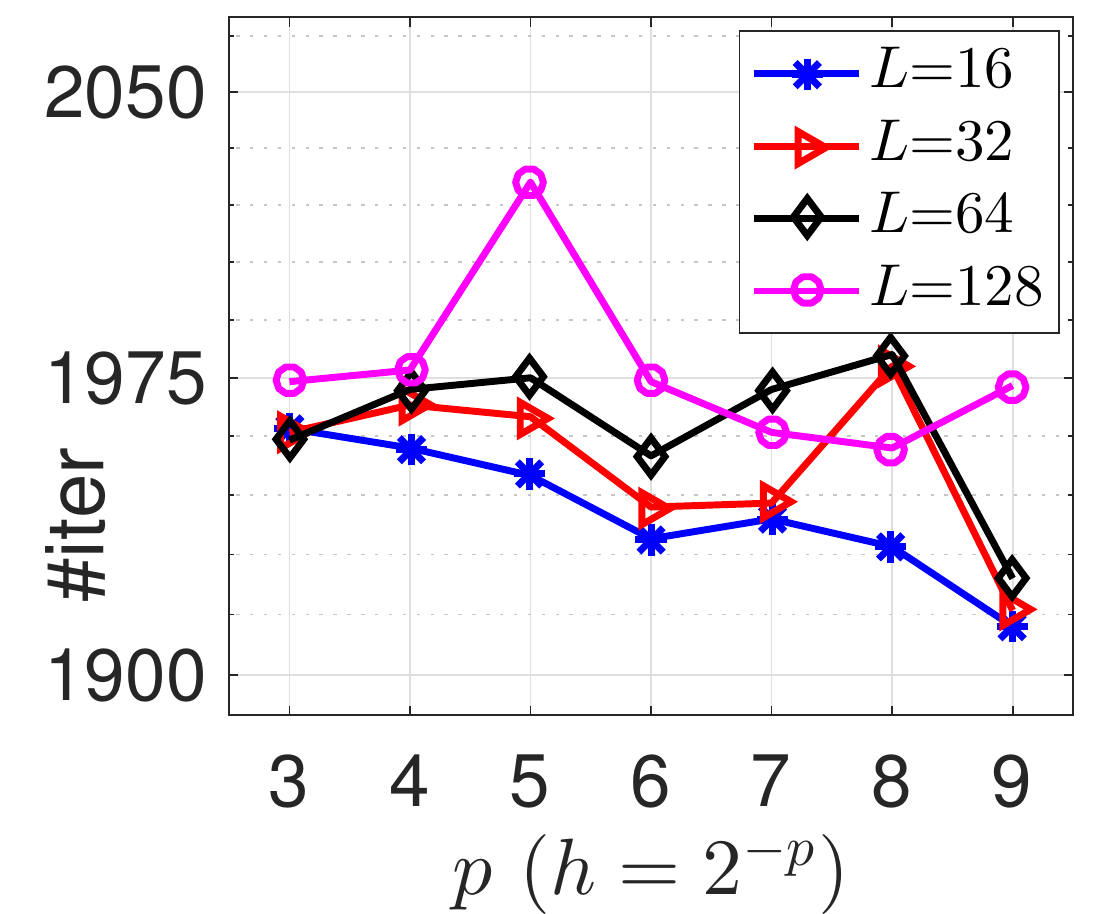}&
\includegraphics[height=3.6cm,width=4cm]{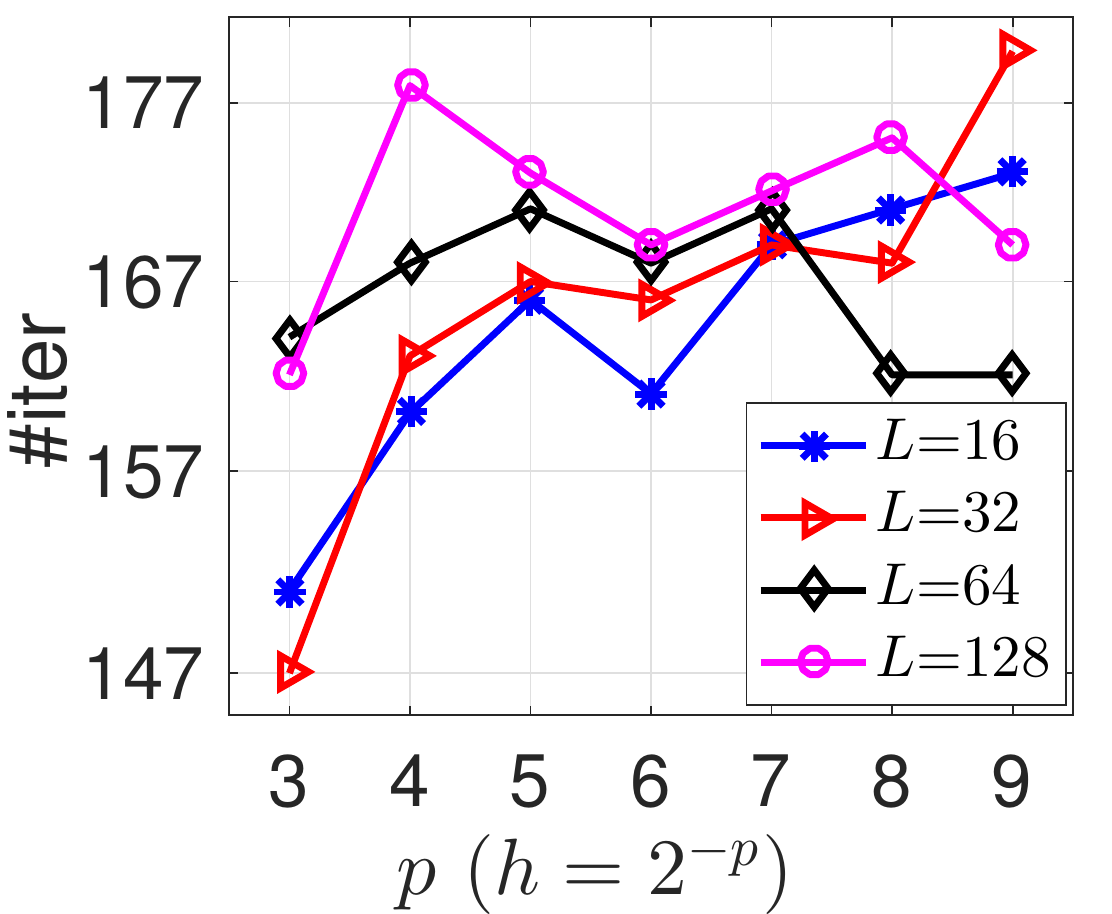}&
\includegraphics[height=3.6cm,width=4cm]{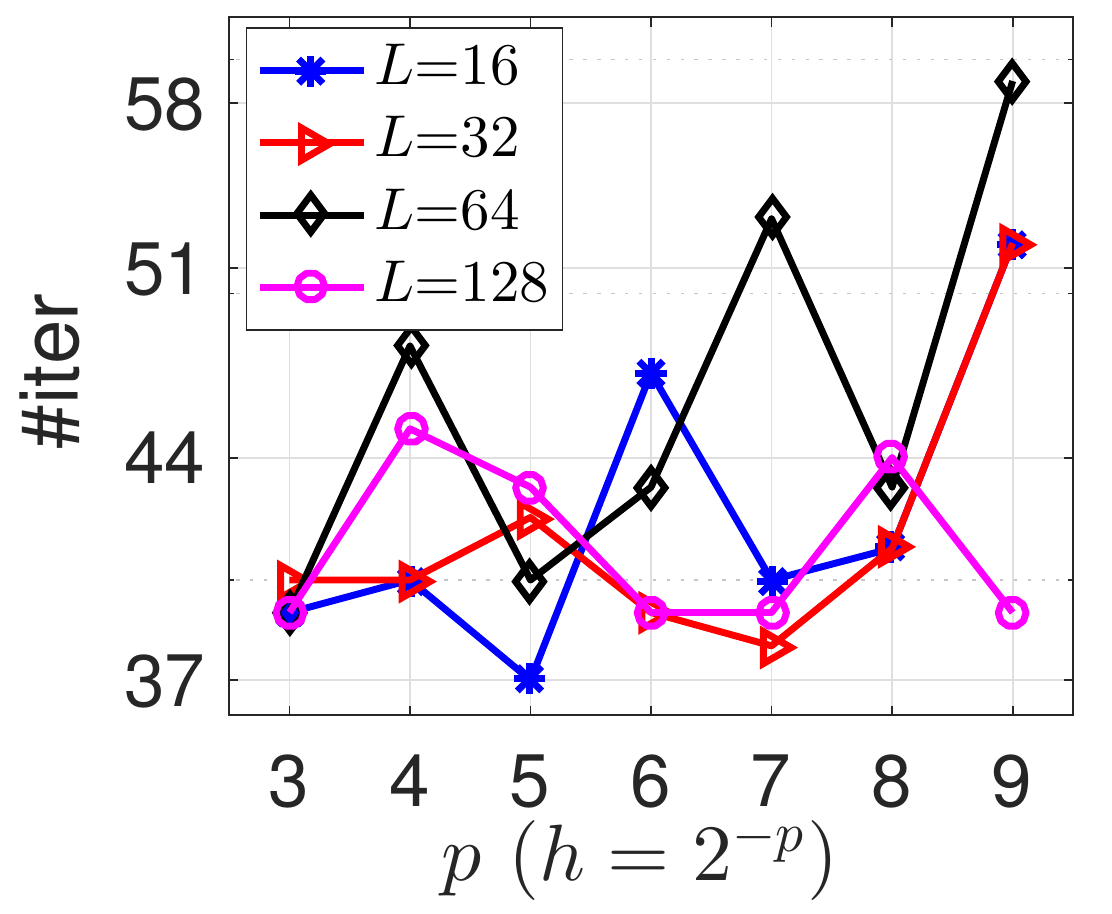}\\
  \end{tabular}
  \caption{Example \ref{eg:1DComIters}. Number of iterations to converge for {\rm BE}$_{\nu}$, {\rm PG}$_{\nu}$ and {\rm PCG}$_{\nu}$
for $\nu=I, V, \Delta, \textrm{C}$, vs. the mesh
refinement $h$.}
\label{fig:IterNum_1D_Compare_Diff_Meth}
\end{figure}

\begin{exmp}
  \label{eg:1DComIters} {\rm In this example, we compare the
    performance of {\rm PG}$_{\nu}$, {\rm PCG}$_{\nu}$ and BE$_{\nu}$
    ($\nu=I, \Delta,V,\textrm{C}$) with respect to different domain
    and mesh sizes.  To this end, we fix $\eta=250$.
    Fig. \ref{fig:IterNum_1D_Compare_Diff_Meth} shows the total
    iteration number for these schemes with different $L$ and $h$.
    From this figure and additional numerical results not shown here
    for brevity, we see that: (i) Preconditioned solvers outperform
    unpreconditioned solvers; (ii) The potential preconditioner $P_V$
    \eqref{PV} makes the solver mainly depend on the spatial
    resolution $h$, while the kinetic potential preconditioner
    $P_\Delta$ \eqref{PDelta} prevents the deterioration as $h$
    decreases for a fixed $L$,
    consistent with the theoretical analysis in subsection
    \ref{sectionPreconditioners}; (iii) The deterioration is less marked for
    Krylov-based methods (BE and PCG) than for the PG method, because Krylov methods
    only depends on the square root of the condition number (iv) The combined preconditioner
    $P_{\rm C}$ \eqref{combinedP} makes the solvers almost independent
    of both the parameters $h$ and $L$, although we theoretically
    expect a stronger dependence. We attribute this to the fact that
    we start with a specific initial guess that does not excite the
    slowly convergent modes enough to see the dependence on $h$ and
    $L$; (v) For each solver, the combined preconditioner $P_{\rm C}$
    performs best.  Usually, PCG$_{\rm C}$ is the most efficient
    algorithm, followed by PG$_{\rm C}$, and finally BE$_{\rm C}$.  }

\end{exmp}

\begin{figure}[h!]
\centerline{
\includegraphics[height=4.5cm,width=8.1cm]{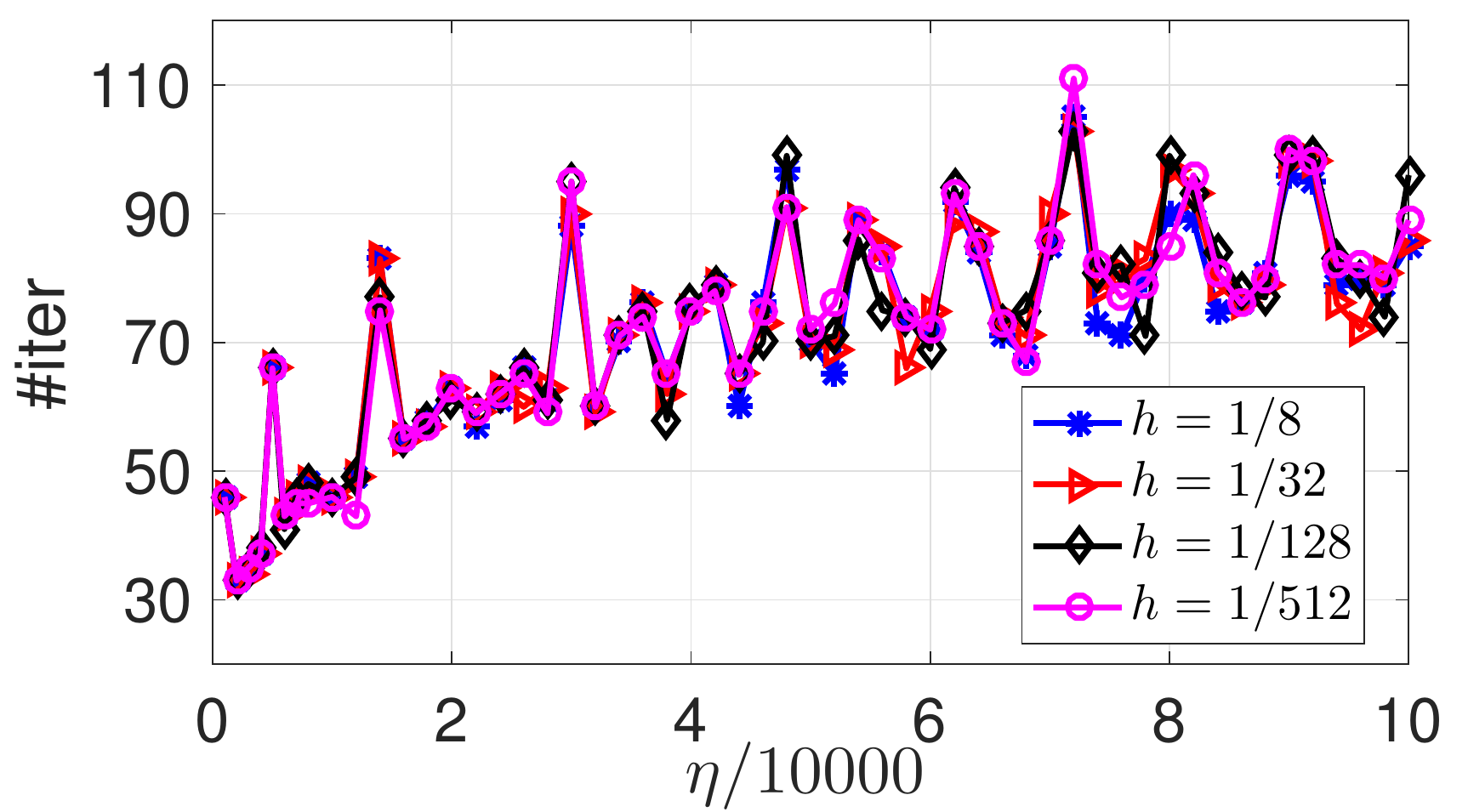}\;
\includegraphics[height=4.5cm,width=8.1cm]{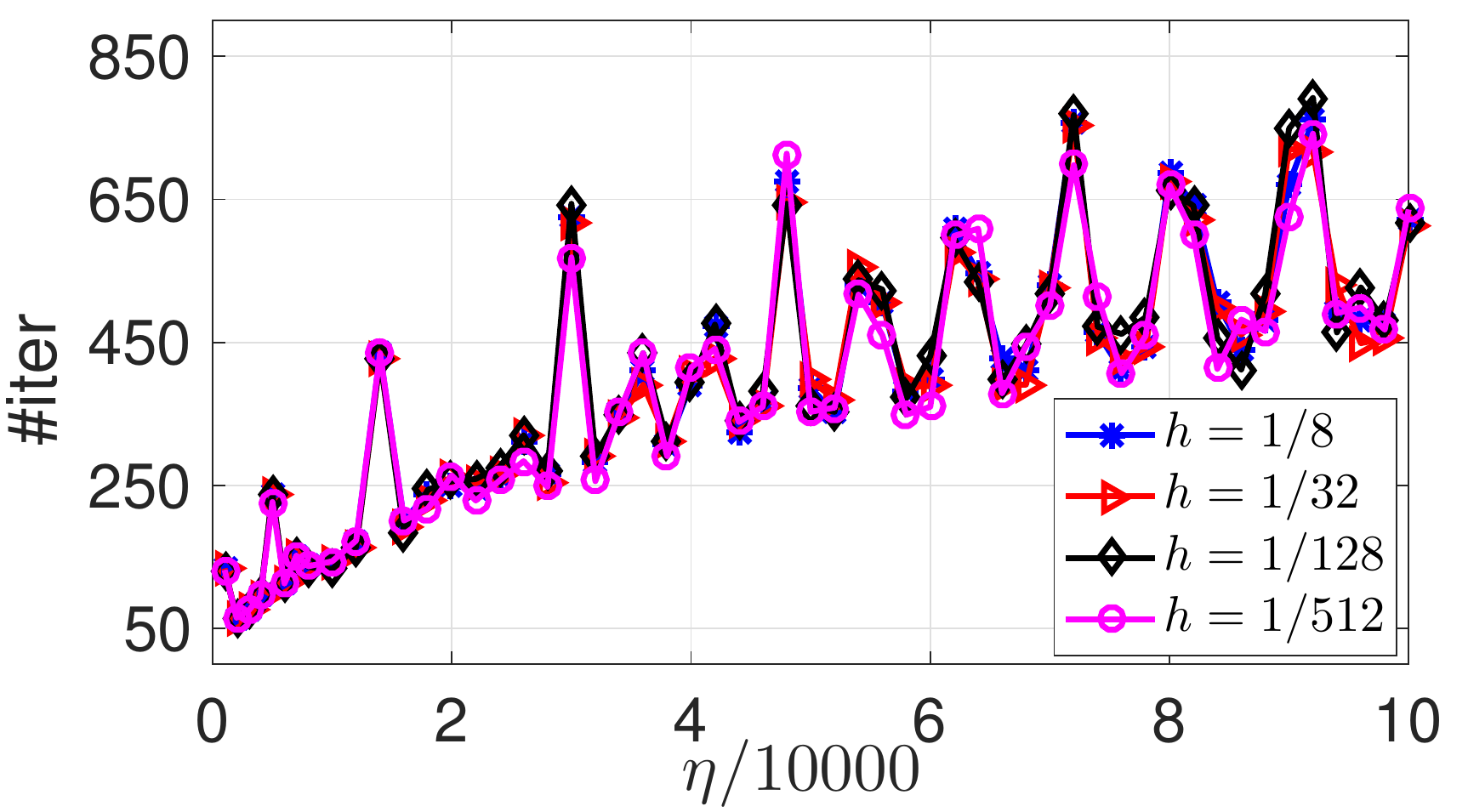}
}
\vspace{-0.2cm}
\caption{Example \ref{eg:1DComp_PCG_PG_SymPrecon}. Number of
  iterations to converge for  PCG$_{\rm C}$ (left) and
  PG$_{\rm C}$ (right) with $L=128$,  and various values of $h$ and
  $\eta$.}
\label{fig:1DComp_PCG_PG_SymPrecon}
\end{figure}

\begin{exmp}
\label{eg:1DComp_PCG_PG_SymPrecon}
{\rm 
We investigate further  the performance of 
PCG$_{\rm C}$ and PG$_{\rm C}$  with different nonlinear interaction strenghts $\eta$.
To this end, we take $L=128$ and different discretization parameters $h$. We vary the nonlinearity from $\eta=0$
to $\eta=10^5.$ Fig. \ref{fig:1DComp_PCG_PG_SymPrecon} depicts the corresponding iteration numbers to converge. 
We could clearly see that:
(i) The iteration counts for both methods are almost independent on $h$, but both depend on
the nonlinearity $\eta$; PCG$_{\rm C}$  depends slightly on $\eta$
while PG$_{\rm C}$ is more sensitive; 
(ii) For fixed $\eta$ and $h$,  PCG$_{\rm C}$ converges much faster than PG$_{\rm C}$.
}

\end{exmp}

From examples \ref{eg:1D_PCG_Diff_Precon}--\ref{eg:1DComp_PCG_PG_SymPrecon}, we see that PCG$_{\rm  C}$, i.e. 
PCG with combined symmetric preconditioner $P_{\rm C}$ \eqref{combinedP}
is the best solver. Hereafter, unless stated, we use PCG$_{\rm C}$ 
as the default solver to compute the ground states.

\subsection{Numerical results in 2D}\label{Section2DNumerics}
Here,  we choose $V(\bx)$ as the harmonic plus quartic potential \eqref{Quart_Poten} with $\gm_x=\gm_y=1$, $\alpha=1.2$
and $\kappa=0.3.$     The computational domain and mesh sizes are chosen respectively as  $\mathcal{D}=[-16, 16]^2$ and $h=\fl{1}{16}.$ 

\begin{figure}[h!]
\centerline{
\includegraphics[height=4.3cm,width=10cm]{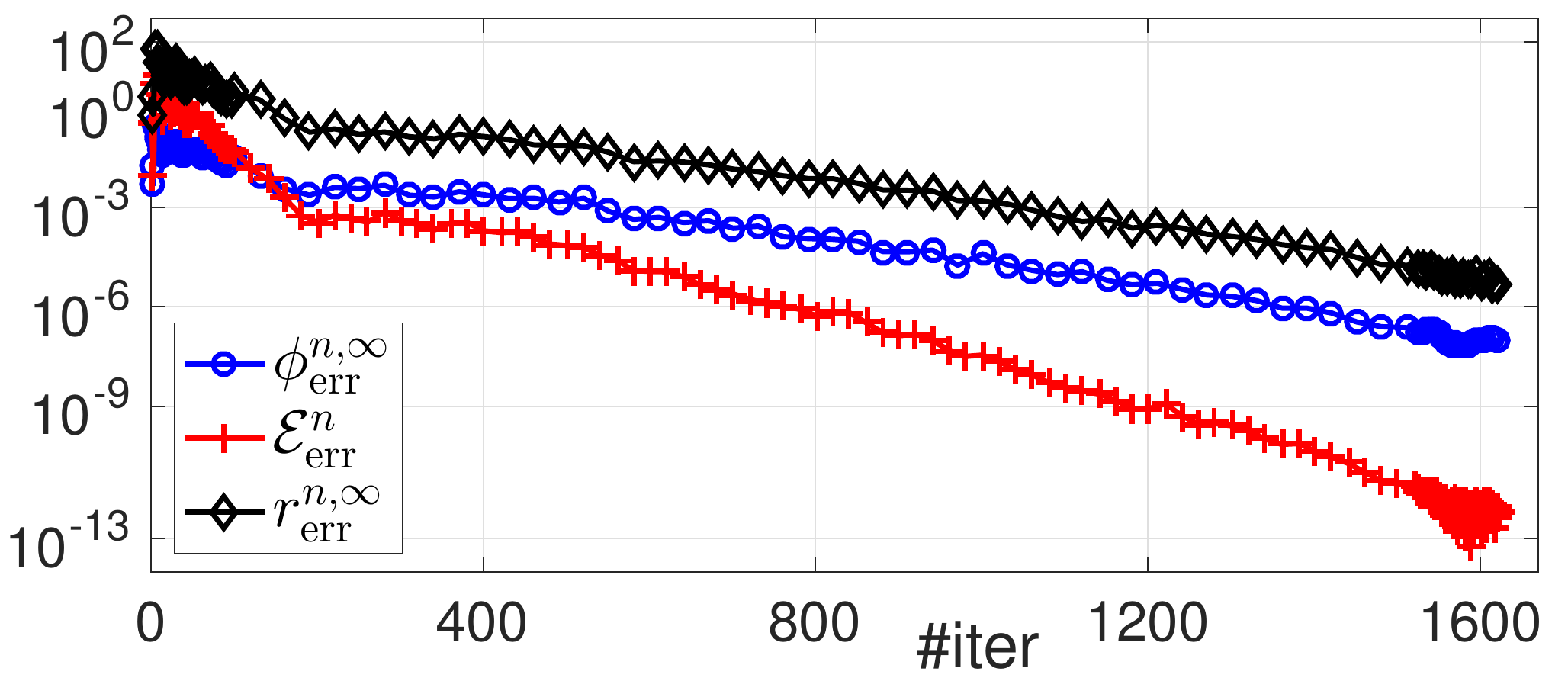}
}
\vspace{-0.2cm}
\caption{Evolution of the errors vs. the  total  number  of iterations.}
\label{fig:2D_stopping_ctiteria}
\end{figure}

First, we test the evolution of the three errors
\eqref{Stop_Max}--\eqref{Stop_Energy} as the algorithm progresses.  To
this end, we take $\eta=1000$ and $\og=3.5$ as
example. Fig. \ref{fig:2D_stopping_ctiteria} plots the
$\phi^{n,\infty}_{\rm err}$, $r^{n,\infty}_{\rm err}$ and
$\mathcal{E}^{n}_{\rm err}$ errors with respect to the iteration
number.  We can see clearly that
$\mathcal{E}^{n}_{\rm err}$ converges faster the other two indicators,
as expected.  Considering $\phi^{n,\infty}_{\rm err}$ or
$r^{n,\infty}_{\rm err}$ with an improper but relative large tolerance
would require a very long computational time to converge even if the
energy would not change so much. This is most particularly true for
large values of $\og$.
In all the examples below, unless stated, we fix  $\mathcal{E}^{n}_{\rm err}$  \eqref{Stop_Energy}  with  $\eps=10^{-12}$ to terminate the code. 
   
\begin{figure}[h!]
\centerline{
\includegraphics[height=5.2cm,width=8cm]{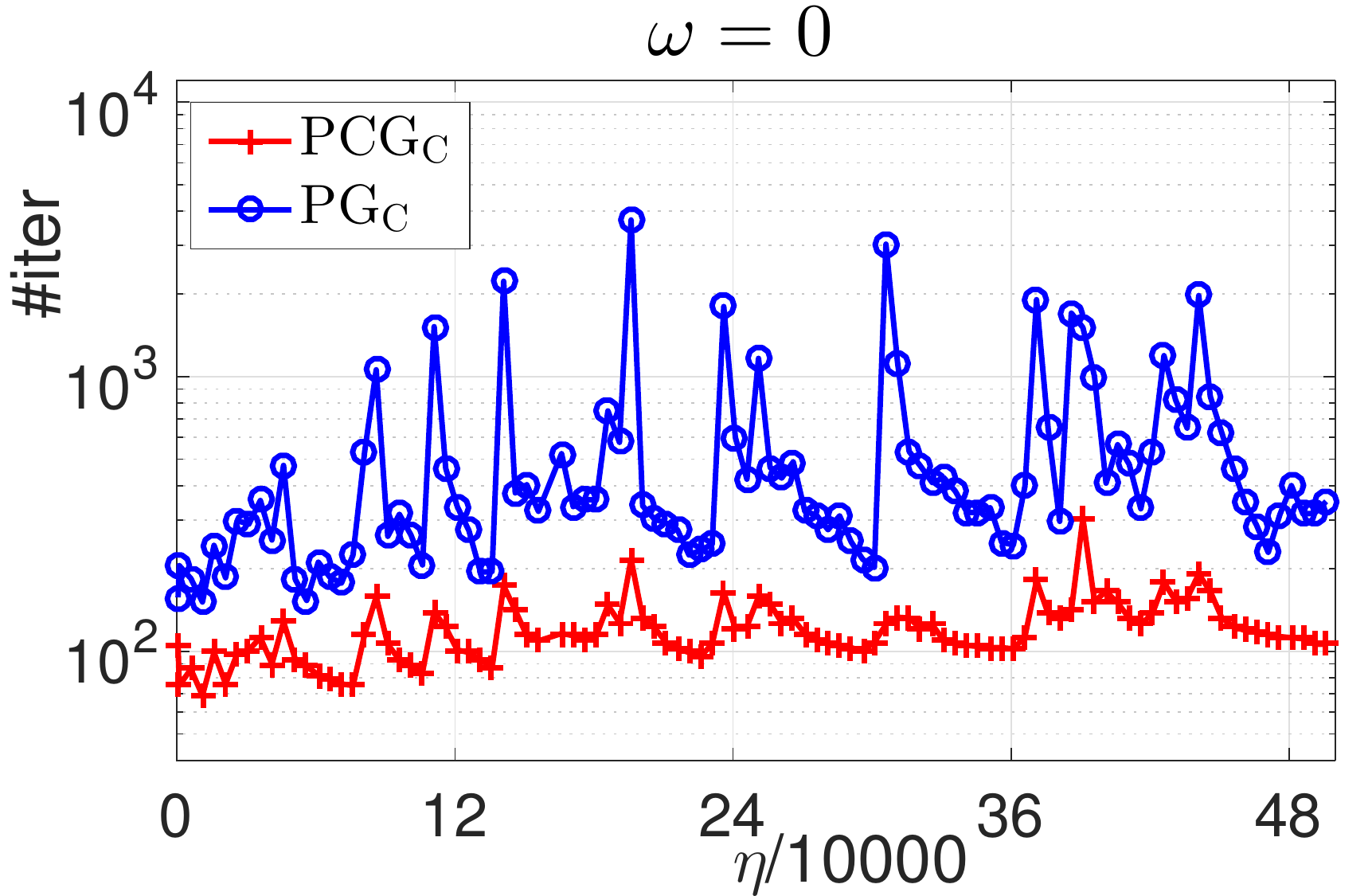}
\qquad
\includegraphics[height=5.2cm,width=8cm]{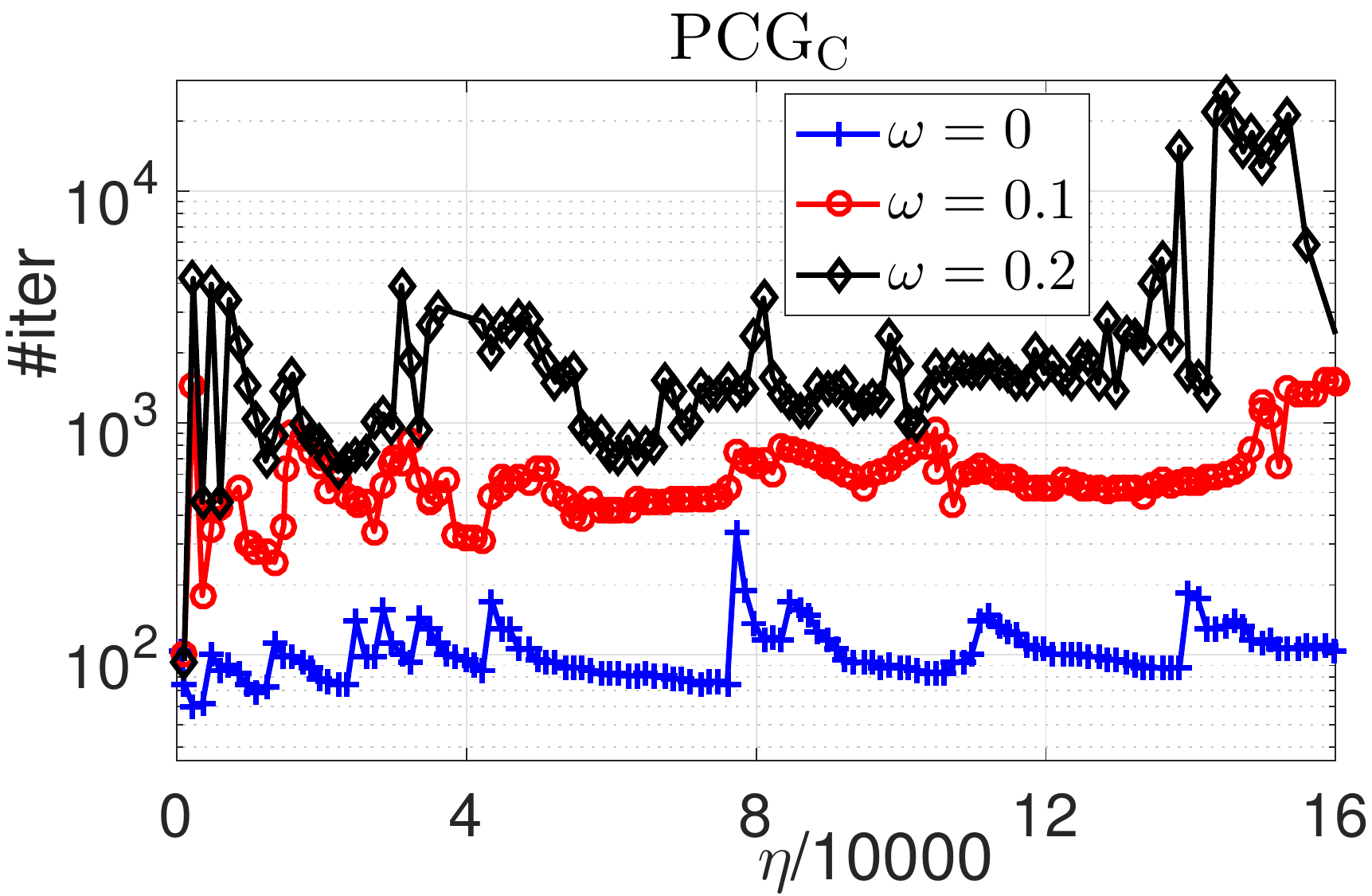}
}
\vspace{-0.2cm}
\caption{Exemple \ref{eg:2D_IterNum_PG_PCG_Og0}. Number of iterations for {\rm PCG}$_{\rm C}$ and {\rm PG}$_{\rm C}$
 for $\og=0$ (left) and  {\rm PCG}$_{\rm C}$  for different $\og$ (right)  vs. $\eta.$}
\label{fig:og0_2d_Diff_Beta_Comp_PG_PCG}
\end{figure}

\begin{exmp}
\label{eg:2D_IterNum_PG_PCG_Og0}
{\rm
In this example, we compare 
 the  performance of  {\rm PCG}$_{\rm C}$ and {\rm PG}$_{\rm C}$ for the 2D rotating case.  To this end,  $V(\bx)$ is chosen as
 the  harmonic plus lattice potential  \eqref{Lattice_Poten}
with $\gm_x=\gm_y=1$, $k_x=k_y=25$ and $q_x=q_y=\fl{\pi}{2}$.  The computational domain and mesh sizes are chosen respectively as  $\mathcal{D}=[-32, 32]^2$ and
 $h=\fl{1}{8}$.   
Fig. \ref{fig:og0_2d_Diff_Beta_Comp_PG_PCG} (left) shows the iteration number 
of   {\rm PCG}$_{\rm C}$ and {\rm PG}$_{\rm C}$ vs. different values of $\eta$ for $\og=0$, while  Fig. \ref{fig:og0_2d_Diff_Beta_Comp_PG_PCG} (right) reports the  number of iterations of 
{\rm PCG}$_{\rm C}$ with respect to  $\eta$ and $\og$.
From this figure, we can see that: (i) Similarly to the 1D case, {\rm PCG}$_{\rm C}$  outperforms {\rm PG}$_{\rm C}$;
   (ii) For $\og=0$, the iteration number for
{\rm PCG}$_{\rm C}$  would oscillate  in a small regime, which indicates the very slight 
dependence with respect to the nonlinearity strength $\eta$. When $\og$ increases, the number of iterations increases for a fixed $\eta$. Meanwhile, the dependency on $\eta$ 
 becomes stronger as $\og$ increases. Let us remark here that it would be extremely interesting  to 
build a robust preconditioner  including the rotational effects
to get a weaker $\og$-dependence in terms of convergence.
}

\end{exmp}

\begin{figure}[h!]
\centerline{
\includegraphics[height=4.5cm,width=6cm]{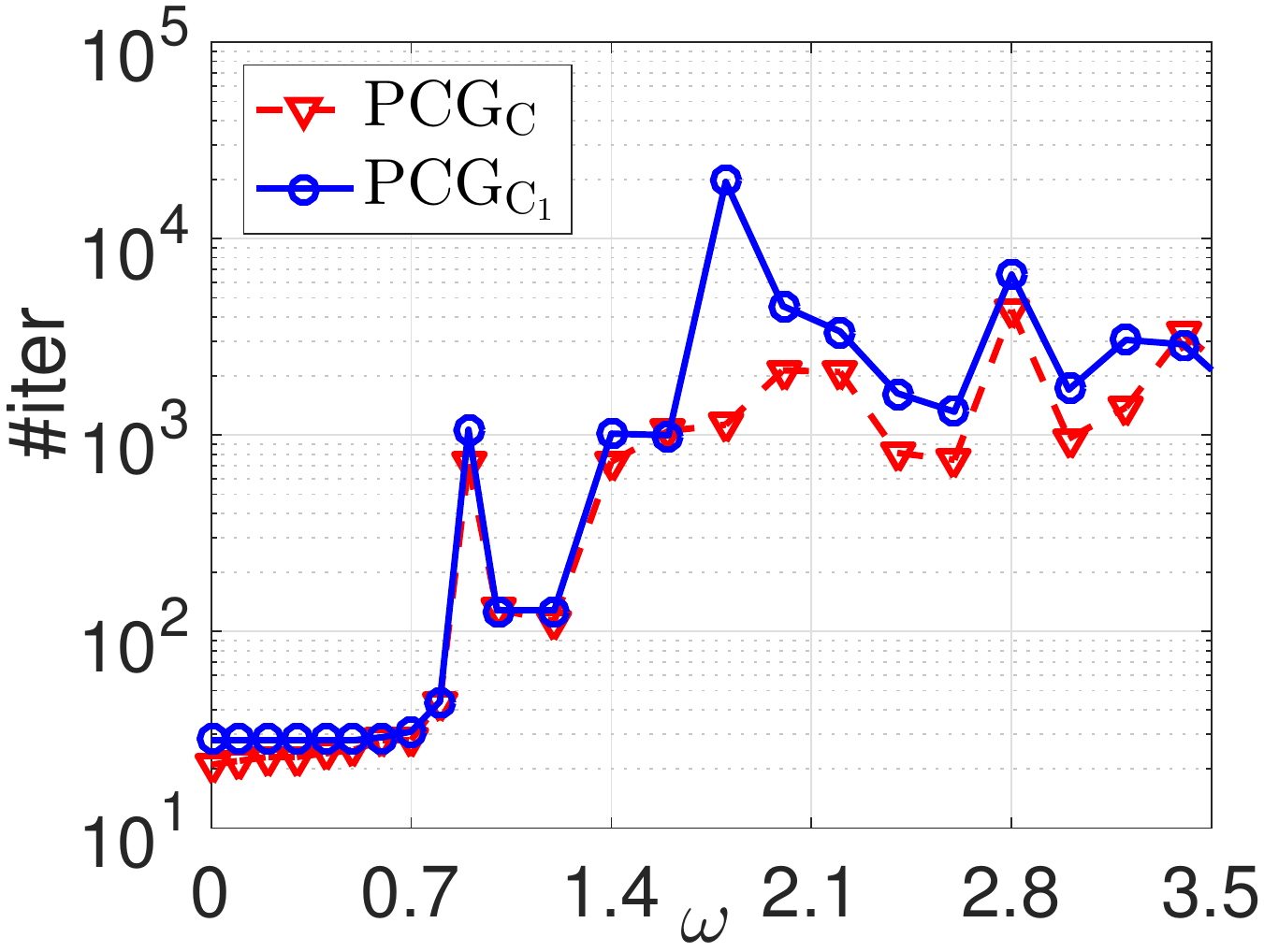}\qquad\qquad
\includegraphics[height=4.5cm,width=6cm]{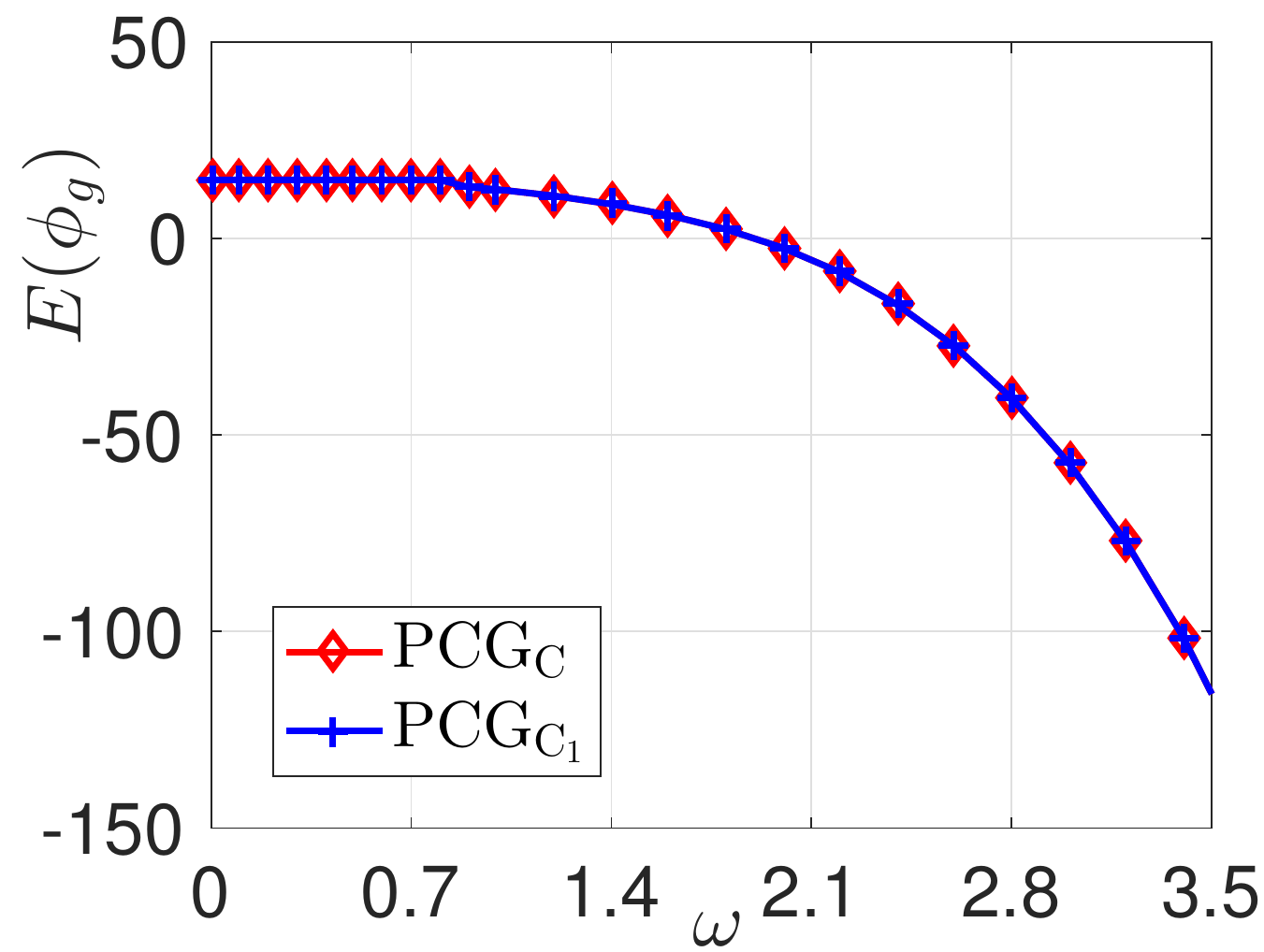}
}
\caption{Exemple \ref{eg:2D_omega}. Number of iterations for  {\rm PCG}$_{\rm C}$ and
  {\rm PCG}$_{\rm C_1}$   (left) and its corresponding total energies (right)  vs. $\og$.}
\label{fig:2D_Diff_Preconditioner_IterNum_vs_Og}
\end{figure}

\begin{exmp}
  \label{eg:2D_omega}
 {\rm
Following the previous example, here we compare the performance of 
 PCG$_{\rm C}$ and PCG$_{\rm C_1}$  for different values $\og$.
 To this end, we  fix $\eta=1000$ and vary $\og$ from 0 to 3.5. 
  Fig. \ref{fig:2D_Diff_Preconditioner_IterNum_vs_Og} 
 illustrates the number of iterations  of  these method  vs. different values of
 $\og$ and there corresponding energies.  From this figure and other
 experiments now  shown here, we see that
  (i)  All the methods converge to the stationary state with same
  energy;  (ii) The symmetrized preconditioner has more stable
  performance than the non-symmetric version, a fact we do not
  theoretically understand.
   }

 \end{exmp}

\begin{exmp}
  \label{eg:difficult_problem}
{\rm
In this example,  we apply PCG$_{\rm C}$ to solve some more difficult problems. We compute the ground states $\phi_{g}$ of 
rotating BECs with large values of $\eta$ and $\og$. To this end, we take $L=20$, $h=1/16$ and 
set the stopping tolerance  in \eqref{Stop_Energy} to  $\eps=10^{-14}$. Table  \ref{tab:CPU_Diff_Og} lists 
the CPU times for the PCG$_{\rm C}$ solver  
to converge while Fig. \ref{fig:2D_density_diff_og_quartic_poten} shows the  contour plot of the 
   density function $|\phi_g(\bx)|^2$ for different $\og$  and $\eta$.
We can see that  the PCG$_{\rm C}$ method
   converges very fast to the stationary states. Let us remark that, to the best of our knowledge,
   only a few results were reported 
   for such  fast rotating BECs with highly nonlinear (very large $\eta$) problems, although they are actually more relevant for 
   real physical problems.  Hence, PCG$_{\rm C}$ can tackle efficiently difficult realistic problems on a laptop.
   }

\begin{table}[h!]
\tabcolsep 12pt  
\caption{CPUs  time ({\it seconds}) for PCG$_{\rm C}$ to compute the  ground states of the GPE with various $\og$ and $\eta$.}
\label{tab:CPU_Diff_Og}
\begin{center}\vspace{-1.5em}
\def\temptablewidth{0.95\textwidth}
{\rule{\temptablewidth}{1pt}}
\begin{tabularx}{\temptablewidth}{@{\extracolsep{\fill}}p{0.75cm}|cccccccc}
  $\eta$    &  $\og$=1  & 1.5 & 2& 2.5  &  3 & 3.5 & 4 & 4.5   \\
\hline
1000    &      493      &   551    &   560   &   2892   &    2337    &   720    &  966     &  3249 \\[0.05cm]
5000    &        1006   &   1706    &   867   & 6023     &     1144   &    1526   &  12514   &19248   \\[0.05cm]
10000    &  4347         & 21525      &   5511   &15913  &15909    &   6340     &     16804   & 32583   \\
  \hline
\end{tabularx}
{\rule{\temptablewidth}{1pt}}
\end{center}
\end{table}

\begin{figure}[h!]
\centerline{
\includegraphics[height=3.5cm,width=3.5cm]{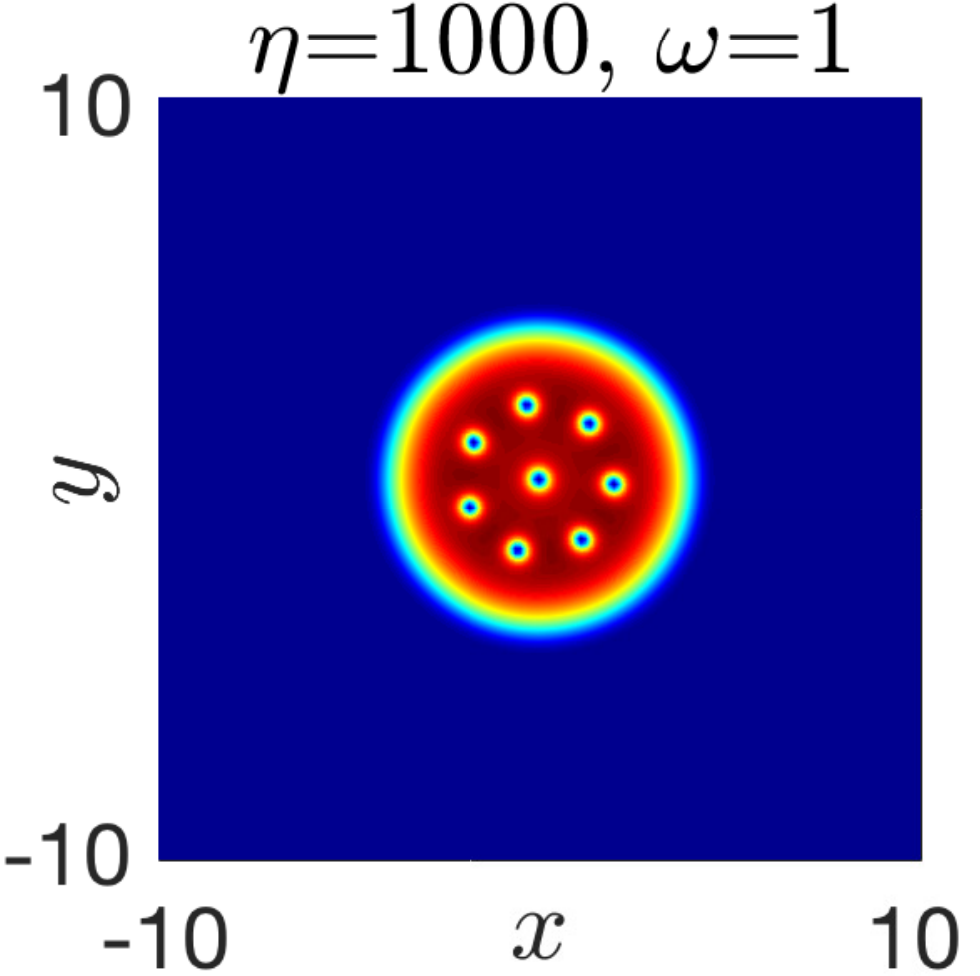}
\includegraphics[height=3.5cm,width=3.5cm]{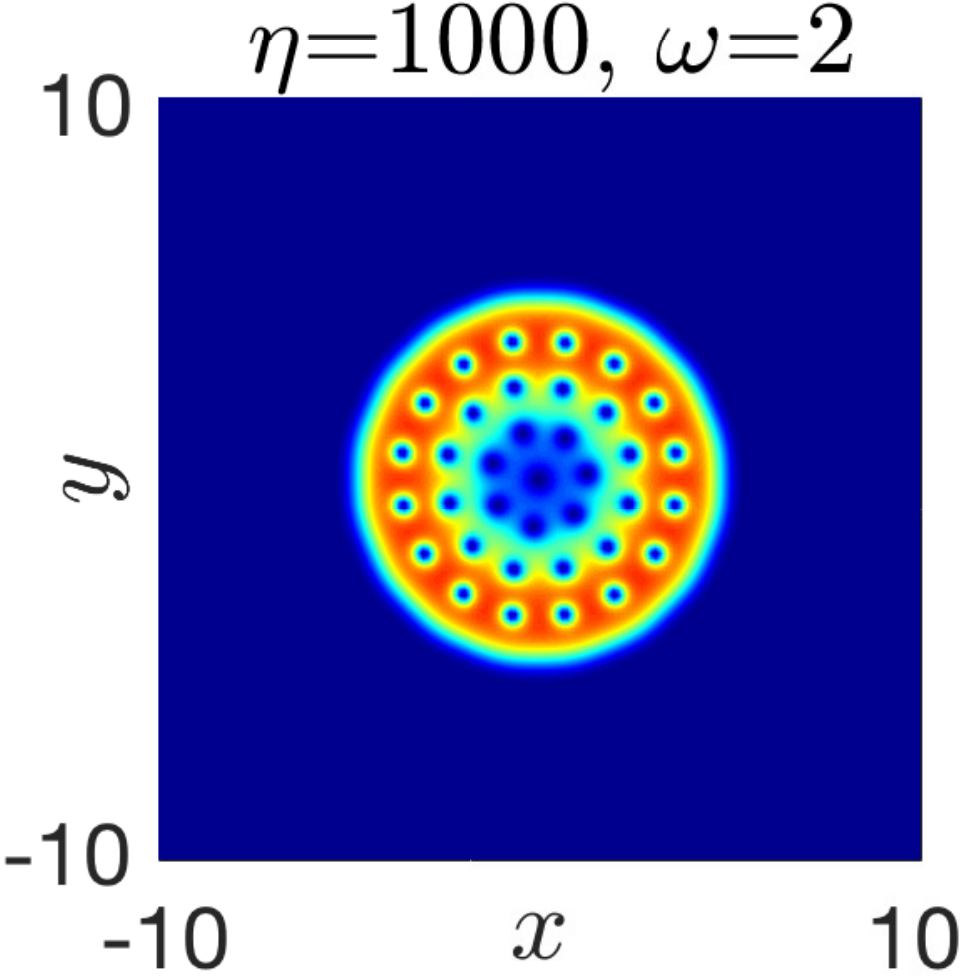}
\includegraphics[height=3.5cm,width=3.5cm]{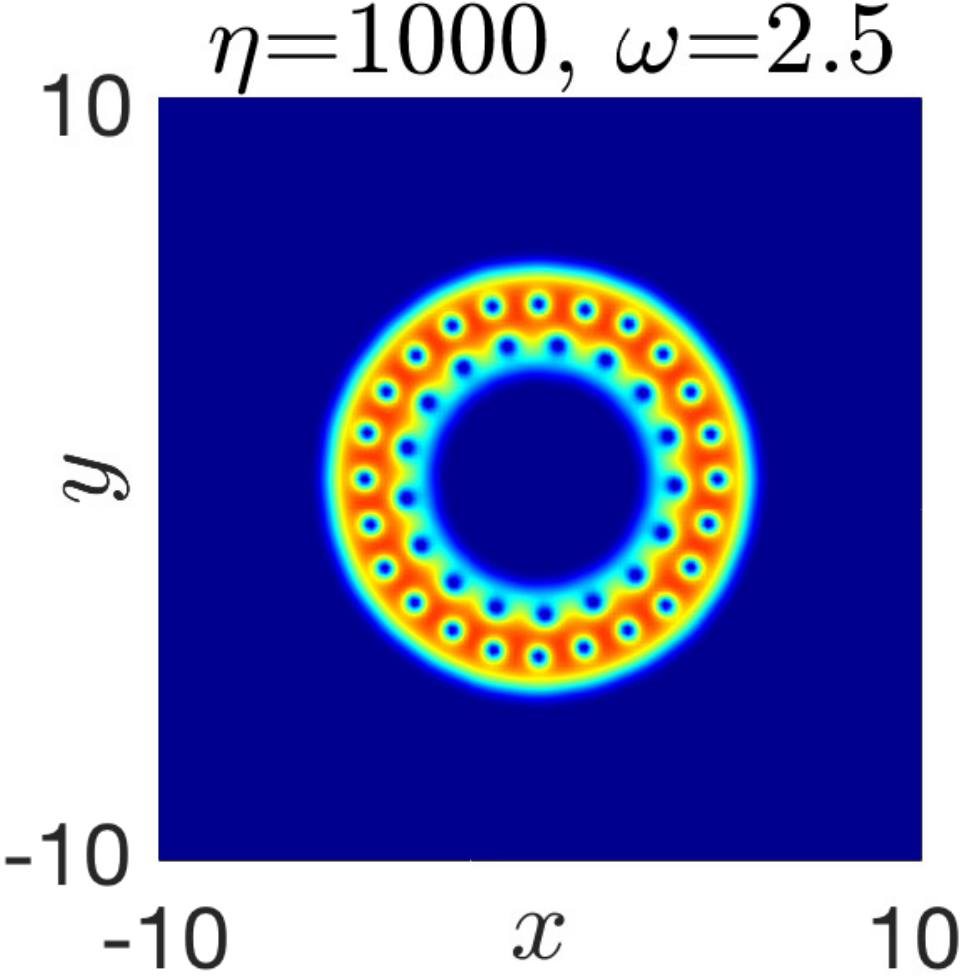}
\includegraphics[height=3.5cm,width=4.8cm]{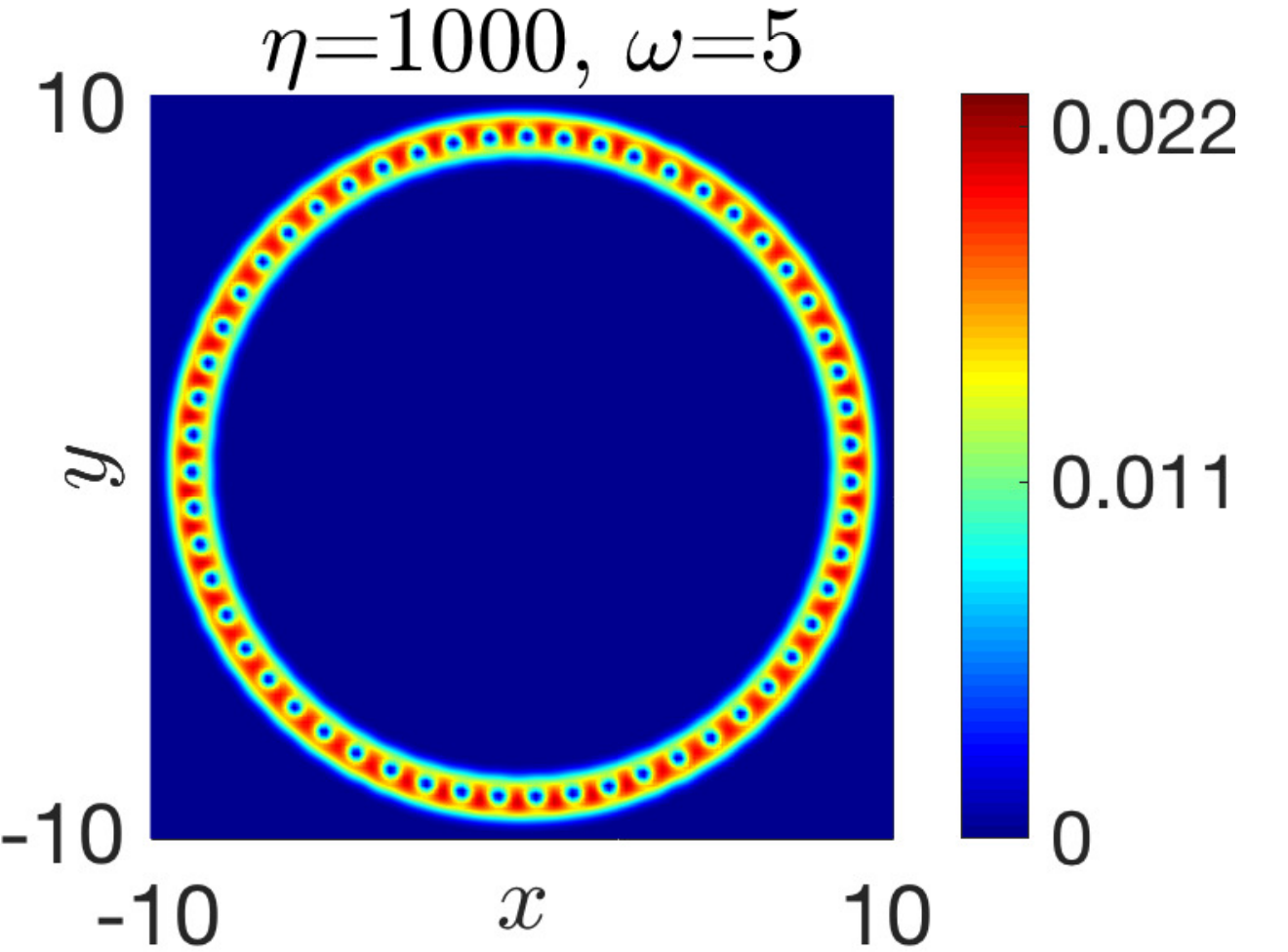}
}
\vspace{0.3cm}
\centerline{
\includegraphics[height=3.5cm,width=3.5cm]{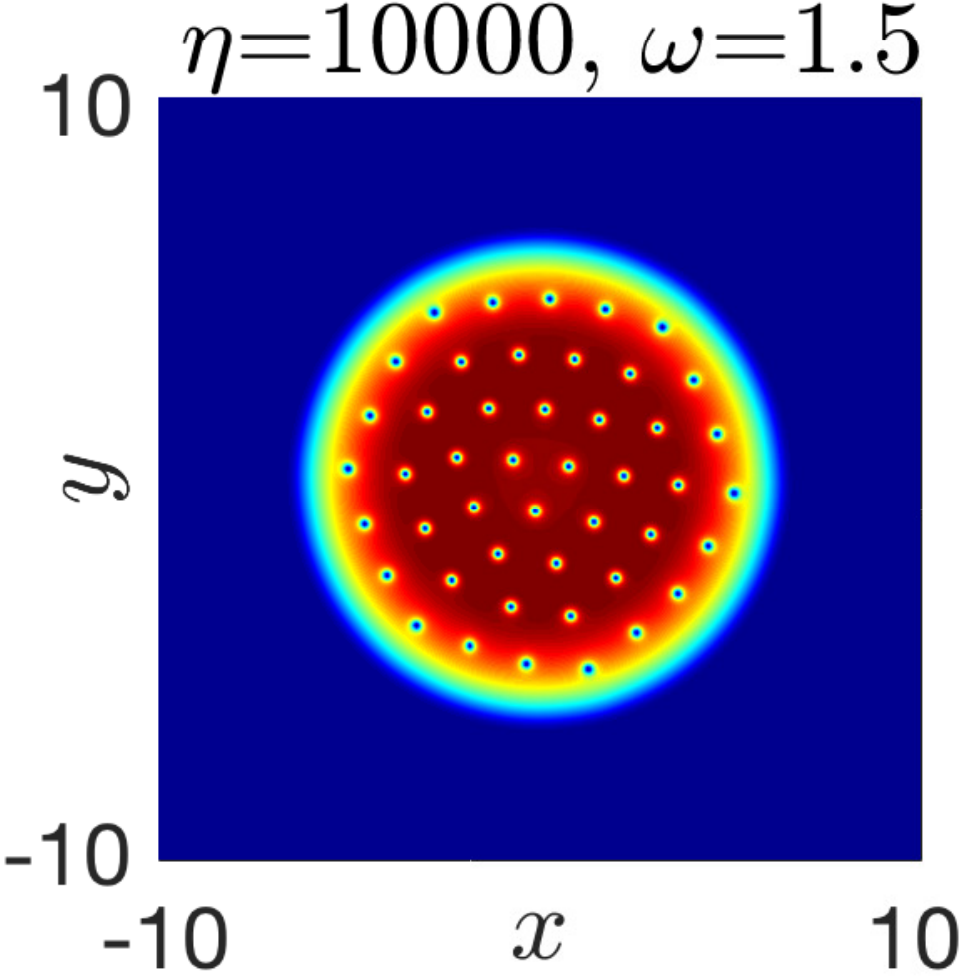}
\includegraphics[height=3.5cm,width=3.5cm]{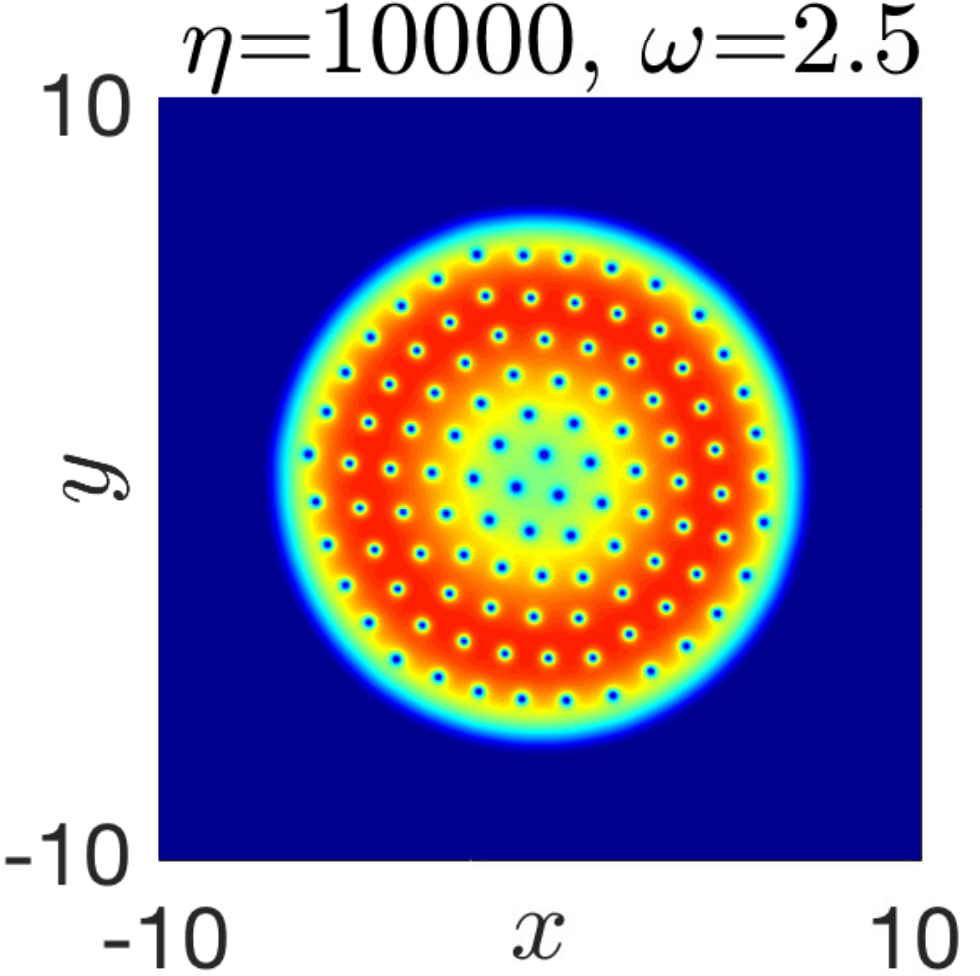}
\includegraphics[height=3.5cm,width=3.5cm]{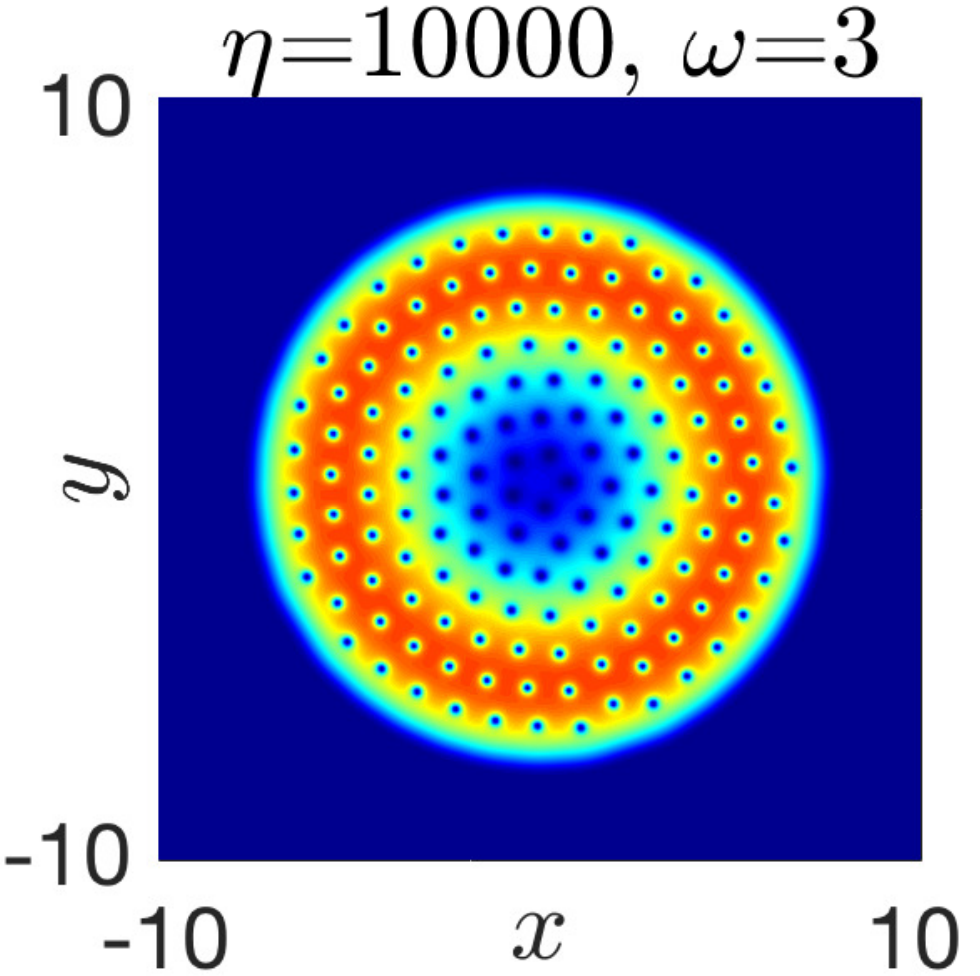}
\includegraphics[height=3.5cm,width=4.8cm]{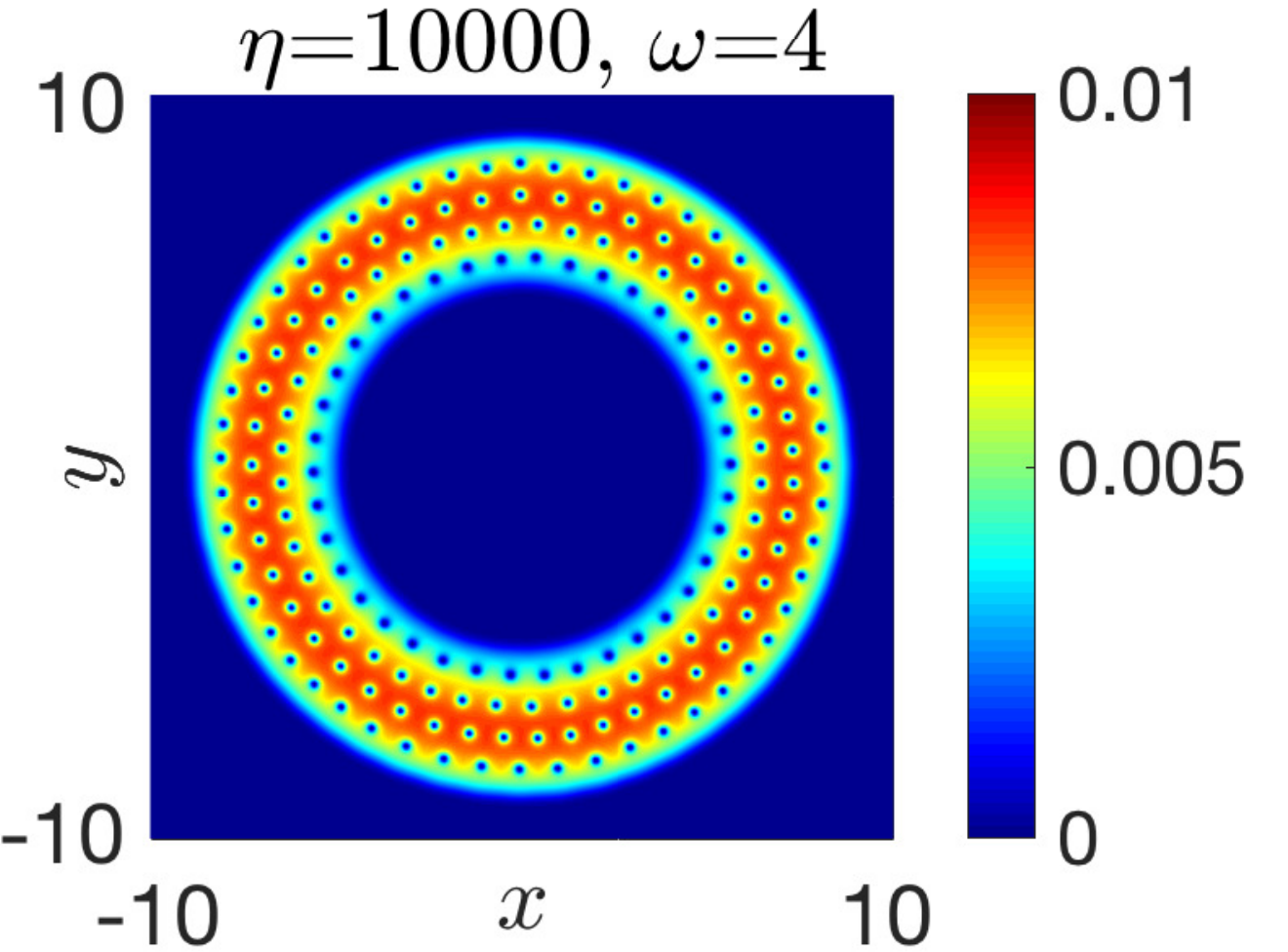}
}
\vspace{-0.1cm}
\caption{Exemple \ref{eg:difficult_problem}. Corresponding contour plots of the density function $|\phi_g(\bx)|^2$ of Table \ref{tab:CPU_Diff_Og}.}
\label{fig:2D_density_diff_og_quartic_poten}
\end{figure}

 \end{exmp}

\begin{exmp}
\label{eg:diff_initial_fix_grid}
{\rm 
The  choice of the initial data  also affects the final converged stationary states.   Since
all the algorithms we discussed are local minimisation algorithms, inappropriate initial guess might lead to local minimum. 
To illustrate this claim, we take $\eps=10^{-14}$, $V(\bx)=\fl{|\bx|^2}{2}$, $\eta=500$ and compute the ground 
states of the rotating GPE for  different  $\og$ and 10 types of frequently used initial data
\bea
&&  (a) \ \phi_a(\bx)=\fl{1}{\sqrt{\pi}} e^{-(x^2+y^2)/2},   \qquad  (b) \ \phi_b(\bx)=(x+iy)\phi_a(\bx),
\qquad (\bar{b}) \  \phi_{\bar{b}}(\bx)=\bar{\phi}_{b}(\bx), \\
&&  (c)\ \phi_c=\fl{(\phi_a(\bx)+\phi_b(\bx))/2}{\|(\phi_a(\bx)+\phi_b(\bx))/2\|},   \qquad\qquad
 (d) \ \phi_d=\fl{(1-\og)\phi_a(\bx)+\og\phi_b(\bx)}{\|(1-\og)\phi_a(\bx)+\og\phi_b(\bx)\|},  \\
&& (e) \ \phi_e=\fl{\og\phi_a(\bx)+(1-\og)\phi_b(\bx)}{\og\phi_a(\bx)+(1-\og)\phi_b(\bx)\|},\qquad
 (\bar{c})\  \phi_{\bar{c}}(\bx)=\bar{\phi}_{c}(\bx),  \qquad
  (\bar{d})\  \phi_{\bar{d}}(\bx)=\bar{\phi}_{d}(\bx), \\
&&    (\bar{e})\  \phi_{\bar{e}}(\bx)=\bar{\phi}_{e}(\bx),    \qquad\qquad\qquad\qquad
(f)\ {\rm Thomas\ Fermi\ approximation}\; \eqref{TF_App}.
\eea
Another approach to prepare some initial data is as follows: we first consider one of the above initial guess (a)-(f),
next compute the ground state on a coarse spatial  grid, say with a  number of grid points $M_p\times M_{p}$ (with $M_{p}=2^{p}$),
and then denote the corresponding stationary state by $\phi_g^{p}$. We next subdivide the grid with 
 $M_{p+1}\times M_{p+1}$ points,  interpolate  $\phi_g^{p}$ on the refined grid $M_{p+1}\times M_{p+1}$ to get a new initial data and launch the algorithm
 at level $p+1$, and so on  until the finest grid with $M\times M$ points where the converged solution is still denoted by
 $\phi_{g}$.
Similarly to \cite{weizhu2015}, this  multigrid  technique is  applied here with the
 coarsest grid based with $M_{6}=2^{6}$ and ends with the finest
grid $M=2^{9}$. We use the tolerance parameters  $\eps=10^{-14}$ for $M=2^{9}$, and $\eps=10^{-12}$ for $p=6, 7, 8$. 

Tables  \ref{tab:2d_diff_init_data_fix_grid} and \ref{tab:diff_init_data_multi_grid} list  the energies obtained by {\rm PCG}$_{\rm C}$
\textit{via} the fixed and multigrid approaches, respectively,  for different initial data and $\og$.  The stationary states $\phi_g(\bx)$ with lowest energies are marked by underlines and the corresponding CPU times are listed in the same Table.  Moreover,  
Fig. \ref{fig:2D_density_harmonic_diff_initial_data_fix_grids} shows the contour plots 
 $|\phi_g(\bx)|^2$ of the  converged solution with lowest energy obtained by  the multigrid approach.
 Now, let us denote by $E^{p}_{n}:=E(\phi_{n}^p)$ the evaluated energy at step $n$ for a discretization level
 $p=6, 7, 8$, and let $E_{g}=E(\phi_{g})$ the energy for the converged stationary state for the finest grid. 
 Then, we represent on  
  Fig. \ref{fig:2D_Multi_grid_Error_vs_CPU} 
  the evolution of $\log_{10}(|E^{p}_{n}-E_{g}|)$ vs. the CPU time for a rotating velocity  $\og=0.95$.
  For comparison, we also show the corresponding evolution  obtained by the fixed grid approach. 
The contour plots of $|\phi^{p}_g(\bx)|^2$ obtained for each
 intermediate coarse
  grid for $p=6,7,8$, and the initial guess   are also reported. 
  
 From these Tables and Figures, we can see that:  (i) Usually, the PCG$_{\rm C}$ algorithm with an initial data of  type ($d$) or ($\bar{d}$)  converges to the stationary state
  of lowest energy; 
 (ii) The multigrid approach is more robust  than the fixed grid approach in terms of CPU time and possibility to obtain a stationary 
 state with lower energy.

   }
   
\begin{table}[h!]
\tabcolsep 12pt  
\caption{Exemple \ref{eg:diff_initial_fix_grid}. Fixed grid  approach  (with $M=2^{9}$):  
converged energies and the CPU times ({\it seconds}) for the solution with lowest energy  (which is underlined). }
\label{tab:2d_diff_init_data_fix_grid} 
\begin{center}
\vspace{-1.8em}
\def\temptablewidth{0.98\textwidth}
{\rule{\temptablewidth}{1pt}}
\setlength{\tabcolsep}{0.5em}
\begin{tabularx}{\temptablewidth}{c|cccccccccc|c}
$\og$   & (a)       & (b)      & (b2)      &  (c)  &    (c2)    &   (d)  &  (d2)   &   (e) 	&  (e2)    &  (f) &  CPU   \\ \hline  	      
0.5 &    8.5118    &  8.2606   & 9.2606  &  8.0246 &    8.0197&   8.0246  &  $\underline{8.0197}$   &   8.0246  &    8.0197& 8.0246   &176.0    \\[0.05cm]
0.6 &    8.5118      & 8.1606   &  9.3606 &7.5845   & 7.5910   &7.5845    &  $\underline{7.5845 }$    &  7.5845   &    7.5910& 7.5845&310.7     \\[0.05cm]
0.7 &     8.5118     &8.0606     &9.4606   &6.9754    & $\underline{6.9731 }$   & 6.9792    &  6.9754   &6.9754    & 6.9792  & 6.9767& 542.4    \\[0.05cm]
0.8 &     8.5118    &7.9606     &9.5606   &6.1016   &$\underline{ 6.0997}$   &6.1031    &6.1031     &  6.1040   & 6.1019  & 6.1016& 417.0    \\[0.05cm]
0.9 &     8.5118    &7.8606     & 9.6606  & 4.7777  & 4.7777    & 4.7777     &$\underline{4.7777}$    &4.7777     & 4.7777   & 4.7777& 1051.1    \\[0.05cm]
0.95 &   8.5118   &  7.8106    &  9.7106  &3.7414   &3.7414    &$\underline{3.7414 }$ &3.7414        & 3.7414     & 3.7414    &3.7414 &  3280.5  \\[0.05cm]
\end{tabularx}
{\rule{\temptablewidth}{1pt}}
\end{center}
\end{table}

\begin{table}[h!]
\tabcolsep 12pt  
\caption{ Exemple \ref{eg:diff_initial_fix_grid}. Multigrid approach  (starting from the coarsest level $p=6$ to the finest level $p=9$):
converged energies and the CPU times ({\it seconds}) for the solution with  lowest energy (underline). }
\label{tab:diff_init_data_multi_grid} 
\begin{center}
\vspace{-1.8em}
\def\temptablewidth{0.98\textwidth}
{\rule{\temptablewidth}{1pt}}
\setlength{\tabcolsep}{0.5em}
\begin{tabularx}{\temptablewidth}{c|cccccccccc|c}
$\og$   & (a)       & (b)      & (b2)      &  (c)  &    (c2)    &   (d)  &  (d2)   &   (e) 	&  (e2)    &  (f) &  CPU   \\ \hline  	      
0.5 &    8.0246    &  8.0197   & 8.0197  &  8.0197 &    8.0197 &  8.0197&   $\underline{8.0197}$    &   8.0197  &    8.0197& 8.0257   & 29.5  \\[0.05cm]
0.6 &    7.5845     & 7.5845    &  7.5910 &7.5845   & 7.5910   &7.5890     &  $\underline{7.5845 }$  &  7.5845   &    7.5910& 7.5845& 32.3  \\[0.05cm]
0.7 &     6.9767     &$\underline{6.9726 }$    &6.9792   &6.9754    & 6.9731    & 6.9731    &  6.9731   &6.9757    & 6.9731  & 6.9731& 53.3   \\[0.05cm]
0.8 &     6.1019    &6.1031     &6.1019   &$\underline{6.0997}$   & 6.1016   &6.1016     &6.1016     &  6.1019   & 6.1016   & 6.0997& 75.2   \\[0.05cm]
0.9 &     4.7777    &4.7777     & 4.7777  & 4.7777  & 4.7777    & $\underline{4.7777}$   &4.7777     &4.7777     & 4.7777   & 4.7777&  238.1  \\[0.05cm]
0.95 &     3.7414   &3.7414     &3.7414   &3.7414   &3.7414   &$\underline{3.7414 }$  &3.7414        &3.7414     & 3.7414   &3.7414 & 621.9  \\[0.05cm]
\end{tabularx}
{\rule{\temptablewidth}{1pt}}
\end{center}
\end{table}

\begin{figure}[h!]
\centerline{
\includegraphics[height=3.5cm,width=4.8cm]{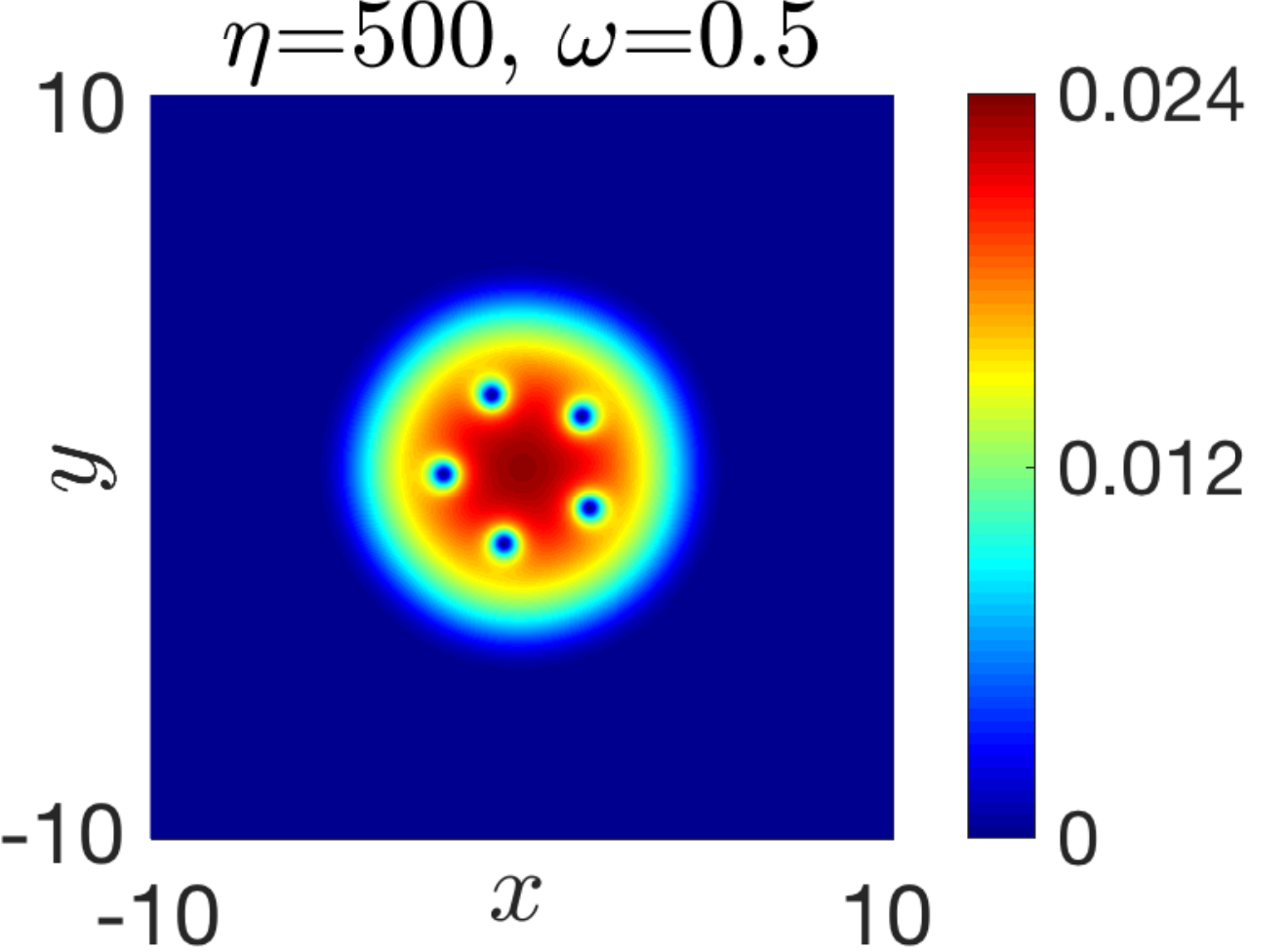}\quad
\includegraphics[height=3.5cm,width=4.8cm]{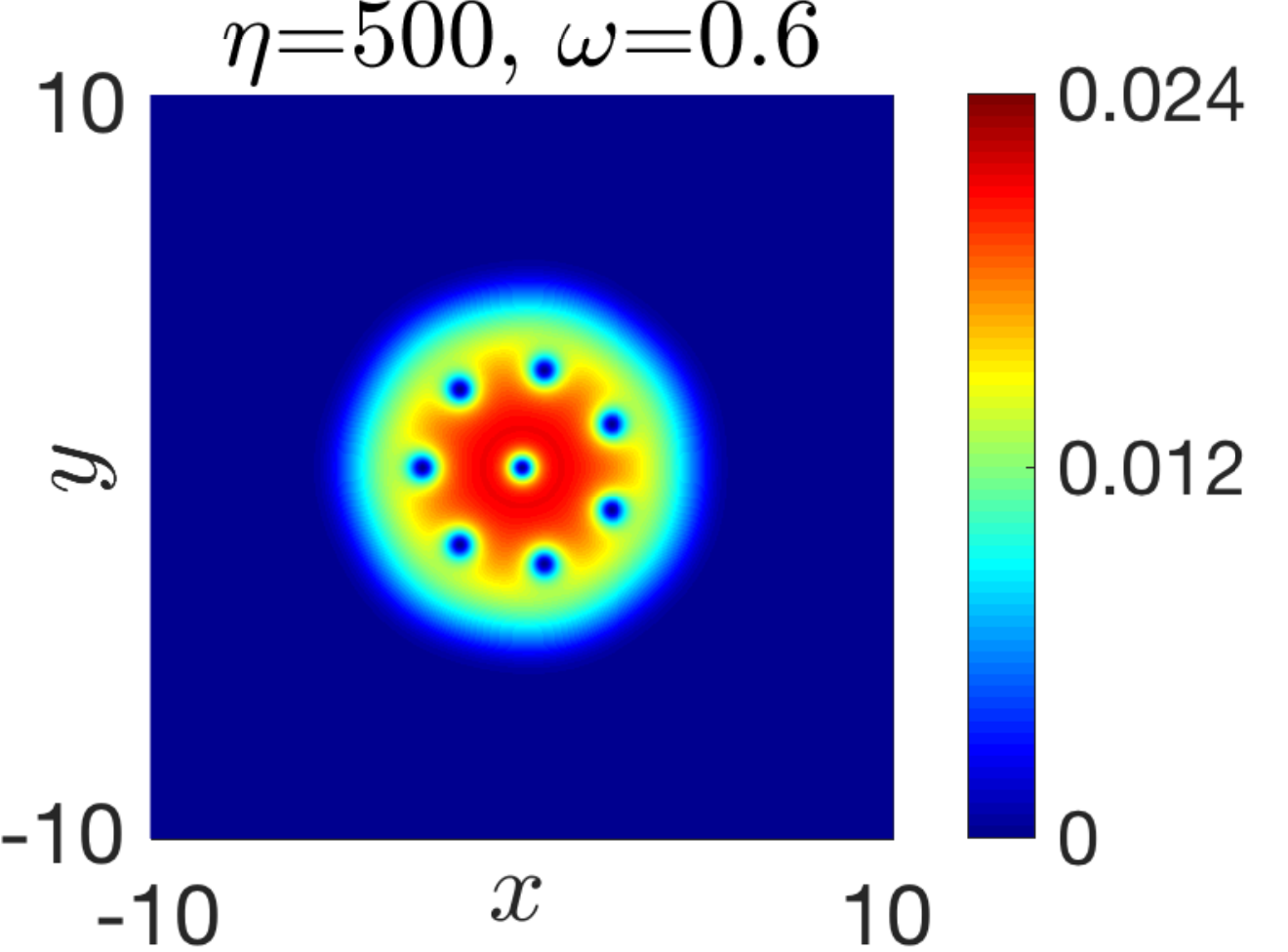}\quad
\includegraphics[height=3.5cm,width=4.8cm]{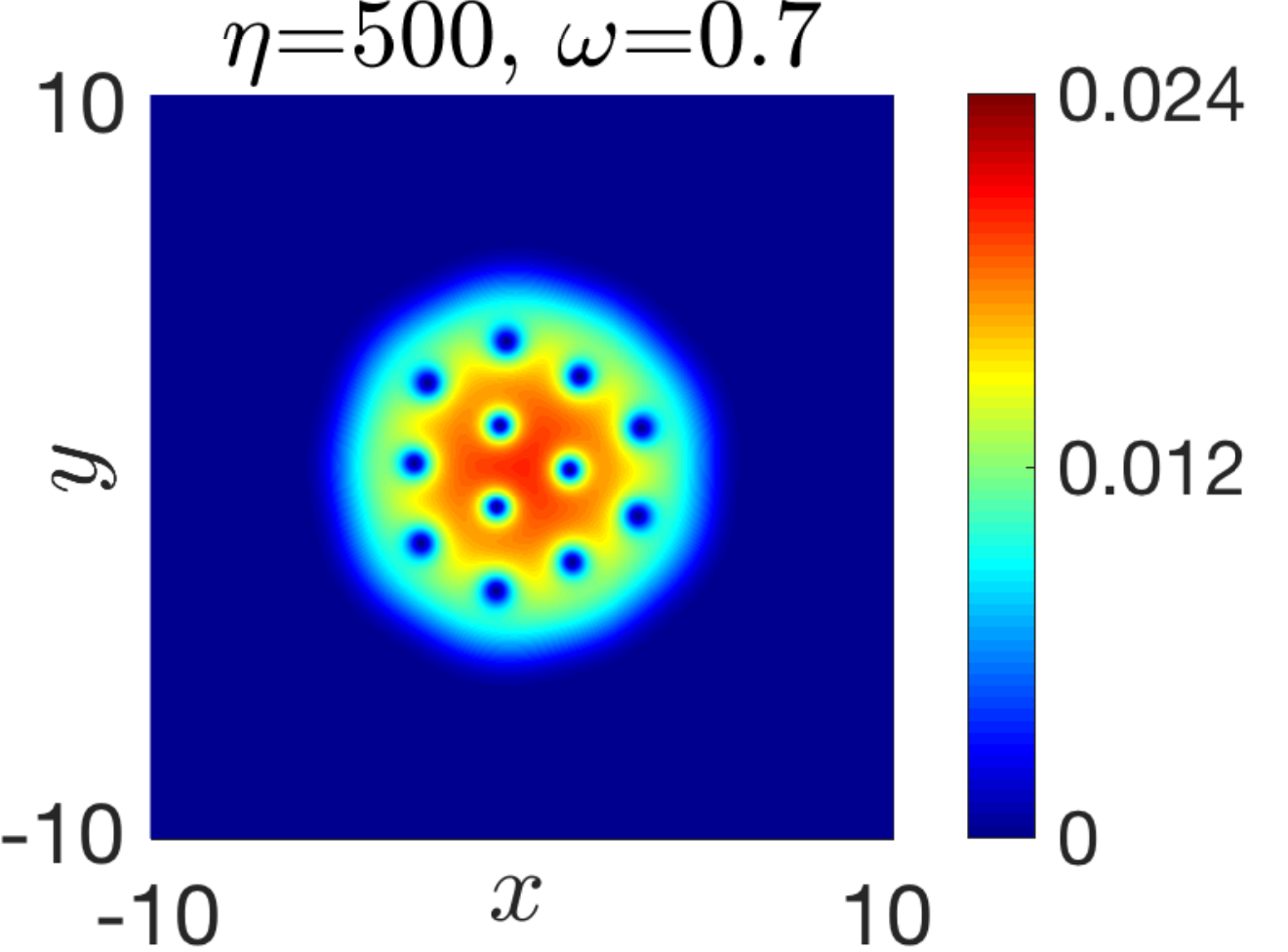}
}
\vspace{0.4cm}
\centerline{
\includegraphics[height=3.5cm,width=4.8cm]{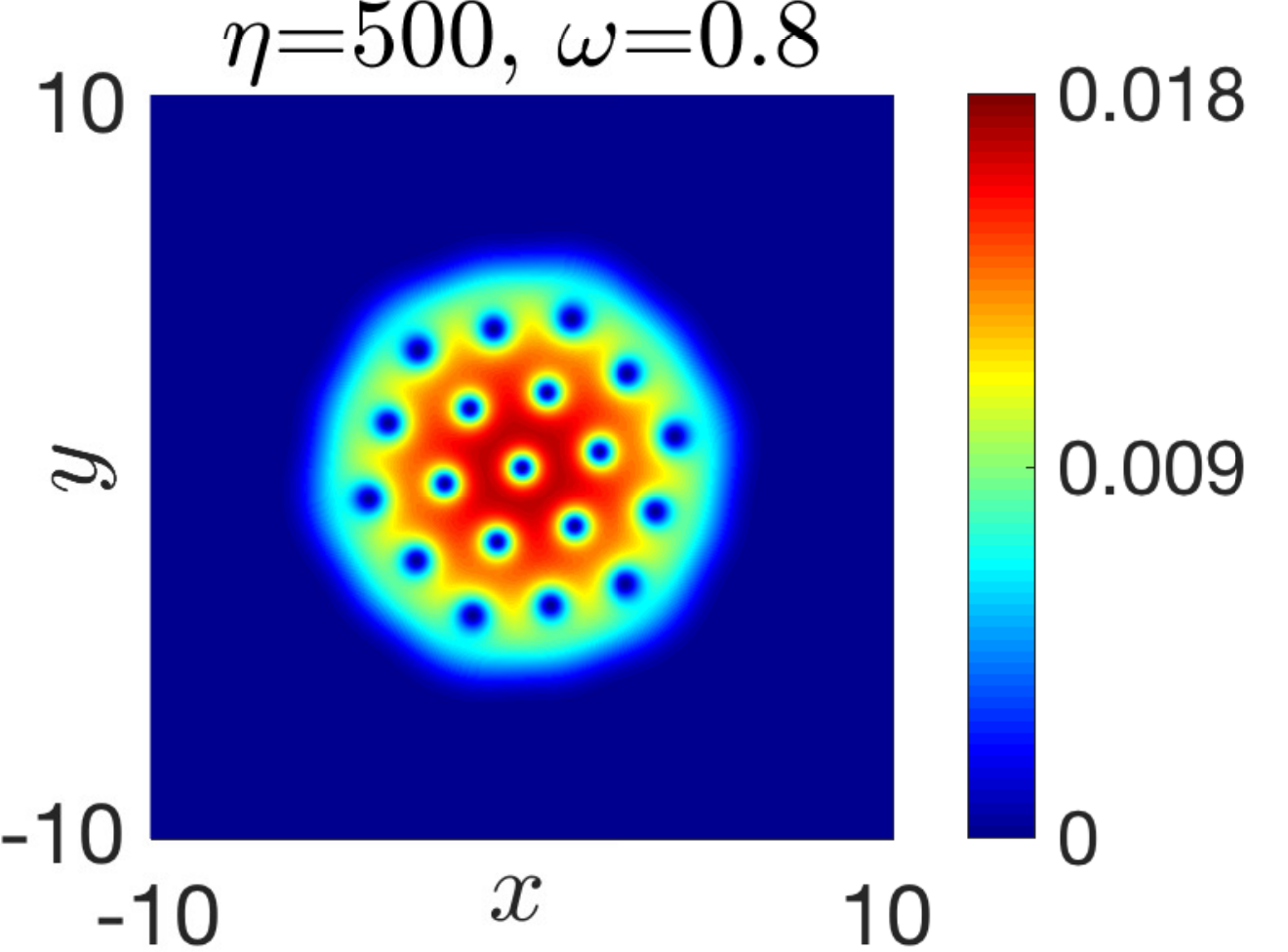}\quad
\includegraphics[height=3.5cm,width=4.8cm]{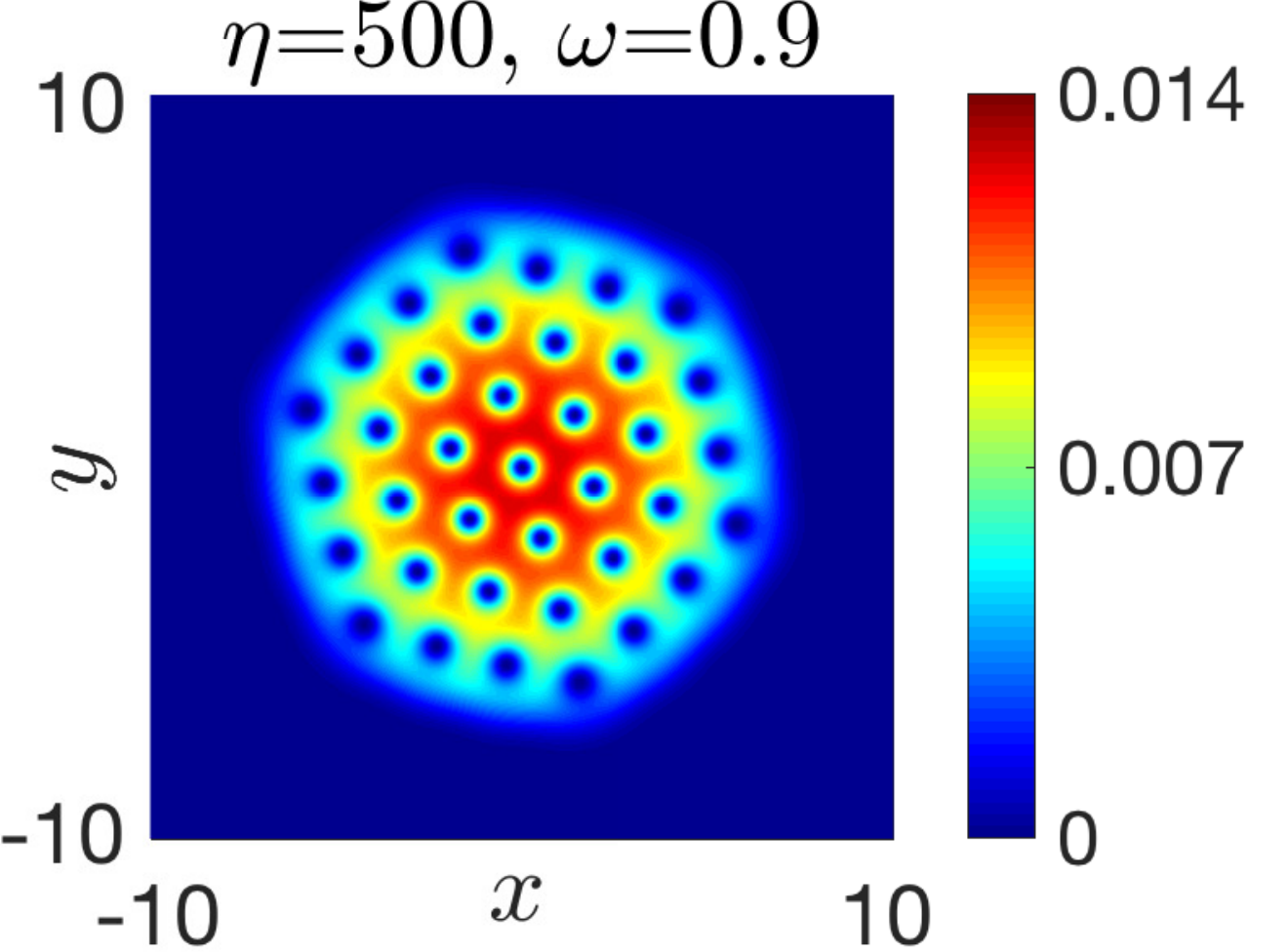}\quad
\includegraphics[height=3.5cm,width=4.8cm]{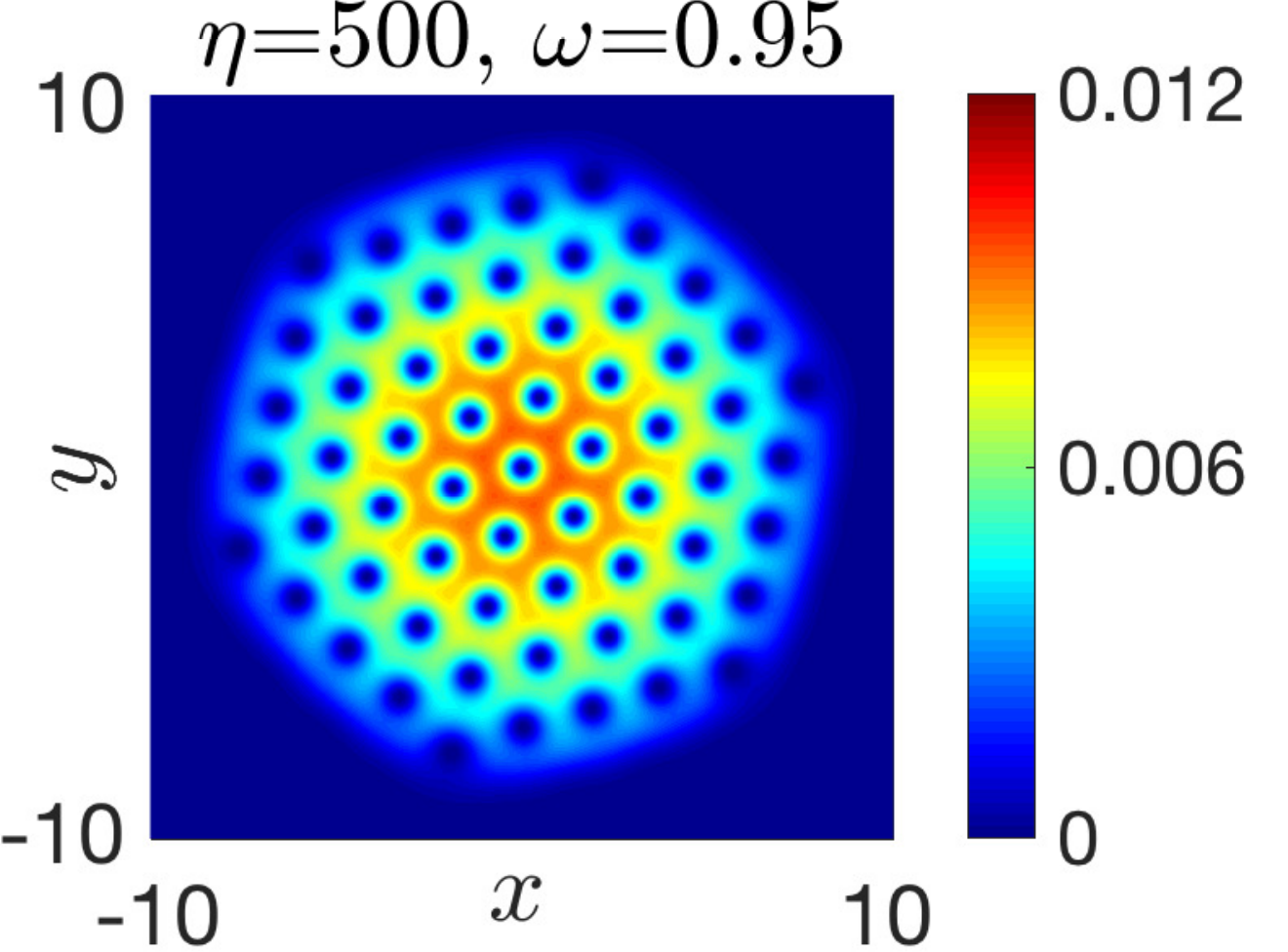}
}
\caption{Contour plots of $|\phi_g(\bx)|^2$  corresponding to the lowest energy levels  in Table  \ref{tab:diff_init_data_multi_grid}.}
\label{fig:2D_density_harmonic_diff_initial_data_fix_grids}
\end{figure}

\begin{figure}[h!]
\centerline{
\includegraphics[height=4.2cm,width=7.1cm]{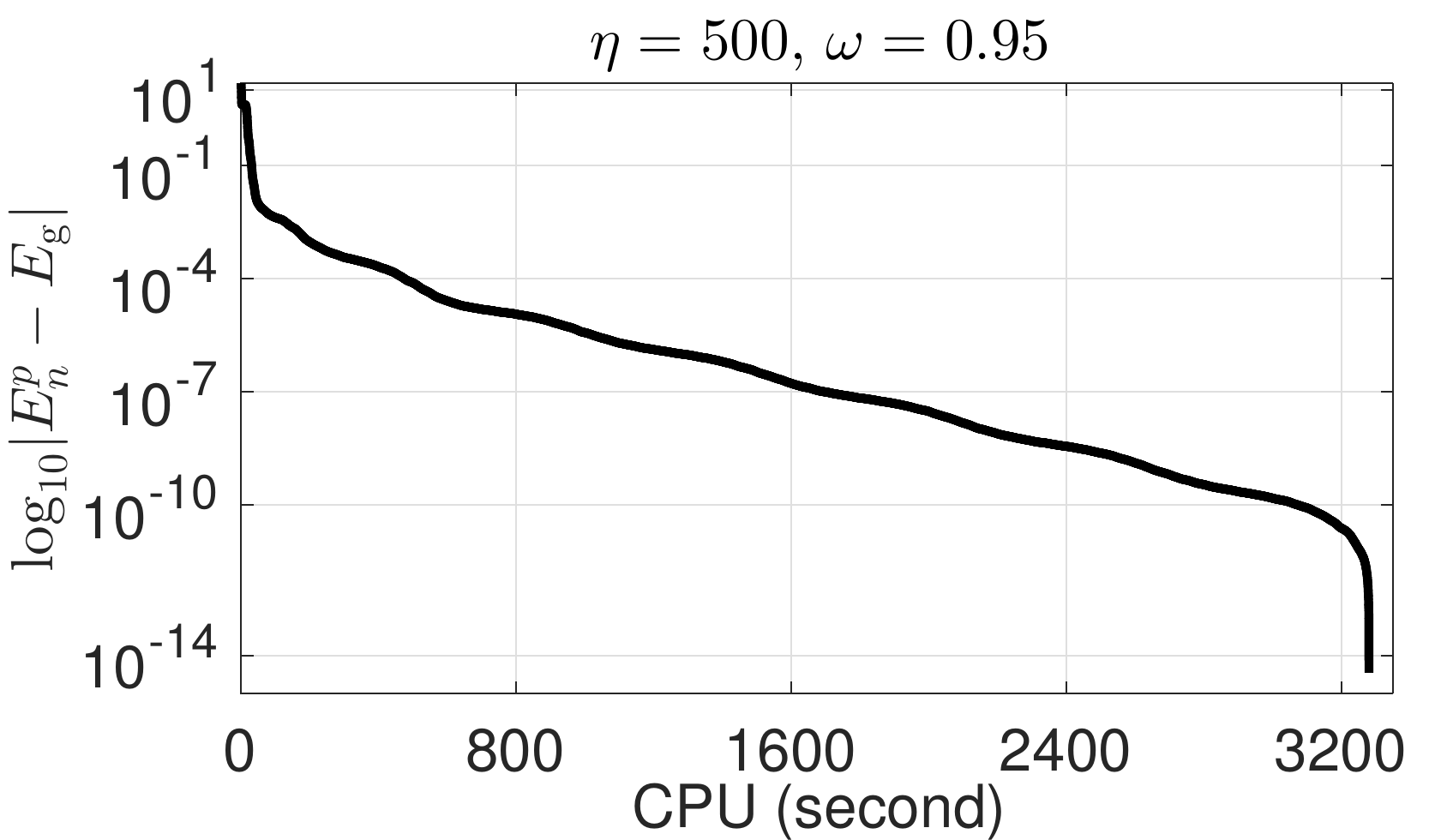}\quad
 \includegraphics[height=4.2cm,width=7.1cm]{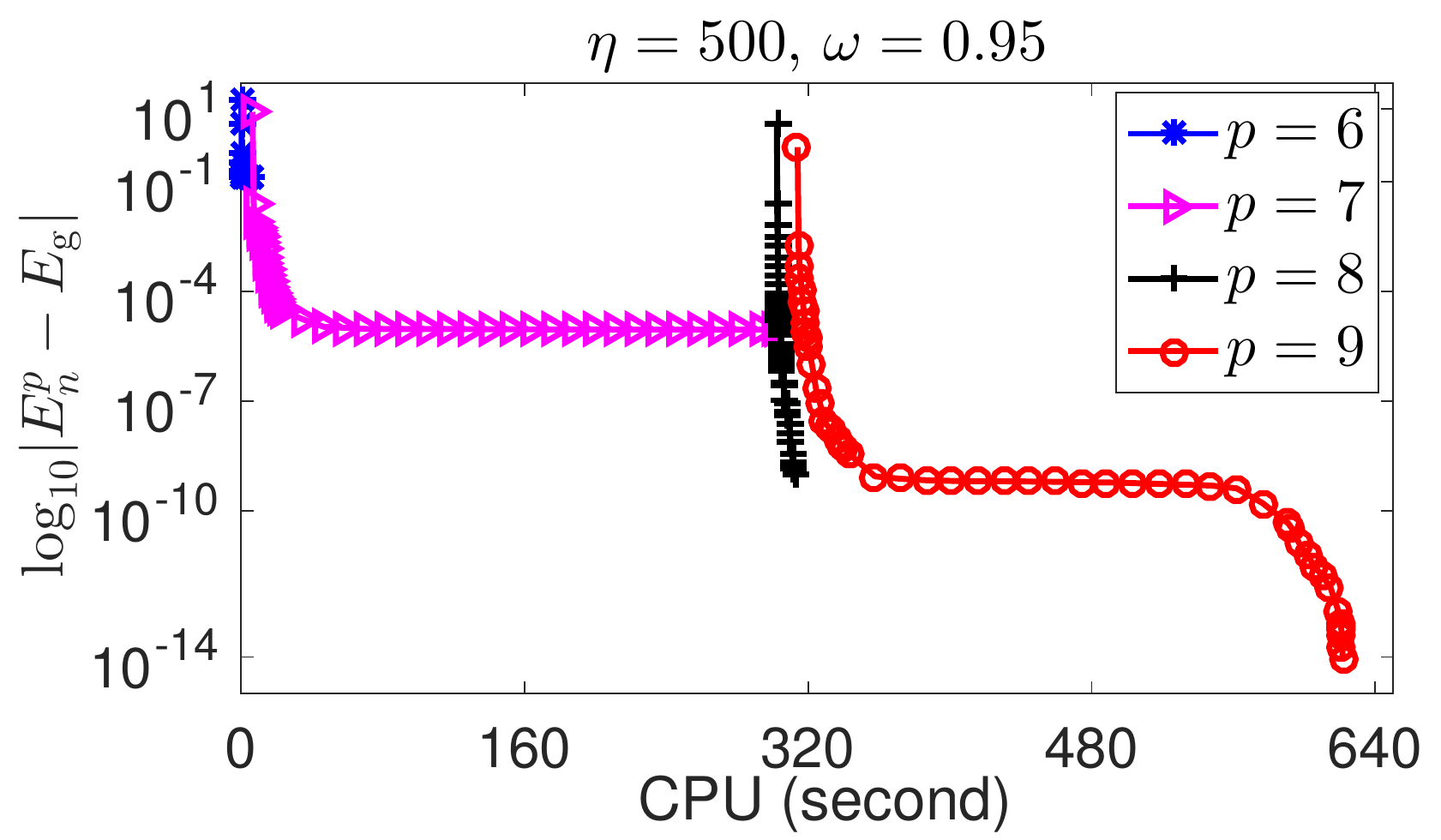}
}
\vspace{0.4cm}
\centerline{
\includegraphics[height=3.45cm,width=4.2cm]{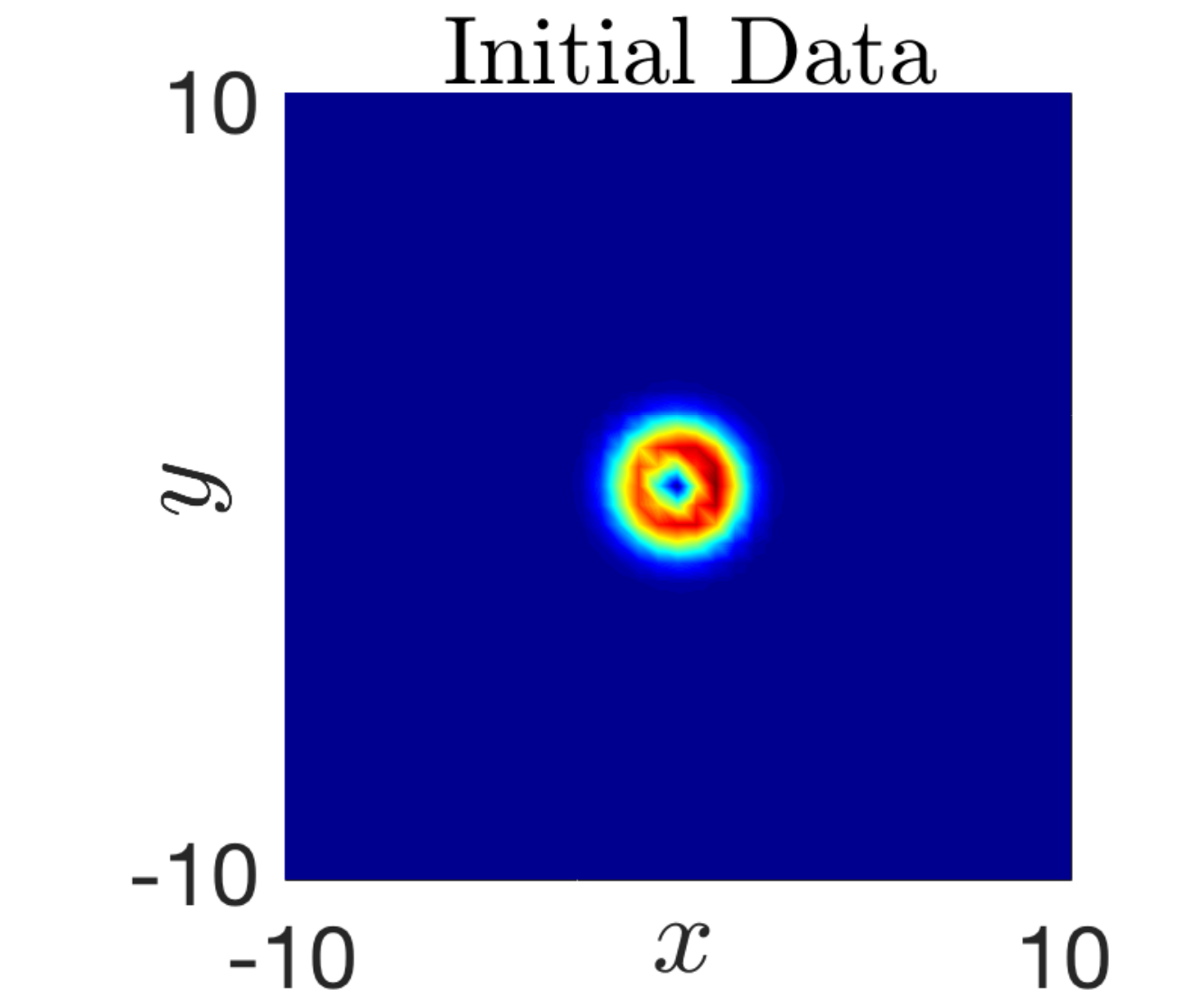}
\includegraphics[height=3.5cm,width=3.6cm]{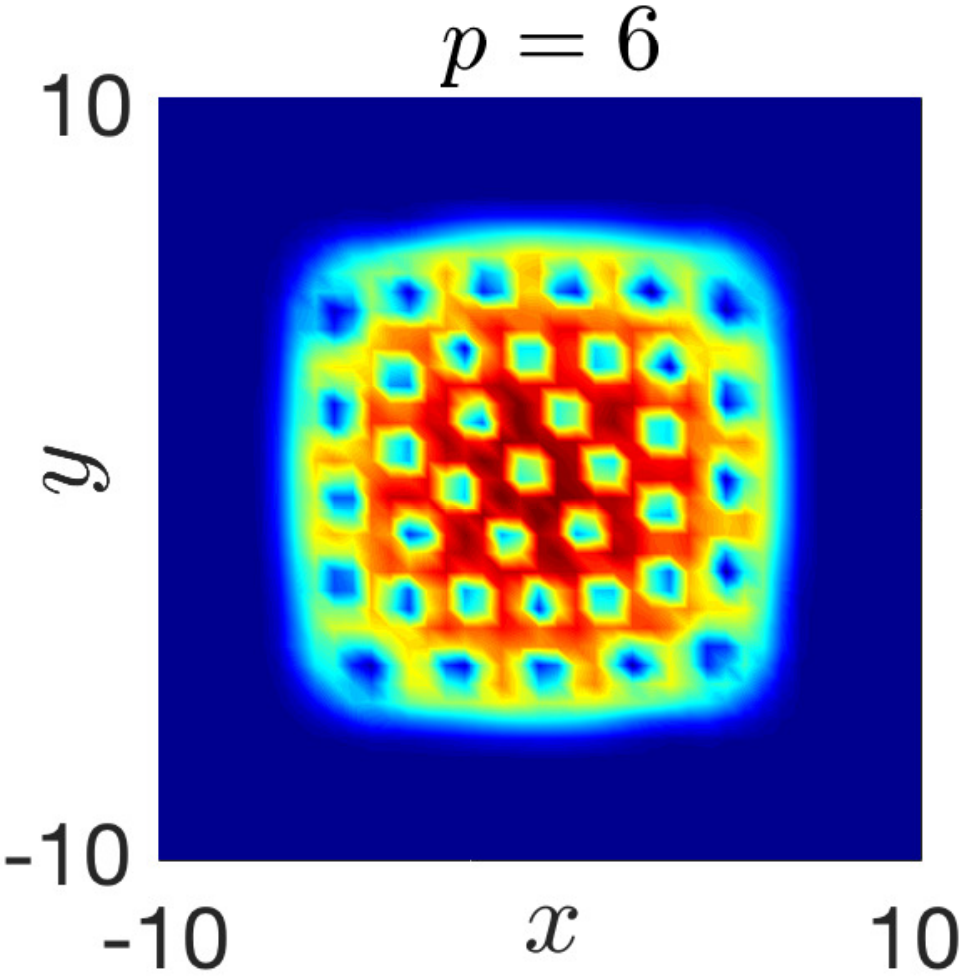}
\includegraphics[height=3.45cm,width=4.2cm]{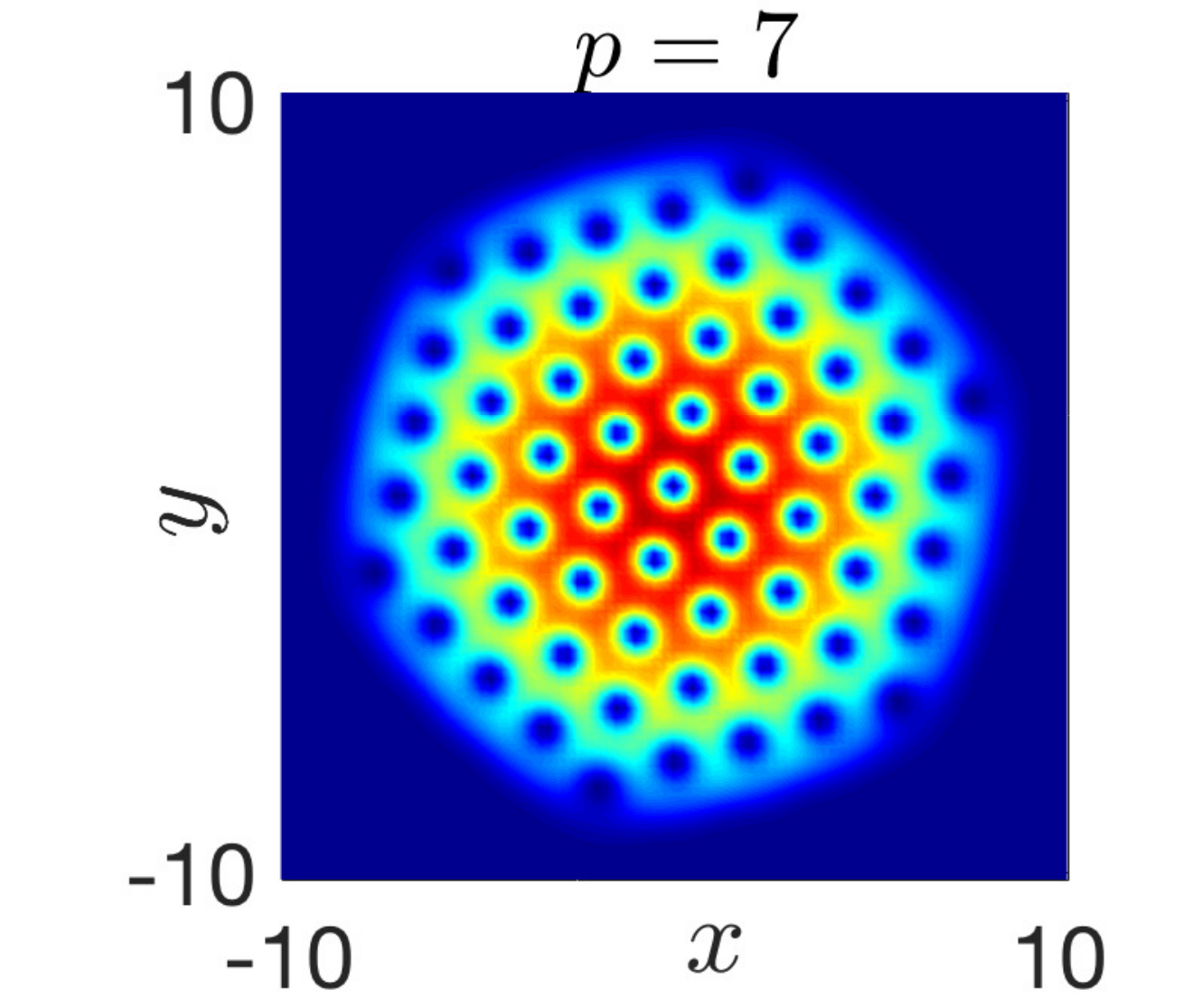}
\includegraphics[height=3.5cm,width=3.6cm]{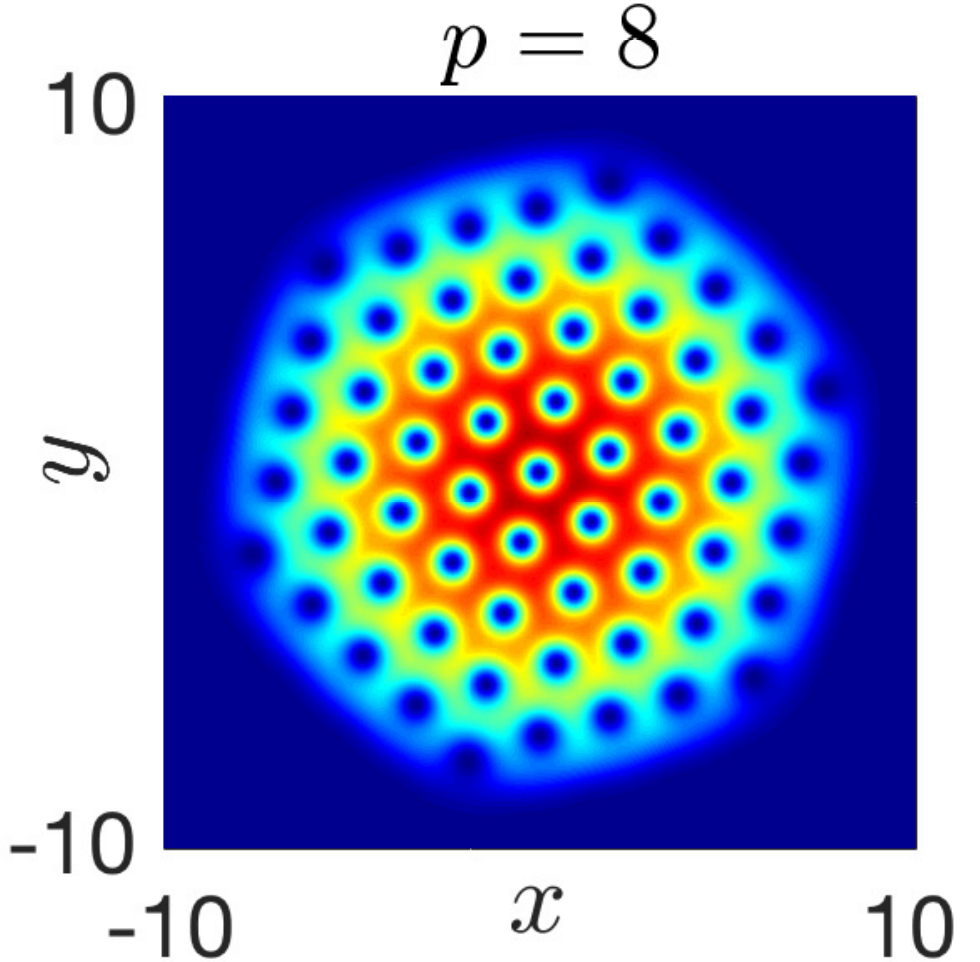}
}
\vspace{0.01cm}
\centerline{\quad\;
\includegraphics[height=0.8cm,width=3.2cm]{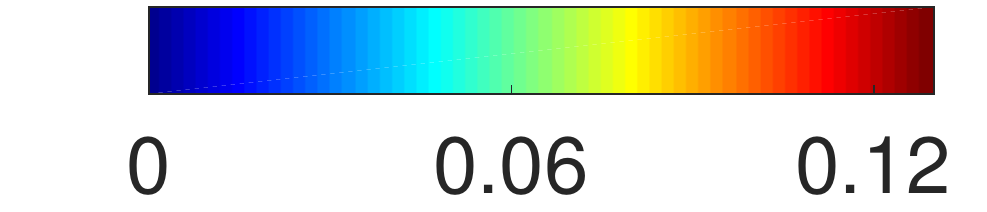}
\quad\quad\;\;
\includegraphics[height=0.8cm,width=11cm]{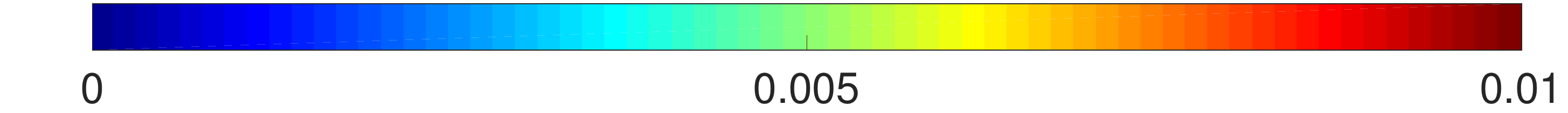}
}
\vspace{-0.2cm}
\caption{Exemple \ref{eg:diff_initial_fix_grid}. Energy error $\log_{10}(|E^{p}_{n}-E_{g}|)$  vs. the accumulated CPU time
 for $\og=0.95$ with initial data (d)   in Table \ref{tab:2d_diff_init_data_fix_grid} ($p=9$, upper left)
 and respectively  Table  \ref{tab:diff_init_data_multi_grid} (upper right) as well as  the 
stationary state obtained at each intermediate level (lower, $p=6, 7, 8$)}
\label{fig:2D_Multi_grid_Error_vs_CPU}
\end{figure}

\end{exmp}

\subsection{Numerical results in 3D}\label{Section3DNumerics}
 \begin{exmp}
 \label{eg:3D}
 {\rm 
Here, we apply the {\rm PCG}$_{\rm C}$ algorithm to compute some realistic 3D  challenging   problems.  To this end, 
$V(\bx)$ is chosen as the harmonic plus quartic potential \eqref{Quart_Poten}, with $\gm_x=\gm_y=1$, $\gm_z=3$,
 $\alpha=1.4$ and $\kappa=0.3.$    The computational domain is $\mathcal{D}=[-8,  8]^3$
   and the mesh size is:  $h=\fl{1}{8}$.
We test four cases:  (i) $\eta=100$, $\og=1.4$;  (ii) $\eta=100$, $\og=1.8$;
(iii) $\eta=5000$ and $\og=3$; (iv) $\eta=10000$ and $\og=3$. The initial guess
is always taken as the Thomas-Fermi initial data and the multigrid algorithm is used.
Fig. \ref{fig:3D_density_harmonic_quartic} shows the isosurface  $|\phi_{g}(\bx)|^2=10^{-3}$ and
the surface plot of $|\phi_{g}(x,y,z=0)|^2$ for the four cases. The CPU times  for these four cases are respectively
  2256 s, 1403 s, 11694 s and  21971 s.
  }

\begin{figure}[h!]
\centerline{
\includegraphics[height=3.2cm,width=3.6cm]{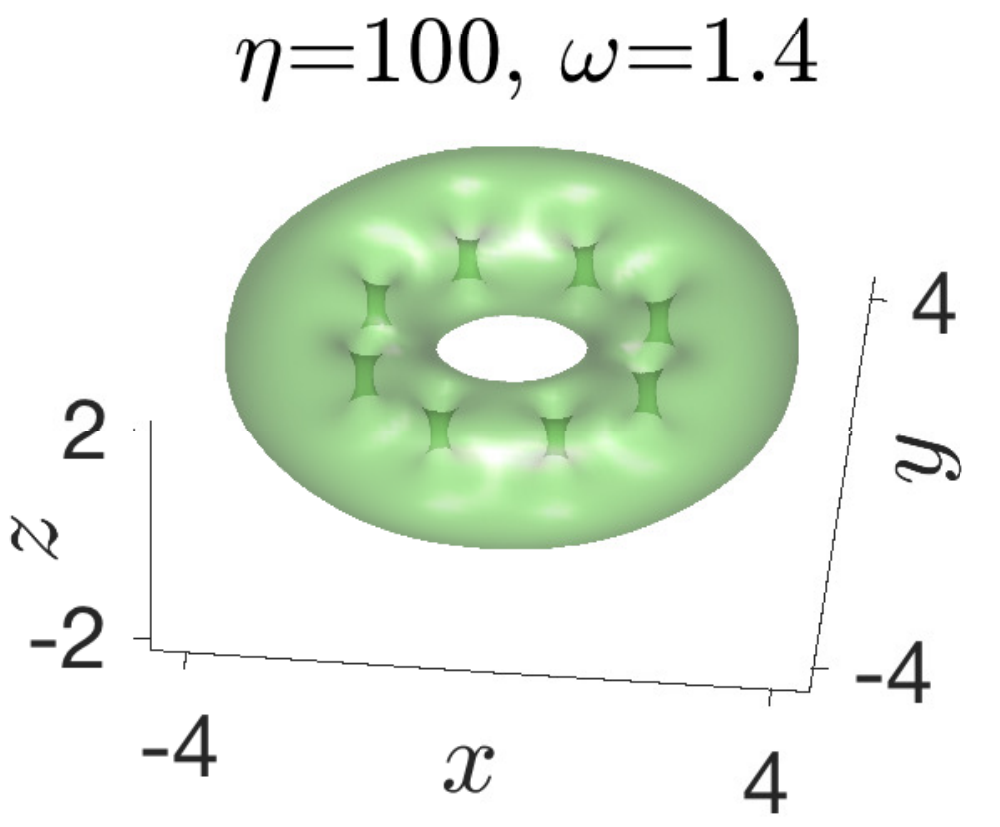}
\includegraphics[height=3.2cm,width=3.6cm]{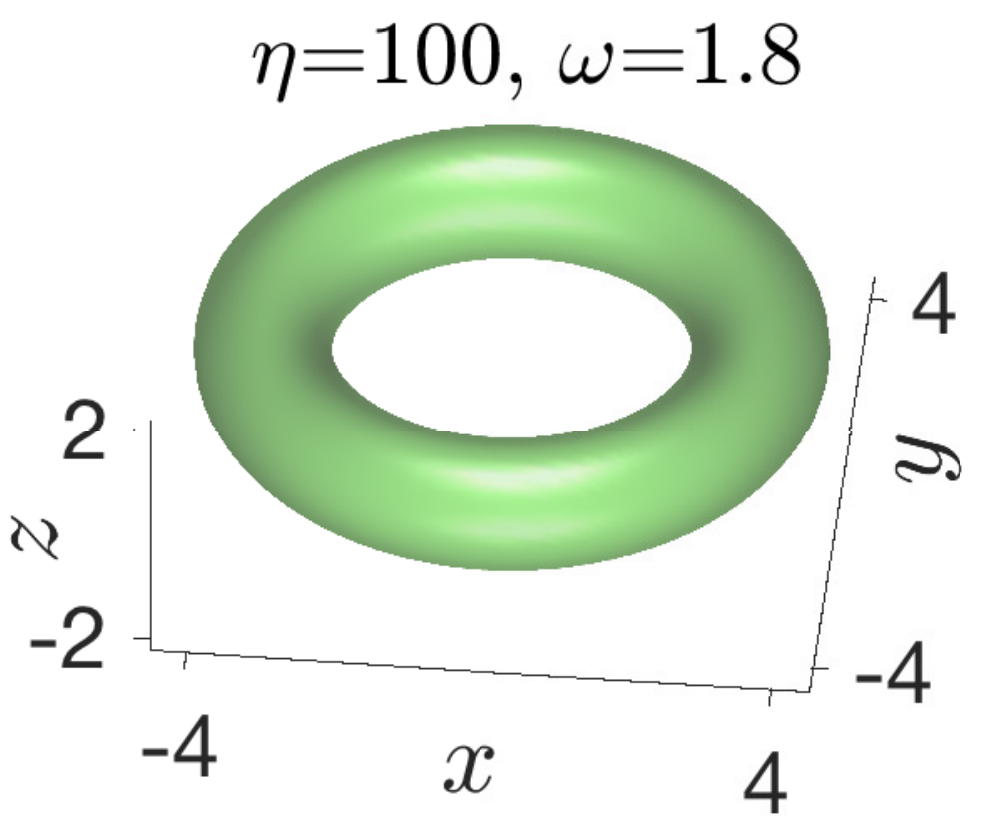}
\includegraphics[height=3.2cm,width=3.6cm]{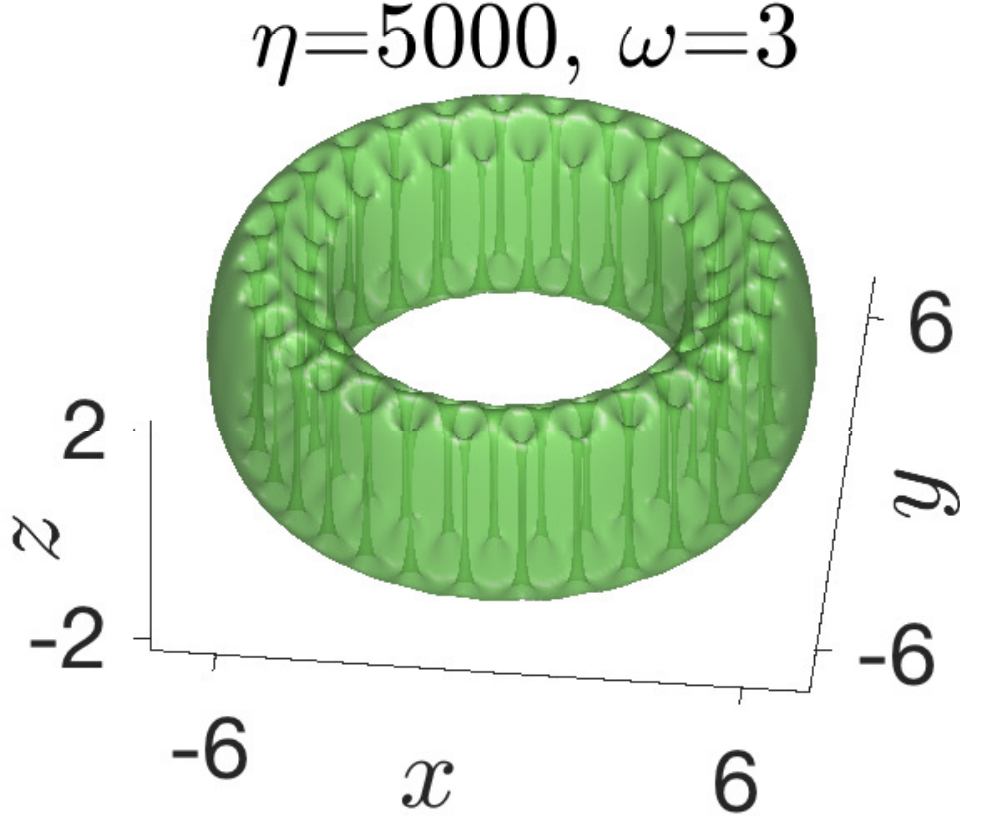}
\includegraphics[height=3.2cm,width=3.6cm]{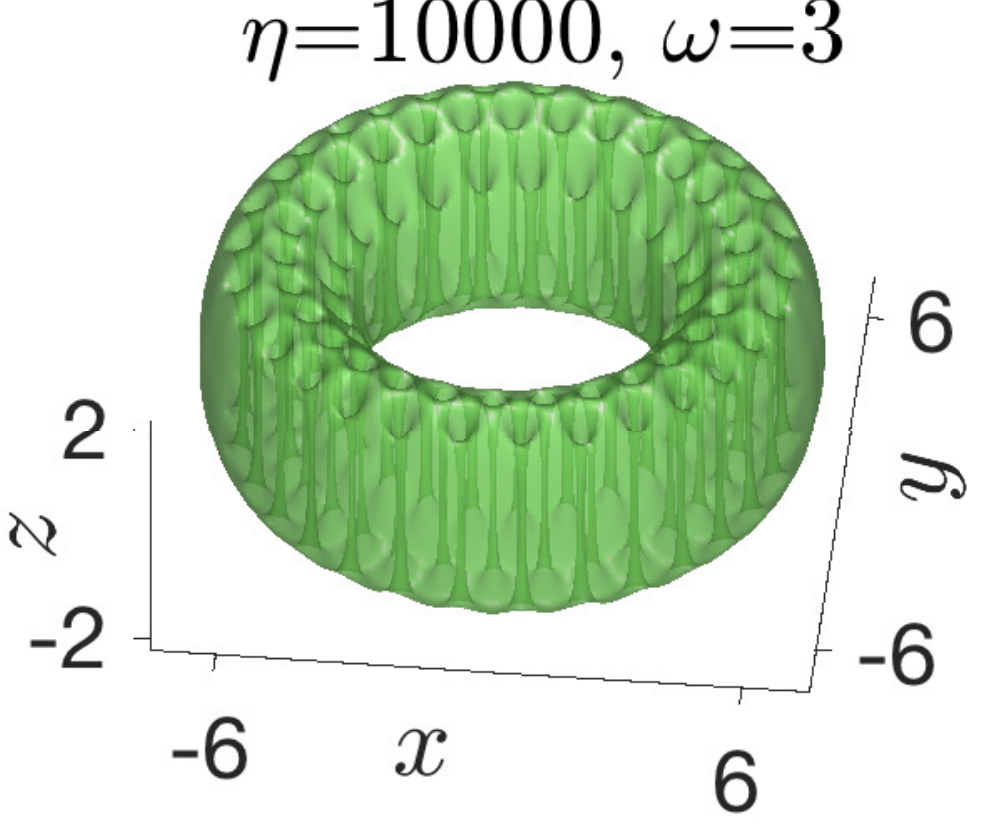}
}
\vspace{0.5cm}
\centerline{
\includegraphics[height=3.2cm,width=3.6cm]{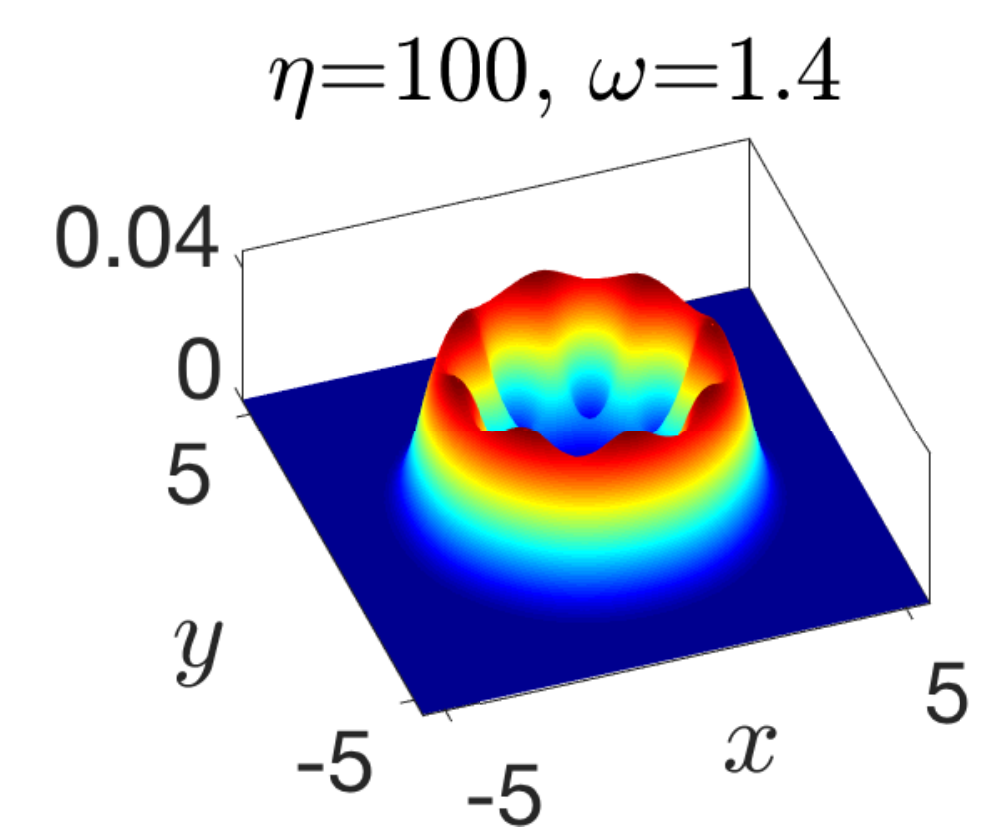}
 \includegraphics[height=3.2cm,width=3.6cm]{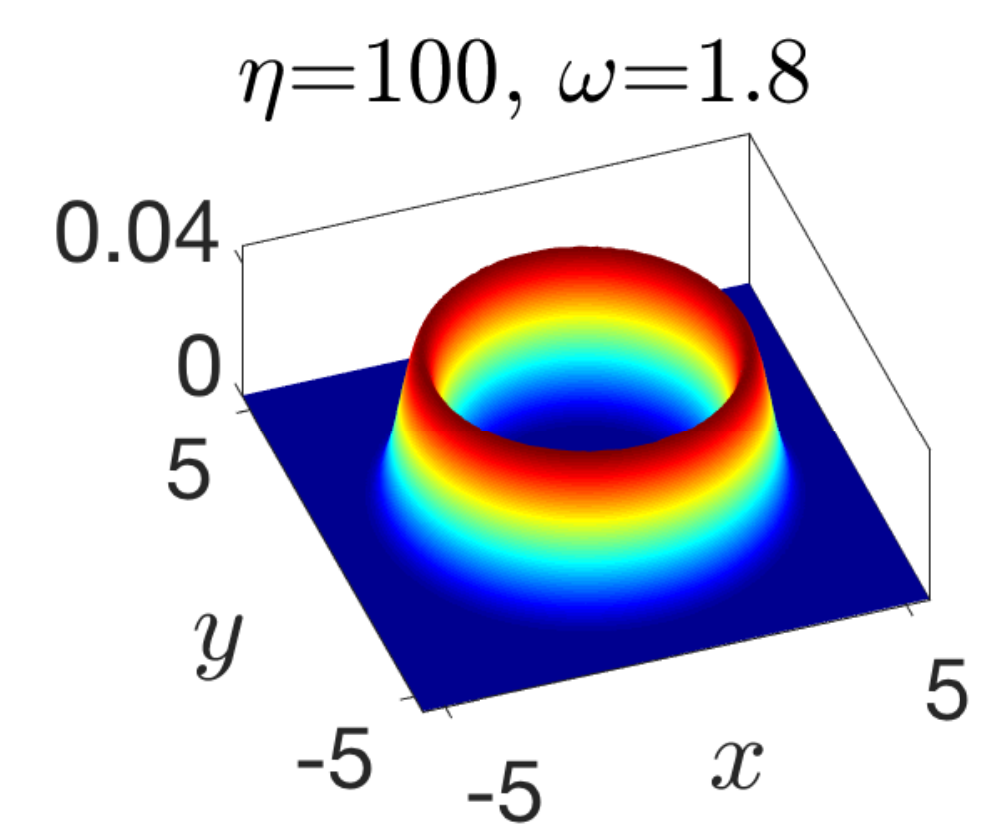}
 \includegraphics[height=3.2cm,width=3.6cm]{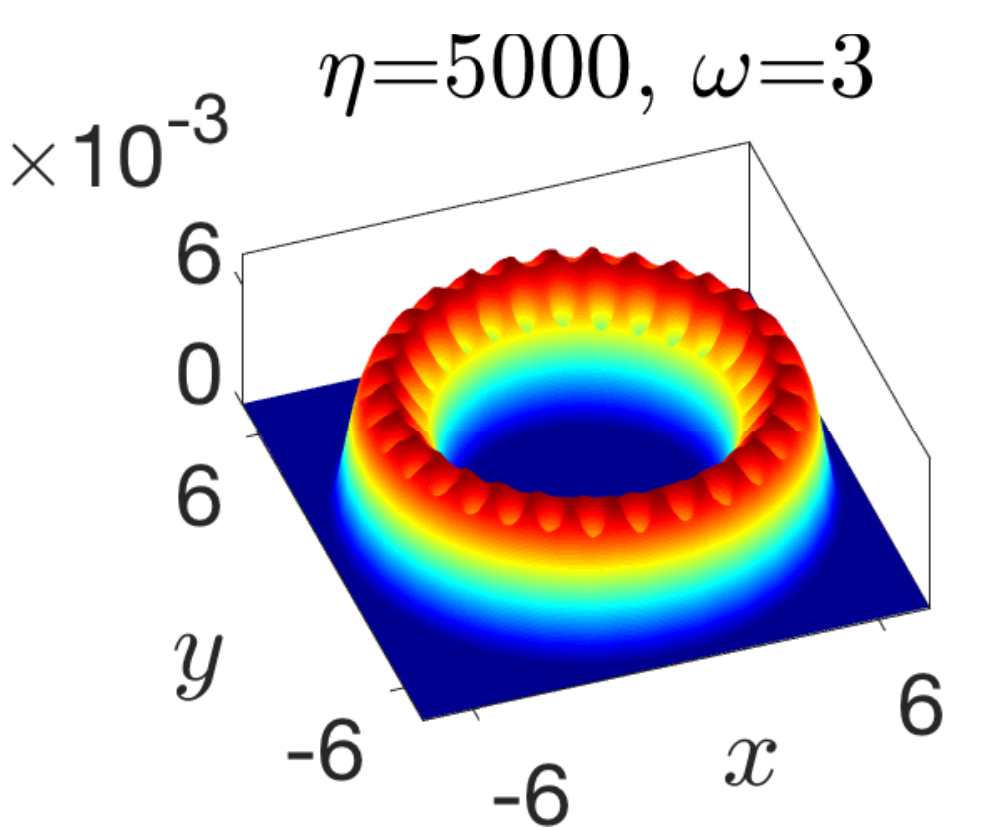}
  \includegraphics[height=3.2cm,width=3.6cm]{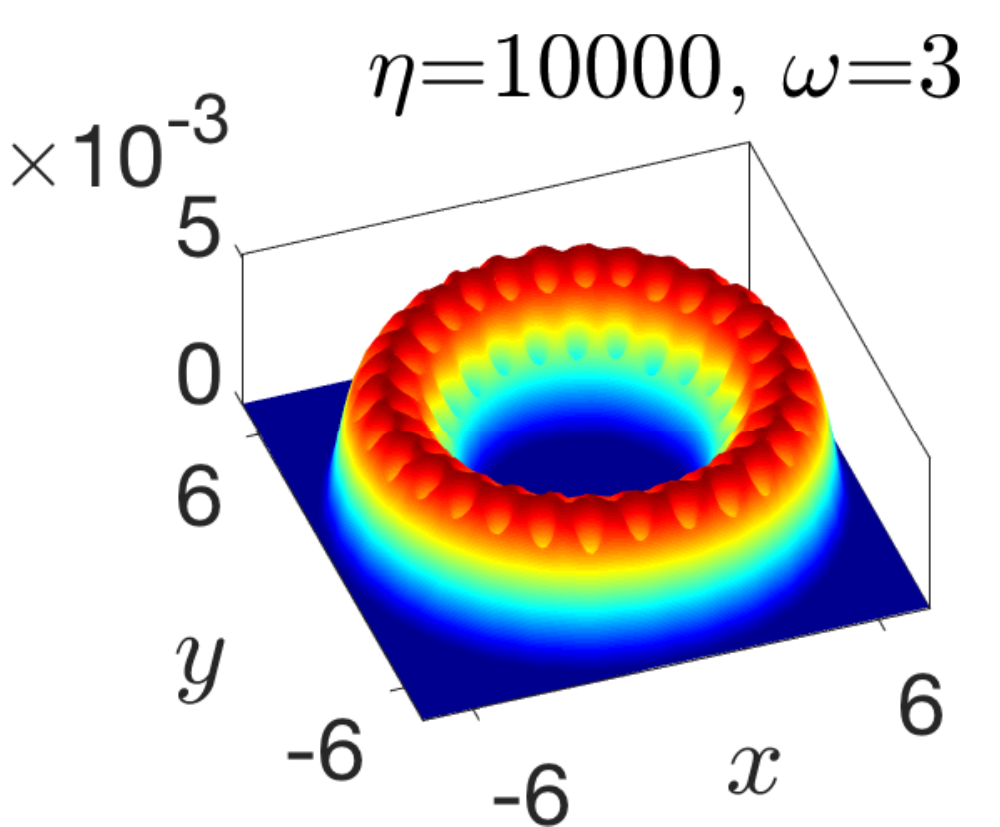}
}
\caption{Exemple \ref{eg:3D}. Isosurface $|\phi_{g}(\bx)|^2=10^{-3}$ (upper) and surface plot of $|\phi_{g}(x,y,z=0)|^2$ (lower) in example \ref{eg:3D}.
The CPU cost for these four cases are respectively
  2256 (\it s), 1403 ({\it s}), 11694 ({\it s}) and  21971  ({\it s}).}
\label{fig:3D_density_harmonic_quartic}
\end{figure}

 \end{exmp}

\section{Conclusion}\label{SectionConclusion}
We have introduced a new preconditioned nonlinear conjugate gradient
algorithm to compute the stationary states of the GPE with fast
rotation and large nonlinearities that arise in the modeling of
Bose-Einstein Condensates.  The method, which is simple to implement,
appears robust and accurate. In addition, it is far more efficient
than standard approaches as shown through numerical examples in
1D, 2D and 3D.  Furthermore, a simple multigrid approach can still
accelerates the performance of the method and leads to a gain of
robustness thanks to the initial data.  The extension to much more
general systems of GPEs is direct and offers an interesting tool for
solving highly nonlinear 3D GPEs, even for very large rotations.

\section*{Acknowledgements}
X. Antoine and Q. Tang thank the support of the French ANR grant  ANR-12-MONU-0007-02 BECASIM
(``Mod\`eles Num\'eriques'' call). 


\bibliographystyle{plain}
\bibliography{CGGPEbiblio}

\end{document}